\documentclass[11pt]{amsart}
\usepackage[latin1]{inputenc}
\usepackage{amssymb}

\usepackage[colorlinks, bookmarks=true]{hyperref}

\usepackage{nameref}
\usepackage{thmtools}

\usepackage[noabbrev,capitalize]{cleveref}
\usepackage{autonum}

\usepackage{amsfonts,amscd,amsmath,amssymb,enumerate}

\usepackage{multicol}

\usepackage[T1]{fontenc}
\usepackage{mathptmx}
\usepackage{stmaryrd}

\theoremstyle{plain} 
\newtheorem{theorem}{Theorem}[section]
\newtheorem{proposition}[theorem]{Proposition}
\newtheorem{lemma}[theorem]{Lemma}
\newtheorem{corollary}[theorem]{Corollary}
\theoremstyle{definition}
\newtheorem{remark}[theorem]{Remark}
\newtheorem{definition}[theorem]{Definition}

\numberwithin{equation}{section}

\newcommand{\unit}{\mathbf{1}}

\newcommand{\PP}{\mathcal{P}}

\begin{document}
	
\setlength{\abovedisplayskip}{9pt}
\setlength{\belowdisplayskip}{9pt}
\setlength{\abovedisplayshortskip}{9pt}
\setlength{\belowdisplayshortskip}{9pt}

	
\title[A Mackey-Gleason-Bunce-Wright theorem for JBW$^*$-algebras]{The Mackey-Gleason-Bunce-Wright problem for vector-valued measures on projections in a JBW$^*$-algebra}

\author[G.M. Escolano]{Gerardo M. Escolano}

\address{Departamento de An{\'a}lisis Matem{\'a}tico, Facultad de Ciencias, Universidad de Granada \hyphenation{Gra-nada}, 18071 Granada, Spain.}
\email{gemares@ugr.es}

\author[A.M. Peralta]{Antonio M. Peralta}

\address{Departamento de An{\'a}lisis Matem{\'a}tico, Facultad de Ciencias, Universidad de Granada, 18071 Granada, Spain.
Instituto de Matem{\'a}ticas de la Universidad de Granada (IMAG).}
\email{aperalta@ugr.es}

\author[A.R. Villena]{Armando R. Villena}

\address{ Departamento de An{\'a}lisis Matem{\'a}tico, Facultad de Ciencias, Universidad de Granada, 18071 Granada, Spain.
Instituto de Matem{\'a}ticas de la Universidad de Granada (IMAG).}
\email{avillena@ugr.es}



\begin{abstract} Let $\PP (\mathfrak{J})$ denote the lattice of projections of a JBW$^*$-algebra $\mathfrak{J}$, and let $X$ be a Banach space. A bounded finitely additive $X$-valued measure on $\PP(\mathfrak{J})$ is a mapping $\mu: \PP(\mathfrak{J}) \rightarrow X$ satisfying 
	\begin{enumerate}[$(a)$]
		\item $\mu(p +q) = \mu(p) + \mu(q)$, whenever $p \circ q = 0$ in $\PP (\mathfrak{J})$,
		\item $\sup \{ \| \mu(p)\| \, : \, p \in \PP(\mathfrak{J})\} < \infty$.
	\end{enumerate} In this paper we establish a Mackey-Gleason-Bunce-Wright theorem by showing that if $\mathfrak{J}$ contains no type $I_2$ direct summand, every bounded finitely additive measure $\mu: \PP(\mathfrak{J}) \rightarrow X$ admits an extension to a bounded linear operator from $\mathfrak{J}$ to $X$. This solves a long-standing open conjecture. 
\end{abstract}

\maketitle

\section{Introduction}

A \emph{countably additive positive measure} on the lattice of all closed subspaces of a complex Hilbert space $H$ is a mapping $\mu$ which assigns to each closed subspace $K_i\subseteq H$ a non-negative real number and enjoys the additional property of being countably additive, that is, $\displaystyle \mu\left(\overline{\hbox{span}}\left(\bigcup_{n} K_n\right)\right) =  \sum_{n} \mu\left(K_n\right)$, for every countable collection  $\{K_n\}_{n}$ of mutually orthogonal subspaces of $H$. It is much simpler if we replace closed subspaces with orthogonal projections in $B(H)$, and an example can be easily constructed by restricting any positive normal functional to the lattice of projections in $B(H)$. By addressing a question by G. Mackey on the mathematical foundations of quantum mechanics, A. Gleason established in \cite{Gleason1957} that for each separable Hilbert space $H$ of dimension at least three, every quantum measure on the lattice $\PP(B(H))$, of all projections in $B(H)$, extends to a positive normal functional on $B(H)$. A counterexample due to Kadison shows that the conclusion doesn't hold if dim$(H)=2$. Gleason's theorem was successfully applied in ruling out the possibility of certain classes of hidden variables in quantum theory and proving the essential difference between classical and quantum physics (cf. \cite[\S 7.3]{HamhalterBook2003}).\smallskip

Gleason posed the problem of determining all measures on the projections in a von Neumann factor, initiating in this way the study on the so-called Mackey-Gleason problem. The problem can be stated for more general measures on the wider setting of the lattice of projections in a JBW$^*$-algebra.  
\begin{definition}\label{def bounded finitely additive measures} Let $\PP(\mathfrak{J})$ be the lattice of all projections (i.e. symmetric idempotents) of a JBW$^*$-algebra $\mathfrak{J}$, and let $X$ be a Banach space. A bounded finitely additive $X$-valued measure on $\PP(\mathfrak{J})$ is a mapping $\mu: \PP(\mathfrak{J}) \rightarrow X$ satisfying the following axioms:
	\begin{enumerate}[$(a)$]
		\item $\mu(p +q) = \mu(p) + \mu(q)$, whenever $p \circ q = 0$ in $\PP (\mathfrak{J})$,
		\item $\sup \{ \| \mu(p)\| \, : \, p \in \PP(\mathfrak{J})\} < \infty$. 
	\end{enumerate} 
\end{definition} We can assume without loss of generality that $\sup \{ \| \mu(p)\| \, : \, p \in \PP(\mathfrak{J})\}=1$.\smallskip

A measure $\mu$ as above is called positive if $X = \mathbb{R}$ and $\mu (p)\geq 0$ for all $p\in \PP (\mathfrak{J}).$ In this case, boundedness is automatic. The so-called \emph{Mackey-Gleason problem} asks whether a bounded finitely additive $X$-valued measure on the lattice of projections of a von Neumann algebra, or more generally a JBW$^*$-algebra, with no type $I_2$ direct summand is the restriction of a bounded linear operator from $\mathfrak{J}$ to $X$. A complete positive answer to the Mackey-Gleason problem for positive measures on the lattice of projections of an arbitrary von Neumann algebra with no type $I_2$ direct summand was completed in a series of papers by E. Christensen in \cite{Christensen82} and F.W. Yeadon \cite{Yeadon1984}. L.J. Bunce and J.D.M. Wright showed in \cite{BunceWright1985} that the Mackey-Gleason problem admits a positive solution in the case of countably additive positive measures on the lattice of projections of each JBW$^*$-algebra with no type $I_2$ direct summand.\smallskip

An ingenious procedure, via Hahn--Banach theorem, allows to the reduce the Mackey-Gleason problem for vector measures to the case of $\mathbb{R}$-valued measures (see \cite{BunceWright1994}). A. Paskiewicz successfully considered complex measures on $\sigma$-finite von Neumann algebra factors of type different from $I_2$ (see \cite{Paskiewicz}). A definitive solution to the Mackey-Gleason problem for bounded vector-valued finitely additive measures on the lattice of projections of a von Neumann algebra with no type $I_2$ summands was obtained by L.J. Bunce and J.D.M. Wright in \cite{BunceWright1992,BunceWright1994}. We also acknowledge significant contributions by J.F. Aarnes \cite{Aarnes1969}, J. Gunson \cite{Gunson72}, S. Maeda \cite{Maeda1989}, and R. Cooke, M. Keane, W. Moran \cite{CookeKeaneMoran85}, and M.S. Matvechuk \cite{Matveichuk1988,Matveichuk1995}.  The references \cite{Matveichuk1988,Matveichuk1995} considered the Mackey-Gleason problem for AW$^*$- and JW$^*$-algebras, however, as expressed in \cite{BunceWright1994} the arguments contain some less clear aspects, or insufficiently justified tools. It was an admittedly open problem until now whether the Mackey-Gleason problem admits a positive solution for bounded vector-valued finitely additive measures on the lattice of projections of an arbitrary JBW$^*$-algebra without type $I_2$ summands (see for example \cite{Cufaro-Petroni}).\smallskip

The main goal of this paper is to provide a definitive argument to solve the Mackey-Gleason problem in the case of JBW$^*$-algebras without type $I_2$-summands. The result is stated in Theorems~\ref{theo: M_G_JBW_star} and \ref{t BWMG theorem for vector-valued measures}. We have developed a whole package of new tools in the Jordan setting to address this long-standing-open problem in full detail.\smallskip

This paper is structured as follows: this first section is complemented with several subsections containing the basic theory of JBW$^*$-algebras, and the references where the results can be consulted or extended. Section~\ref{sec: prelim} is devoted to study the quasi-linear extension determined by each bounded real-valued finitely additive measure $\mu$ on the lattice of projections of a JBW$^*$-algebra $\mathfrak{J}$. As we shall see below, the measure $\mu$ determines, uniquely, a mapping $\overline{\mu}: \mathfrak{J}\to \mathbb{C},$ called the quasi-linear extension of $\mu$, whose restriction to each self-adjoint associative JBW$^*$-subalgebra of $\mathfrak{J}$ is a self-adjoint bounded linear functional on this JBW$^*$-subalgebra. The quasi-linear extension of $\mu$ enjoys some other algebraic and geometric properties (see Proposition~\ref{prop_quasilin}), and our goal will consist in proving that this quasi-linear extension is the liner extension needed in the statement of the Mackey-Gleason theorem for JBW$^*$-algebras.\smallskip

In section~\ref{sec: modular type In} we prove that every bounded finitely additive real-measure on the lattice of projections of a JBW$^*$-algebra without type $I_2$ summand is uniformly continuous (see \Cref{prop: cont_mu}). We apply this result to give a complete positive solution to the Mackey-Gleason problem in the case of JBW$^*$-algebras of type $I_n$ with $2\neq n < \infty$ (cf. \Cref{theo: lin_type_In}). The case of JBW$^*$-algebras is studied in  \Cref{sec: properly non-modular}. It should be noted that in the different new technical tools employed in our arguments, the main difference with respect to the case of von Neumann algebras is that our arguments rely on the Bunce-Wright equivalence of projections in JBW$^*$-algebras, a notion strictly weaker that the usual Jordan equivalence of projections. The main conclusion in this section shows that  every bounded finitely additive real-measure on the lattice of projections of a properly non-modular JW$^*$-algebra $\mathfrak{J}$ extends to a linear functional on $\mathfrak{J}$ (see \Cref {theo: mu_lin_prop_inf}).\smallskip

Bounded finitely additive real-valued measures on the lattice of projections of a modular JBW$^*$-algebra are finally studied in \Cref{sec: modular JBW-algebras}. The novelties in the setting of modular JBW$^*$-algebras include, among other things, an \emph{Intermediate value property for centre-valued traces} JW$^*$-algebras of type $II_1$.  As in the case of von Neumann algebras, we prove that if $\tau: \mathfrak{J} \rightarrow Z(\mathfrak{J})$ denotes the normal centre-valued faithful unital trace on a JW$^*$-algebra of type $II_1$, for each $p \in \PP(\mathfrak{J})$ and each $w \in Z(\mathfrak{J})$ with $0\leq w \leq \tau(p),$ there exist $q \in \PP(\mathfrak{J})$ satisfying $q \leq p$ and $\tau(q) = w$ (see \Cref{prop: trace_surject}). In this highly technical section we elaborate suitable Jordan versions of results by Christensen \cite{Christensen82}, Maeda \cite{Maeda1989}, and Bunce and Wright \cite{BunceWright1994}, to present a detailed proof of the Mackey-Gleason theorem for bounded real-valued finitely additive measures on the lattice of projections of any JBW$^*$-algebra with no type $I_2$ summand (see \Cref{theo: mu_lin_type_II_1}).\smallskip

We devote the final paragraphs of this paper to show that for every type $I_2$ JBW$^*$-algebra $\mathfrak{J}$, there exists a positive finitely additive measure $\mu: \PP (\mathfrak{J})\to \mathbb{R}$ which does not admit an extension to a bounded linear functional on $\mathfrak{J}$.

\subsection{Basic notions and definitions}\label{subsec:background} \ \smallskip

This section is entirely devoted to recall some structure theory of JBW$^*$-algebras and their lattices of projections. The technical arguments in subsequent sections heavily depend on structure theory classifying all JBW$^*$-algebras in terms of the properties of their projection lattices, and the so-called ``dimension theory'' \cite[\S 5]{HOS} and \cite[\S 3]{AlfsenShultz2003}. We begin by recalling that a (real or complex) \emph{Jordan-Banach algebra} is a (real or complex) Banach space $\mathfrak{J}$ equipped with a bilinear mapping $(a,b) \rightarrow a \circ b$ (called the Jordan product) satisfying the following axioms: 
\begin{itemize}
    \item[(J1)] $\| a \circ b\| \leq \|a\|\|b\|$, for all $a,b \in \mathfrak{J}$; 
    \item[(J2)] $a\circ b = b\circ a$, for all $a, b \in \mathfrak{J} \ \ $  (commutativity);
    \item[(J3)] $(a^2 \circ b)\circ a = (a \circ b)\circ a^2$, for all $a,b \in \mathfrak{J} \ \ $ (Jordan identity).
\end{itemize}

A natural example is provided by any associative algebra $A$ equipped with the natural Jordan product given by $a\circ b := \frac{1}{2}(ab + ba)$. Any linear subspace of an associative algebra which is closed under the Jordan product is a Jordan algebra, such Jordan algebras are called \emph{special}.   Jordan algebras which can not be embedded as a Jordan subalgebras of an associative algebra are called \emph{exceptional}. A widely known example of an exceptional Jordan algebra is the algebra $H_3(\mathbb{O})$ of all Hermitian $3\times 3$ matrices with entries in the complex octonions (see \cite[Corollary 2.8.5]{HOS} for more informtaion).\smallskip


A Jordan algebra $\mathfrak{J}$ is called \emph{unital} if there exists an element $\mathbf{1} \in \mathfrak{J}$ (called the unit of $\mathfrak{J}$) such that $\mathbf{1}\circ a = a$ for all $a \in \mathfrak{J}$. Several properties of Jordan algebras, like invertibility, are characterized in terms of $U$-maps. Let us define these very special operators in the Jordan setting. Given elements $a,c \in \mathfrak{J}$, the symbol $U_{a,b}$ will stand for the linear mapping on $\mathfrak{J}$ defined by $$U_{a,c}(b) := (a\circ b)\circ c + (b\circ c)\circ a - (a\circ c)\circ b, \ \  (b\in \mathfrak{J}).$$ We simply write $U_a$ for $U_{a,a}$.  For each $a \in \mathfrak{J}$ we denote by $M_a$ the Jordan multiplication operator by the element $a$, that is, $M_a : \mathfrak{J} \rightarrow \mathfrak{J}, \, M_a (b) = a\circ b$.\smallskip

There are two closely related types of Jordan-Banach algebras, which are defined by algebraic--geometric axioms, and are known as JB-algebras and JB$^*$-algebras. A \emph{JB-algebra} is a real Jordan-Banach algebra $\mathfrak{J}$ in which the norm satisfies the following two additional conditions: 
\begin{itemize}
    \item[(JB1)] $\|a^2\| = \|a\|^2$  for all $a \in \mathfrak{J}$; 
    \item[(JB2)] $\|a^2\| \leq \|a^2 + b^2\|$ for all $a,b \in \mathfrak{J}$.
\end{itemize}
If $\mathfrak{J}$ is unital, it is clear from (JB1) that $\|\textbf{1}\| = 1$. The self-adjoint part of a C$^*$-algebra is a JB-algebra with respect to the natural Jordan product associated with every associative algebra.\smallskip

The Jordan analogue to C$^*$-algebras is constituted by JB$^*$-algebras, a model introduced by I. Kaplansky in 1976. A complex Jordan-Banach algebra $\mathfrak{J}$ equipped with an involution $^*$ is said to be a \emph{JB$^*$-algebra} if the following axiom is satisfied:
\begin{itemize}
\item[(JB$^*$1)] $\|a\|^3 = \| U_a(a^*)\|$ for all $a \in \mathfrak{J}$.
\end{itemize}

It is known that the involution of each JB$^*$-algebra $\mathfrak{A}$ is isometric, that is, $\|a^*\| = \|a\|$ for all $a \in \mathfrak{J}$ (cf. \cite[Lemma 4]{Youngson78}). If $\mathfrak{J}$ is unital, it can be easily seen that $\textbf{1}^* = \textbf{1}$.\smallskip

Along this paper, given a JB$^*$-algebra $\mathfrak{J}$, we shall make use of the triple product on $\mathfrak{J}$ defined by $\{a,b,c\} := U_{a,c} (b^*)$ ($a,b,c\in \mathfrak{J}$).\smallskip

JB- and JB$^*$-algebras are mutually linked in the following way: the set $\mathfrak{J}_{sa}$ of all \emph{self-adjoint} elements in a JB$^*$-algebra $\mathfrak{J}$, i.e. $\mathfrak{J}_{sa} : = \{a \in \mathfrak{J} \ : \, a^* = a\}$, is a JB-algebra (see \cite[Proposition 3.8.2]{HOS}). Conversely, by a deep result due to J.D.M. Wright \cite{Wright1977}, each JB-algebra corresponds to the self-adjoint part of a (unique) JB$^*$-algebra.\smallskip

A \emph{JBW$^*$-algebra} (resp., a \emph{JBW-algebra}) is a JB$^*$-algebra (resp., a JB-algebra) which is also a dual Banach space. Each JBW$^*$-algebra admits  a unique (isometric) predual (cf. \cite[Theorem 4.4.16]{HOS} or \cite[Theorem 2.55]{AlfsenShultz2003}). Thus JBW$^*$-algebras can be considered as the Jordan analogue of von Neumann algebras. It is known that a JB$^*$-algebra $\mathfrak{J}$ is a JBW$^*$-algebra if and only if $\mathfrak{J}_{sa}$ is a JBW-algebra \cite{Edwards80}. It is also established in the just quoted reference (see also \cite[Lemma 2.2]{BoHamKal2017}) that  the following assertions hold: 
    \begin{enumerate}[\rm(i)]\label{lemma:relJBW}
        \item[\rm(i)] $\mathfrak{J}_{sa}$ is weak$^*$-closed in $\mathfrak{J}$.
        \item[\rm(ii)] The operator $\phi: (\mathfrak{J}_*)_{sa} \rightarrow (\mathfrak{J}_{sa})_* $ defined by $\phi(\omega) = \omega |_{\mathfrak{J}_{sa}}$ is an onto linear isometry of real Banach spaces, where $$(\mathfrak{J}_*)_{sa} = \{\varphi\in \mathfrak{J}_{*} :\varphi (a^*) =\overline{\varphi(a)},\ \forall a\in \mathfrak{J} \}$$ is the self-adjoint part of the predual, $\mathfrak{J}_*$, of $\mathfrak{J}$.
        \item[\rm(iii)] The operator $\psi : \mathfrak{J}_{sa} \times \mathfrak{J}_{sa} \rightarrow \mathfrak{J}$ defined by $\psi(x,y) = x + iy$ is a onto real-linear weak$^*$-to-weak$^*$ homeomorphism.
    \end{enumerate}
   
JBW-algebras and their classification in terms of the properties of their lattices of projections are deeply studied in \cite{AlfsenShultz2003} and \cite{HOS}. The classification will be reviewed later.\smallskip

Given a JBW$^*$-algebra (respectively, a JB$^*$-algebra) $\mathfrak{J}$ and a subset $\mathcal{S}\subseteq \mathfrak{J}$ we shall denote by $W^*(\mathcal{S})$ (respectively, $J^*(\mathcal{S})$) the JBW$^*$-subalgebra (respectively, the JB$^*$-subalgebra) of $\mathfrak{J}$ generated by $\mathcal{S}$. To simplify the notation, we shall write $W^* (a_{1},\ldots, a_m)$ and $J^* (a_{1},\ldots, a_m)$
when $\mathcal{S}$ is the finite subset $\{a_{1},\ldots, a_m\}$.\smallskip

A \emph{JC-algebra} is a norm-closed Jordan subalgebra of the self-adjoint part of a C$^*$-algebra \cite[see Proposition 1.35]{AlfsenShultz2003}. There are examples of JB-algebras which are not JC-algebras, for instance the algebra, $H_3(\mathbb{O}),$ of all Hermitian $3\times 3$ matrices with entries in the complex octonions \cite[Corollary 2.8.5]{HOS}. A \emph{JC$^*$-algebra} is a JB$^*$-algebra which materialises as a norm-closed self-adjoint Jordan subalgebra of a C$^*$-algebra, equivalently, of some $B(H)$. Along this paper, the JBW$^*$-algebra obtained by complexifying $H_3(\mathbb{O})$ will be denoted by $H_3(\mathbb{O}^{\mathbb{C}}).$ \smallskip

A \emph{JW-algebra} is a weak$^*$-closed real Jordan subalgebra of the self-adjoint part of a von Neumann algebra. JW-algebras were first studied by D.M. Topping \cite{Topping65} and E. St{\o}rmer \cite{Stormer66}. A \emph{JW$^*$-algebra} is a JC$^*$-algebra which is also a dual Banach space, or equivalently, a weak$^*$-closed JB$^*$-subalgebra of some von Neumann algebra.\smallskip

Let $\mathfrak{A}$ and $\mathfrak{B}$ be pair of JB$^*$-algebras. A map $\varphi: \mathfrak{A} \rightarrow \mathfrak{B}$ is called a \emph{Jordan $^*$-homomorphism}, if it is $\mathbb{C}$-linear, preserves the involution (i.e. $\varphi (a^*) = \varphi(a)^*$) and the Jordan product $\varphi(a\circ b ) = \varphi(a)\circ \varphi(b)$ for all $a,b \in \mathfrak{A}$. A  \emph{Jordan $^*$-isomorphism} is a {Jordan $^*$-homomorphism} which is also bijective.\smallskip

An element $s$ in a unital JB$^*$-algebra $\mathfrak{J}$ is said to be a \emph{symmetry} if $s=s^*$ and $s^2= \mathbf{1}$. It is known that for each symmetry $s \in \mathfrak{J}$ the map $U_s$ is a Jordan $^*$-automorphism (i.e. a linear mapping preserving Jordan products and involution, see \cite[Proposition  2.34]{AlfsenShultz2003}). 

\subsection{Positive elements and the centre}\ \smallskip

According to the standard notation, elements $a$ and $b$ in a Jordan algebra $\mathfrak{J}$ are said to \emph{operator commute} if the operators $M_a, M_b$ commute (i.e. if $(a\circ c)\circ b = a \circ (c\circ b)$ for every element $c \in \mathfrak{J}$). The reader should be warned that operator commutativity of $a$ and $b$ is not, in general, related to the property that $a$ and $b$ generate a commutative and associative subalgebra  of $\mathfrak{J}$ (cf. \cite[2.5.1 and Example 2.5.2]{HOS}). The \emph{centre} of a Jordan algebra $\mathfrak{J}$ (denoted by  $Z(\mathfrak{J})$) consists of all elements $z\in \mathfrak{J}$ that operator commute with every element in $\mathfrak{J}$. Elements in the centre are called central. The centre of a JB$^*$-algebra $\mathfrak{J}$ is a commutative C$^*$-algebra, and contains the identity of $\mathfrak{J}$ if it exists (see \cite[Proposition 1.52]{AlfsenShultz2003}). The centre of a JBW$^*$-algebra is a commutative von Neumann algebra (see \cite[Proposition 2.36]{AlfsenShultz2003}, \cite{Edwards80}).\smallskip

Given a unital JB$^*$-algebra (respectively, a JB-algebra) $\mathfrak{J}$, the set $\mathfrak{J}^2 = \{ a^2 \, : \, a \in \mathfrak{J}\}$,  known as the \emph{cone of positive elements} in $\mathfrak{J},$ is a proper convex closed cone that induces a partial order on  $\mathfrak{J}$, making the latter a (norm) complete order-unit space whose distinguished order unit is the unit element of $\mathfrak{J}$, if any. The partial ordering on $\mathfrak{J}$ is concretely given by $a \leq b$ if $b-a \in \mathfrak{J}^2$ for $a, b \in \mathfrak{J}$ (see \cite[Theorem 2.1]{AlfShulStor79GelfandNeumark} or \cite[\S 1.2 and 3.3]{HOS}). Moreover, for each self-adjoint element $a$ in a JB$^*$-algebra $ \mathfrak{J}$, there exist unique positive elements $a_+$ and $a_-$ in $\mathfrak{J}$ such that $a = a_+ - a_-$, $a_+ \circ a_- = 0$, and the elements $a_+,$ and $ a_-$ belong to the JB$^*$-subalgebra generated by $a$, and $\| a\| = \max\{\|a_+\|, \|a_-\|\}$. Furthermore, each JB$^*$-subalgebra of $\mathfrak{J}$ generated by a single hermitian element is isometrically Jordan $^*$-isomorphic to a commutative C$^*$-algebra (see \cite[Theorem 3.2.2]{HOS}). 

\subsection{Structure theory for \texorpdfstring{JBW$^*$-}{JBW*-}algebras}\label{subsec: structure}\ \smallskip

An element $p$ of a JBW-algebra is called a \emph{projection} if $p^2= p \circ p = p$. Similarly, an element $p$ in a JBW$^*$-algebra $\mathfrak{J}$ is called a projection if $p^* = p = p\circ p$, i.e. $p$ is a projection in $\mathfrak{J}$ if and only if $p$ is a projection in $\mathfrak{J}_{sa}$. Two projections $p, q$ in $\mathfrak{J}$ are called \emph{orthogonal} ($p\perp q$ in short) if $p\circ q = 0$ (see  \cite[Lemma 4.2.2]{HOS} for equivalent reformulations). The symbol $\PP (\mathfrak{J})$ will stand for the orthocomplemented lattice of all projections in $\mathfrak{J}$ with the partial order given in the previous section. Note that $p\leq q$ in $\PP(\mathfrak{J})$ if, and only if, $p\circ q = p$.  The lattice $(\PP(\mathfrak{J}), \leq)$ is an orthomodular lattice for the order reversing map $p\mapsto p^{\perp}:=\unit-p$ \cite[Proposition 2.25]{AlfsenShultz2003}, that is, given $p,q \in \PP(\mathfrak{J})$, the supremum $p\vee q$ and infimum $p\wedge q$ of $p$ and $q$ exist in $\PP(\mathfrak{J})$ and the following properties hold:\begin{multicols}{2}
	\begin{enumerate}[$(i)$]
	\item $p^{\perp\perp} = p$.
	\item $p\leq q \Rightarrow p^{\perp}\geq q^{\perp}$.
	\item $p\vee p^{\perp} = \unit$ and $p\wedge p^{\perp}= 0.$
	\item If $p\leq q$, then $q = p \vee (q\wedge p^{\perp}).$
\end{enumerate}
\end{multicols} It is also known that $$(p\vee q)^{\perp} = p^{\perp}\wedge q^{\perp} \hbox{ and } (p\wedge q)^{\perp} = p^{\perp}\vee q^{\perp}, \hbox{ for all $p,q\in  \PP(\mathfrak{J})$.}$$ Let $p$ be a projection in $\mathfrak{J}$. The smallest central projection $c\in \PP(\mathfrak{J})$ satisfying $p \leq c$ is called the \emph{central cover of $p$}, and it is denoted by $c(p)$ (see \cite[Lemma 2.37 and Definition 2.38]{AlfsenShultz2003}). \smallskip

Two projections $p,q$ in a JBW-algebra $J$ (or in a JBW$^*$-algebra $\mathfrak{J}$) are called (\emph{Jordan}) \emph{equivalent} if there exist symmetries $s_1,\ldots, s_n$ in $J$ (with $n\in \mathbb{N}$) such that $q = U_{s_1}U_{s_2}\cdots U_{s_n} (p)$, and we write $p \sim q,$ or more concretely $p \sim_n q$.\label{def Jordan equivalence} If $n = 1$ we say $p$ and $q$ are \emph{exchanged by a symmetry}. We write $p \lesssim_n q$ if there exists a projection $q_0 \leq q$ such that $p\sim_n q_0$ \smallskip

A projection $p$ in a JBW$^*$-algebra $\mathfrak{J}$ is called \emph{abelian} (respectively, \emph{minimal}) if the algebra $\mathfrak{J}_p := U_p(\mathfrak{J})$ $= \{U_p(x) \, : \, x \in \mathfrak{J}\} $ is associative, i.e. $\mathfrak{J}_p = Z(\mathfrak{J}_p)$ (respectively, $\mathfrak{J}_p = \mathbb{C} p\neq \{0\}$). If $\mathfrak{J}$ is a JC$^*$-algebra regarded inside a C$^*$-algebra $A$, it follows from \cite[1.49]{AlfsenShultz2003} that $p$ is abelian if and only if $\mathfrak{J}_p$ consists of mutually commuting elements in $A$ in the usual sense.\smallskip

A projection $p$ in a JBW$^*$-algebra $\mathfrak{J}$ is called \emph{modular} or \emph{finite} if the projection lattice $[0,p] := \{ q \in \PP(\mathfrak{J}_p) \, : \, 0 \leq q \leq p\}= \{ q \in \PP(\mathfrak{J}) \, : \, 0 \leq q \leq p\}$ is modular, i.e. for every pair of  projections $q,e \in [0,p]$ such that $e \leq q $, the identity 
$$(e \vee q)\wedge r = e \vee (q \wedge r)$$ holds for every $r \in [0,p]$. If $\unit$ is modular, $\mathfrak{J}$ itself is called \emph{modular} (also called \emph{finite} in \cite{BunceWright1985}). We say $\mathfrak{J}$ is \emph{properly non-modular} if $\mathfrak{J}$ has no central modular projections except $0$. Finally, $\mathfrak{J}$ is called \emph{purely non-modular} if $\mathfrak{J}$ contains no modular projections except $0$ (see \cite[5.1.4]{HOS}, \cite[page 29]{Topping65}). A projection $p$ in $\mathfrak{J}$ is called \emph{properly non-modular} if $\mathfrak{J}_p = U_p (\mathfrak{J})$ is a properly non-modular JBW$^*$-algebra. We shall additionally $\mathfrak{J}$ is \emph{locally modular} if every direct summand of it contains a modular projection.\smallskip 

A JBW$^*$-algebra $\mathfrak{J}$ is of \emph{type $I$} if it contains an abelian  projection with central cover $\unit$. We say that $\mathfrak{J}$ is of \emph{type $II$} if and only if there is a modular projection $p$ in $\mathfrak{J}$ with 
$c(p) = \unit$, and $\mathfrak{J}$ contains no nonzero abelian projection. Finally, $\mathfrak{J}$ is of \emph{type $III$} if and only if it contains no nonzero modular projection (cf. \cite[Theorem 5.1.5]{HOS}).\smallskip 

For the purposes in this note we observe that each JBW$^*$-algebra $\mathfrak{J}$ admits a unique decomposition in the form $\mathfrak{J}= \mathfrak{J}_{sp}\oplus^{\infty}\mathfrak{J}_{ex}$, where $\mathfrak{J}_{sp}$ is a (special) JW$^*$-algebra and $\mathfrak{J}_{ex}$ is a purely exceptional JBW$^*$-algebra. Moreover, $\mathfrak{J}_{ex}$ is isometrically Jordan $^*$-isomorphic to $C(\Omega, H_3(\mathbb{O}^{\mathbb{C}}) )$, where $\Omega$ is a {hyperStonean compact Hausdorff} space (cf. \cite[Theorem 3.9]{Shultz79}, \cite[Theorem 2.8]{Wright1977} and \cite[Theorem 3.2]{Edwards80}). Furthermore, the JW$^*$-algebra $\mathfrak{J}_{sp}$ admits a finer decomposition in the form $\mathfrak{J}_{_{I_{mod}}}\oplus \mathfrak{J}_{_{I_{\infty}}}\oplus \mathfrak{J}_{_{II_{1}}}\oplus \mathfrak{J}_{_{II_{\infty}}}\oplus \mathfrak{J}_{_{III}}$, where the subindex refers to the corresponding type of the summand, that is, type ${I_{mod}}$ (type $I$ and modular), type ${I_{\infty}}$ (type $I$, properly non-modular, and locally modular), type ${II_{1}}$ (type $II$ and modular), type  ${II_{\infty}}$ (type $II$, properly non-modular, and locally modular) and type $III$ (purely non-modular).  $\mathfrak{J}$ is called \emph{properly infinite} or \emph{properly non-modular} if its finite or modular part, $\mathfrak{J}_{_{I_{mod}}} \oplus \mathfrak{J}_{_{II_{1}}},$ vanishes (see \cite{HOS,AlfsenShultz2003,Ayupov82,Topping65}).\smallskip

Type $I$ JBW$^*$-algebras admits a more detailed decomposition. A JBW$^*$-algebra $\mathfrak{J}$ is said to be of type $I_n$ if  there is a family $(p_j)_{j\in \Gamma}$ of abelian projections such that $c(p_j)= \unit,$ $\displaystyle \sum_{j\in \Gamma} p_j = \unit$, and card$(\Gamma) = n$. $\mathfrak{J}$ is said to be of type $I_{\infty}$ if it is a direct sum of JBW$^*$-algebras of type $I_n$ with $n$ infinite \cite[5.3.3]{HOS}. Each JBW$^*$-algebra of type $I$ admits a unique decomposition as a direct sum of JBW$^*$-algebras of type $I_n$ (see \cite[Theorem 5.3.5]{HOS}).   
\smallskip

Henceforth, we shall say that a non-zero  projection $p$ in a JBW$^*$-algebra $\mathfrak{J}$ can be halved if there exist projections $q_1,q_2\in \mathfrak{J}$ with $q_1 \perp q_2$, $q_1\sim_1 q_2,$ and $ p = q_1 + q_2$.

\begin{lemma}\label{remark: decomp_finite_alge} Let $\mathfrak{J}$ be a modular JBW$^*$-algebra. Then every non-zero projection $p \in \PP(\mathfrak{J})$ can be written in the form $p = q + r$, where $q$ is a projection that can be halved in case that it is non-zero, and $r$ is a possibly zero abelian projection.
\end{lemma}

\begin{proof} We may assume that $p$ is the unit of $\mathfrak{J}$. Since $\mathfrak{J}$ is modular, it writes as the orthogonal sum of JBW$^*$-algebras of type $I_n$ with $n < \infty,$ and a type $II_1$ JBW$^*$-algebra. The unit of the type $II_1$ summand can be halved by \cite[Lemma 5.2.14]{HOS}. By \cite[5.3.3]{HOS}, the unit of each JBW$^*$-algebra of type $I_n$ with $n\geq 1$ can be written as 
    $\displaystyle \unit =\sum_{i=1}^n p_i$, where $\{p_i\}$ is a family of orthogonal abelian projections such that $c(p_i) = \unit$. Therefore, if $n$ is even, i.e. $n = 2k$, $k \in \mathbb{N}$, the unit can be halved in two parts. If $n$ is odd, i.e. $n = 2k + 1$, $k \in \mathbb{N}$, we write
    $$ \unit = (p_1 + \dots + p_k) + (p_{k+1} + \dots p_{2k}) + p_{2k+1}, $$ which shows that it can be halved in two parts. The summand of type $I_1$ is abelian.\smallskip

    Finally, by glueing all parts together and noting that the sum of mutually orthogonal abelian projections is abelian, we obtain the desired decomposition. 
\end{proof}
    
Our next result is a Jordan version of a classical decomposition of infinite projections in von Neumann algebras (cf. \cite[Proposition 6.3.7]{KadRingrBookvol2}).

\begin{proposition}\label{remark: decomp_central-proj} Let $e$ be a non-modular projection in a JW$^*$-algebra $\mathfrak{J}$. Then there exists a central projection $z \in \PP(\mathfrak{J})$ such that $z \circ e$ is modular and $(\unit - z)\circ e$ is properly non-modular.
\end{proposition}

\begin{proof} Let us consider the JBW$^*$-algebra $\mathfrak{N}= U_e(\mathfrak{J})$, whose centre is precisely $Z(\mathfrak{J})\circ e$ (cf. \cite[Proposition 5.2.17]{HOS}). The set 
    $$ \mathcal{Q} = \Big\{ \{q_{j}\}_{j}\subseteq \PP(Z(\mathfrak{N}))\backslash\{0\} : \hbox{the $q_{j}$'s are mutually orthogonal and modular}\Big\}$$ is inductive with respect to the partial order given by inclusion. Hence, by Zorn's lemma, there exists a maximal family $\{q_{\lambda}\}_{\lambda}$ in $\mathcal{Q}$. We take now the projection $\displaystyle p = \sum_{\lambda}q_{\lambda}\in \mathfrak{N}$, where the series converges with respect to the weak$^*$-topology of $\mathfrak{N}$ (and with respect to the weak$^*$-topology of $\mathfrak{J}$). The projection $p$ is actually central in $\mathfrak{N}$ because $Z(\mathfrak{N})$ is a commutative von Neumann algebra. Since the projections in the family $\{q_{\lambda}\}_{\lambda}$ are central, modular and mutually orthogonal in $\mathfrak{N}$, we deduce from \cite[Corollary 20 or Lemma 22]{Topping65} that $p$ is modular in $\mathfrak{N}$.\smallskip

We claim that $e-p$ is properly non-modular in $\mathfrak{N}$. Indeed, let us take a projection $w \in Z(\mathfrak{N})$ such that $w\circ (e-p) \neq 0$. Note that $w\circ (e-p)\in Z(\mathfrak{N})$. If $z\circ (e-p)$ is modular, then the family $\{z\circ (e-p)\}\cup \{q_{\lambda}\}_{\lambda}$ lies in $\mathcal{Q}$ and contradicts the maximality of $\{q_{\lambda}\}_{\lambda}$. Therefore, $z\circ (e-p)$ is non-modular for every projection $w \in Z(\mathfrak{N})$ such that $w\circ (e-p) \neq 0$, which concludes the proof of the claim.\smallskip

Since $Z(\mathfrak{N}) = e\circ Z(\mathfrak{J})$, we can write $p = z\circ e$ where $z \in Z(\mathfrak{J})$. Having in mind that $e = p + (e-p)$, we deduce that $$z\circ e = z\circ p + z\circ (e-p)= z\circ (z\circ e) + z\circ (e-p) = z\circ e + z\circ (e-p),$$ which implies that $z\circ (e-p)=0$ (i.e. $z\perp (e-p)$). We consequently have $(e-p)= (\unit -z)\circ (e-p)$. Furthermore, $$\begin{aligned}
 (\unit-z) \circ e &= (\unit-z) \circ p + (\unit-z) \circ (e-p) \\
 &= (\unit-z) \circ (z\circ e) + (\unit-z) \circ (e-p) = e-p.
\end{aligned}$$ We have therefore proved that $$ e = z\circ e + (\unit - z)\circ e $$ where $z \in Z(\mathfrak{J})$, $z\circ e = p$ and $e-p= (\unit - z)\circ e$. \smallskip

Having in mind that since $p\leq e$, and hence $U_p(\mathfrak{N}) =  U_p (U_e(\mathfrak{J})) = U_p(\mathfrak{J})$, $p$ being modular in $\mathfrak{N}$ implies that $p$ is modular in $\mathfrak{J}$. Furthermore for each $w\in \PP(Z(\mathfrak{J}))$ with $w\circ (e-p) \neq 0$, by applying that $w$ is central we have $$ (w\circ e) \circ (e-p) = w \circ (e\circ (e-p)) = w \circ (e-p) \neq 0,$$ where $w\circ e \in e\circ Z(\mathfrak{J}) = Z(\mathfrak{N})$. It follows from the fact that $e-p$ is properly non-modular in $\mathfrak{N}$ that $w\circ (e-p)$ must be non-modular in $\mathfrak{J}$. Therefore $e-p = (\unit-z)\circ e$ is properly non-modular in $\mathfrak{J}$.\smallskip  

Consequently, there exists a central projection $z \in \mathfrak{J}$ such that $$ e =  z\circ e + (\unit - z)\circ e,$$  where $z\circ e$ is modular and $(\unit - z) \circ e$ is properly non-modular. 
\end{proof}

\subsection{Traces in \texorpdfstring{JBW$^*$-}{JBW*-}algebras}\label{subsec: traces}\ \smallskip

Following the pioneering works by D.M. Topping and S.A. Ayupov (see \cite[Section 18]{Topping65} and \cite{Ayupov82}), a \emph{trace} of a JB$^*$-algebra $\mathfrak{J}$ is a function  $\tau: \mathfrak{J}^+ \rightarrow [0, \infty]$ satisfying: 
\begin{enumerate}[$(i)$]
\item $\tau(x+y) = \tau(x) + \tau(y),$ for all $x, y \in \mathfrak{J}^+,$
\item $\tau(\lambda x) = \lambda \tau(x),$ for all $x \in \mathfrak{J}^+$, $\lambda \geq 0$;
\item $\tau(U_s(x)) = \tau(x),$ for all $x \in \mathfrak{J}^+$, and every arbitrary symmetry $s$ in $\mathfrak{J}$. 
\end{enumerate} A trace $\tau$ is said to be \emph{faithful} if $\tau(x) > 0$ for all non-zero $x \in \mathfrak{J}$. We say that $\tau$ is \emph{finite} if it is bounded on the set of all positive elements in the closed unit ball of $\mathfrak{J}$. If $\mathfrak{J}$ is a JBW$^*$-algebra, a trace $\tau$ on $\mathfrak{J}$ is called \emph{normal} if $\tau (\sup x_{j}) = \sup \tau(x_{j}),$ for every bounded increasing net $\{x_{j}\}$ with  $x_{j} \in \mathfrak{J}^+$ for all $j$. \smallskip
 
A \emph{centre-valued trace} on a JBW$^*$-algebra $\mathfrak{J}$ is a map $\tau : \mathfrak{J} \rightarrow Z(\mathfrak{J})$ satisfying: 
 \begin{enumerate}[$(i)$]
    \item $\tau(x+y) = \tau(x) + \tau(y),$ for all $x, y \in \mathfrak{J}$;
    \item $\tau(z\circ x) = z\circ \tau(x),$ for all $a \in \mathfrak{J}, z \in Z(\mathfrak{J})$;
    \item $\tau(x) \geq 0,$ if $x  \in \mathfrak{J}^+$;
    \item $\tau(U_s(x)) = \tau(x),$ for all $x \in \mathfrak{J}$, and every symmetry $s$ in $\mathfrak{J}$; 
    \item $\tau(\unit) = \unit$.  
\end{enumerate}  
\emph{Normal} and \emph{faithful} centre-valued traces are defined in similar terms to those given for traces in the previous paragraph.\smallskip

By combining \emph{the comparison theorem} for JBW$^*$-algebras \cite[Theorem 5.2.13]{HOS} with the presence of a  centre-valued trace we get the following Jordan version of \cite[Corollary V.2.8]{Takesaki_1}.

\begin{lemma}\label{lemma: centre_trace} Let $\mathfrak{J}$ be a JBW$^*$-algebra admitting a faithful centre-valued trace $\tau$. Let $p$ and $q$ be two projections in $\mathfrak{J}$. Then the following statements are equivalent:
\begin{enumerate}[$(a)$]
        \item $p \lesssim q$;
        \item $\tau(p) \leq \tau(q)$. 
    \end{enumerate}
\end{lemma}

\begin{proof} $(a)\Rightarrow (b)$ is clear from the properties of $\tau$. We shall prove $(b)\Rightarrow (a)$. Suppose that $\tau(p) \leq \tau(q)$. By the comparison theorem \cite[Theorem 5.2.13]{HOS} there exits a central projection $z \in \mathfrak{J}$ such that 
$$z\circ p \lesssim_1 z\circ q, \quad (\unit-z)\circ q \lesssim_1 (\unit-z)\circ p.$$ It follows from the properties of $\tau$ that 
    $$ (\unit - z)\circ \tau( q) = \tau((\unit - z)\circ q)  \leq \tau((\unit-z)\circ p) =(\unit-z)\circ \tau( p), $$ and since $\tau$ is centre-valued and $\tau(q)- \tau(p)\geq 0$ we get $$0\leq (\unit-z)\circ (\tau(q)- \tau(p) ) = \tau((\unit-z)\circ (q-p)) \leq 0,$$ which implies that $0=\tau((\unit-z)\circ (p-q))$ with $(\unit-z)\circ (p-q)\geq 0$, and hence $(\unit-z)\circ p = (\unit-z) \circ q$.\smallskip

Finally, since $z\circ p \lesssim_1 z\circ q$, we can find a symmetry $s\in \mathfrak{J}$ such that $U_s (z\circ p) \leq  z\circ q.$ Having in mind that $p = z\circ p + (\unit- z) \circ p,$ $ q = z\circ q + (\unit-z)\circ q$, and $s_1 = s\circ z + (\unit -z)$ is a symmetry in $\mathfrak{J}$ with $U_{s_1} (p) \leq q$, we get $p \lesssim q,$ which finishes the proof. 
\end{proof}

It is natural to ask when a JBW$^*$-algebra admits a faithful centre-valued trace. If $\mathfrak{J}$ is a modular JW$^*$-algebra, the existence is guaranteed by a result of D.M. Topping \cite{Topping65}. 

\begin{theorem}\label{theo: topping_trace}{\rm\cite[Theorem 26]{Topping65}}
 A JW$^*$-algebra is modular if, and only if, it possesses a (unique) faithful normal centre-valued trace.
\end{theorem}

The existence of faithful normal centre-valued traces on exceptional JBW$^*$-algebras has not been explicitly treated in the available literature. However, in the finite-dimensional case we find the next result by B. Iochum \cite{Iochum1984}.  

\begin{remark}\label{remark:Iochum}{\rm\cite[Remark V.1.3]{Iochum1984}}
Every finite dimensional JB$^*$-algebra $\mathfrak{J}$ admits a faithful finite  trace, which can be extended to a positive norm-one functional on $\mathfrak{J}$. 
\end{remark}

It follows from \cite[Theorem 3.9]{Shultz79} (see also \cite{Edwards80,Wright1977}) that every exceptional JBW$^*$-algebra $\mathfrak{J}$ can be identified as JBW$^*$-algebra with some $C(K, H_3(\mathbb{O}^{\mathbb{C}}))$, where $K$ is a hyperStonean compact Hausdorff space and $H_3(\mathbb{O}^{\mathbb{C}})$ is the exceptional JBW$^*$-algebra factor of all Hermitian $3\times 3$ matrices with entries in the complex octonions.

\begin{proposition}\label{prop: trace_exceptional} Every exceptional JBW$^*$-algebra admits a faithful normal centre-valued trace. Moreover, if $\mathfrak{F}$ is any finite dimensional JBW$^*$-algebra factor and $K$ is a compact Hausdorff space, the JB$^*$-algebra $C(K, \mathfrak{F})$ admits a faithful centre-valued trace.   
\end{proposition}

\begin{proof} Let us first prove the final statement. By \Cref{remark:Iochum}, there exists a positive norm-one functional $\tau: \mathfrak{F} \rightarrow \mathbb{C} $ whose restriction to $\mathfrak{F}^+$ is a faithful finite trace. Recall that $Z(\mathfrak{F}) = \mathbb{C}\unit_{_{\mathfrak{F}}}$, and  consider the map $\widehat{\tau} : C(K,\mathfrak{F})  \rightarrow C(K) \equiv Z(C(K,\mathfrak{F}))$ defined by $\widehat{\tau}(a) (t) := \tau (a(t))\cdot \unit,$ for every $a \in C(K,\mathfrak{F})$,  $t\in K$. It is not hard to see that $\widehat{\tau}$ is a faithful centre-valued trace on $C(K,\mathfrak{F})$. Namely, the linearity of $\widehat{\tau}$ is clear by definition. An element $a \in C(K,\mathfrak{F})$ is central if, and only if, for every $t \in K$, $a(t) \in \mathbb{C} \unit_{_{\mathfrak{F}}}$. Thus $\widehat{\tau}(a)(t) = \tau(a(t)) = a(t)$ for every $t \in K$, which implies that $\widehat{\tau}(a) = a$. This also proves that $\widehat{\tau}(\unit) = \unit$, and $\widehat{\tau}(z\circ a) = z\circ \widehat{\tau}(a)$ for all $a\in C(K,\mathfrak{F})$, $z\in Z(C(K,\mathfrak{F}))$.\smallskip
    
Note that positivity on $C(K,\mathfrak{F})$ is pointwise characterised, that is $a \in C(K,\mathfrak{F})$ is positive if, and only if, $a(t) \geq 0$ in $\mathfrak{F}$ for every $t \in K$. Therefore $ \widehat{\tau}(a)(t) = \tau(a(t)) \geq 0$ for every $t \in K$, and thus $\widehat{\tau}(a) \geq 0$. Similarly, symmetries in $C(K,\mathfrak{F})$ are pointwise determined, and hence $\widehat{\tau}(U_{{s}} (x)) = \widehat{\tau}(x)$, for all $s,x\in C(K,\mathfrak{F}),$ with $s$ being a symmetry. \smallskip
    
Finally, the faithfulness of $\widehat{\tau}$ is a straightforward consequence of the corresponding property for $\tau$. 
\end{proof}

\section{Quasi-linear maps on \texorpdfstring{JBW$^*$}{JBW*}-algebras}\label{sec: prelim}

\emph{Positive quasi-linear functionals} were originally introduced by J.F. Aarnes in \cite{Aarnes1969} to study the linearity of a physical state an arbitrary C$^*$-algebra. The quoted reference was a first approach to an analogue of Gleason's theorem (see \cite{Gleason1957}) for general von Neumann algebras, and what was later known as the Mackey-Gleason theorem. Our arguments also begin by considering quasi-linear maps associated with bounded finitely additive measures on the lattice of projections of a JBW$^*$-algebra. We introduce first some auxiliary tools. 

\begin{definition}\label{def: alpha and variation} Let $\mathfrak{J}$ be a JBW$^*$-algebra and let $\mu: \PP(\mathfrak{J}) \rightarrow \mathbb{R}$ be a bounded finitely additive measure. We consider the mappings  $\alpha_{_\mu}, V_{_\mu}: \PP(\mathfrak{J}) \rightarrow \mathbb{R}_0^+$ defined by $$\alpha_{_\mu}(p) = \sup \{ \mu(q) \, : \, q \in \PP(\mathfrak{J}), q \leq p\}\geq \mu (0) =0,$$ $$ \hbox{ and }  V_{_\mu}(p) = \sup \{ | \mu(q)| \, : \, q \in \PP(\mathfrak{J}), q \leq p \}\geq 0 ,$$ for all $p \in \PP(\mathfrak{J})$. 
\end{definition}

Note that, by the assumptions after \Cref{def bounded finitely additive measures}, $\alpha_{_\mu} (p), V_{_\mu} (p)\leq 1$ for all $p$ in $\PP(\mathfrak{J})$.

\begin{remark}\label{r properties of alpha and V} Let us observe several properties of the mappings $\alpha_{_\mu}$ and $V_{_\mu}$. We have already commented that both take positive values. 
	\begin{enumerate}[$(a)$]
		\item $\alpha_{_\mu} (p)\leq \alpha_{_\mu}(q)$ for all in $p,q\in \PP(\mathfrak{J})$ with $p\leq q$.
		\item If $z_1,\ldots,z_k\in \PP(Z(\mathfrak{J}))$ we have $\displaystyle \alpha_{_\mu} \left( \sum_{j=1}^k z_j\right) = \sum_{j=1}^k \alpha_{_\mu}(z_j).$ 
	\end{enumerate}
	Property $(a)$ is easy to check since every projection $e\leq p$ satisfies $e\leq q$. Property $(b)$ can be obtained via the arguments employed for von Neumann algebras in \cite[Lemma 2.1]{BunceWright1994}. 
\end{remark}

A functional $\phi$ in the dual space of a JB$^*$-algebra $\mathfrak{J}$ is called self-adjoint if $\phi (a^*) = \overline{\phi(a)}$ for all $a\in \mathfrak{J}$. We write $\mathfrak{J}^*_{sa}$ for the real space of all self-adjoint functionals in $\mathfrak{J}^*$. The assignment $\phi\mapsto \phi|_{\mathfrak{J}_{sa}}$ is a surjective real-linear isometry from  $\mathfrak{J}^*_{sa}$ onto $(\mathfrak{J})_{sa}^*$. Similar statements work for normal self-adjoint functionals, $(\mathfrak{J}_*)_{sa}$, in the predual of a JBW$^*$-algebra.\smallskip

\begin{remark}\label{r symmetries are norming} The set $\hbox{Symm}(\mathfrak{J}),$ of all symmetries in a JBW$^*$-algebra $\mathfrak{J},$ is a norming set for self-adjoint normal functionals in $\mathfrak{J}_*$. Namely, it is well known that $\hbox{Symm}(\mathfrak{J})$ is precisely the set of all extreme points of the closed unit ball of $\mathfrak{J}_{sa}$ (see, for example, \cite[Lemma 4.1]{Neal2000} or the more general result for real JB$^*$-triples in \cite[Lemma 3.3]{IsidroKaupPalacios}). So, the Krein-Milman theorem assures that the norm of every $\phi\in (\mathfrak{J}_{sa})_*$ can be computed as the supreme of the set $\{\phi (s) : s\in \hbox{Symm}(\mathfrak{J})\}$.\smallskip
	
Actually, $\hbox{Symm}(\mathfrak{J})$ is a norming set for self-adjoint functionals in $\mathfrak{J}^*$. To see the latter, let us fix $\phi \in \mathfrak{J}^*_{sa}$ with $\|\phi\| =1$. Fix an arbitrary $\varepsilon >0$ and find a norm-one element $x_{\varepsilon} \in \mathfrak{J}_{sa}$ such that $1-\varepsilon < \phi (x_{\varepsilon})$. Since $\mathfrak{J}_{sa}$ is a JBW-algebra, we can find a finite family of mutually orthogonal projections $\{p_1,\ldots, p_m\}$ in $\mathfrak{J}$, and real numbers $\lambda_1,\ldots, \lambda_m$ such that $\displaystyle \left\|x_{\varepsilon}- \sum_{j=1}^m \lambda_j p_j\right\| <\varepsilon$ (cf. \cite[Proposition 4.2.3]{HOS}). The JBW$^*$-algebra $W^*(\unit,p_1,\ldots,p_m)$ generated by $p_1,\ldots,p_m$ and $\unit$ is finite-dimensional. So, $\phi|_{_{W^*(\unit,p_1,\ldots,p_m)}}\in W^*(\unit,p_1,\ldots,p_m)_*$. By the conclusion in the previous paragraph, there exists a symmetry $s_{_{\varepsilon}}\in W^*(\unit,p_1,\ldots,p_m)$ such that $\| \phi|_{_{W^*(\unit,p_1,\ldots,p_m)}} \| -\varepsilon < \phi (s_{_{\varepsilon}})$. Clearly, $\| \phi|_{_{W^*(\unit,p_1,\ldots,p_m)}} \| > 1-2\varepsilon $, because $\displaystyle \sum_{j=1}^m \lambda_j p_j\in W^*(\unit,p_1,\ldots,p_m)$. By observing that $s_{_{\varepsilon}}$ also is a symmetry in $\mathfrak{J},$ with $\phi (s_{_{\varepsilon}}) > 1-3\varepsilon$, desired conclusion follows from the arbitrariness of $\varepsilon.$ 	
\end{remark}

We can now consider the natural connections with quasi-linear maps. We are inspired by classical arguments in \cite{Aarnes1969,BunceWright1985,BunceWright1992,BunceWright1994}.

\begin{proposition}{\label{prop_quasilin}}
Let $\mathfrak{J}$ be any JBW$^*$-algebra and let $\mu: \PP(\mathfrak{J}) \rightarrow \mathbb{R}$ be a bounded finitely additive measure. Then $\mu$ extends uniquely to a function $\overline{\mu}: \mathfrak{J} \rightarrow \mathbb{C}$ satisfying the following conditions: 
\begin{enumerate}[$(a)$]
    \item The restriction of $\overline{\mu}$ to each associative JBW$^*$-subalgebra $A$ of $\mathfrak{J}$ is a self-adjoint bounded linear functional on $A$. 
    \item $\overline{\mu}(x) \in \mathbb{R}$ and $\overline{\mu}(x+iy) = \overline{\mu}(x) + i \,\overline{\mu}(y)$ for all $x,y \in \mathfrak{J}_{sa}$.
    \item $ \sup\{|\overline{\mu}(x)| :  x = x^*, \|x\| \leq 1\}  = \sup \{\overline{\mu}(s) : s \in \mathfrak{J} \text{\small{ a symmetry}} \} 
         = 2 \alpha_{_{\mu}}(\unit) - \mu(\unit).$
    \item $2 \alpha_{_{\mu}}(\unit) - \mu(\unit)\leq 2 V_{_\mu} (\unit)$.
    \item $\sup\{\overline{\mu}(x): 0 \leq x \leq \unit \} = \alpha_{_{\mu}} (\unit).$
\end{enumerate}
\end{proposition}

\begin{proof} Let $A$ be any associative JBW$^*$-subalgebra of $\mathfrak{J}$. It is known that $A$ is an abelian von Neumann algebra, and obviously $A$ has no direct summand of type $I_2$. Therefore, by the Mackey-Gleason theorem for commutative von Neumann algebras (cf. \cite[Proposition 1.1]{BunceWright1992}), the restriction of $\mu$ to $\PP(A) (\subseteq \PP(\mathfrak{J}))$ extends uniquely to a bounded self-adjoint linear functional $\phi_A : A \rightarrow \mathbb{C}$. We consider now a function $\overline{\mu}: \mathfrak{J} \rightarrow \mathbb{C}$ defined by $$\overline{\mu}(x) = \phi_A (x)\in \mathbb{R},\hbox{ if } x \in \mathfrak{J}_{sa},$$ where $A$ is any associative JBW$^*$-subalgebra of $\mathfrak{J}$ containing the element $x$, and $\phi_{A}: A\to \mathbb{C}$ is the unique bounded self-adjoint functional extending $\mu|_{_{\PP(A)}},$ and $\overline{\mu} (x + i y) = \overline{\mu} (x ) + i \overline{\mu} (y)$, for all $x,y\in \mathfrak{J}_{sa}$.\smallskip 
    
We shall first show that $\overline{\mu}$ is well-defined. Namely, for each $x \in \mathfrak{J}_{sa}$ let $W^*(\unit,x)$ denote the associative JBW$^*$-subalgebra of $\mathfrak{J}$ generated by $x$ and the unit, $\unit\in \mathfrak{J}$, which is a commutative von Neumann algebra. Suppose now that the same element $ x$ belongs to an associative JBW$^*$-subalgebra $A$ of $\mathfrak{J}$. By definition, the restrictions of $\mu$ to $\PP(W^*(\unit,x))$ and to $\PP(A)$ give rise to two bounded linear functionals $\phi_{W^*(\unit,x)}\in W^*(\unit,x)^*$, $\phi_A\in A^*$,  respectively, satisfying  $\phi_{W^*(\unit,x)}|_{\PP(W^*(\unit,x))} = \mu|_{\PP(W^*(\unit,x))}$ and $\phi_A|_{\PP(A)} = \mu |_{\PP(A)}$. Note that $W^*(\unit,x) \subseteq A$, and thus $\phi_{W^*(\unit,x)}$ and $\phi_A$ coincide on $\PP(W^*(\unit,x))$, and consequently on the whole $W^*(\unit,x),$ since elements in the latter von Neumann algebra can be approximated in norm by finite linear combinations of mutually orthogonal projections in $W^*(\unit,x)$. Moreover, having in mind that $x\in W^*(\unit,x)$ we have $\phi_{A} (x) = \phi_{W^*(\unit,x)}$. The arbitrariness of $A$ guarantees that $\overline{\mu}$ is well defined. Statements $(a)$ and $(b)$ and the uniqueness of $\overline{\mu}$ clearly hold true by construction.\smallskip 

$(c)$ We keep the notation above. Fix $x= x^*\in \mathfrak{J}$ with $\|x\|\leq 1$. Let $A$ be any associative JBW$^*$-subalgebra of $\mathfrak{A}$ containing $x$ and $\unit$. It follows from $(a)$ that $\overline{\mu}|_{A}$ is a self-adjoint functional. The second paragraph in \Cref{r symmetries are norming} proves that for each positive $\epsilon$ there exists a symmetry $s_{\varepsilon}\in A$ (and hence in $\mathfrak{J}$) satisfying $ | \overline{\mu} (x)|= |\phi (x)| -\varepsilon \leq \|\phi\| -\varepsilon < \phi (s_{\varepsilon}) = \overline{\mu} (x)$. This proves the first equality in $(c)$.\smallskip
  
It is also known that every symmetry in $\mathfrak{J}$ writes in the form $s = 2 p -\unit$, with $p\in \PP(\mathfrak{J})$.  Therefore, having in mind that every symmetry in $\mathfrak{J}$ lies in an associative JBW$^*$-algebra containing it and the unit element, and the conclusion in $(a)$, we have:  
$$\sup\{\overline{\mu} (s) \, : \, s \in A \hbox{ a symmetry}\} = \sup\{\overline{\mu} (2 p -\unit) \, : \, p \in \PP(\mathfrak{J}) \}\hspace{3cm} $$ 
$$\hspace{4.4cm}  =\sup\{ 2\overline{\mu}(p) - \mu(\unit) \, : \, e \in \PP(\mathfrak{J})\} = 2 \alpha_{_{\mu}} (\unit) - \mu(\unit).
$$

$(d)$ As we have seen above, every symmetry $s\in \mathfrak{J}$ can be written in the form $s = p_1-p_2$, where $p_1,p_2\in \mathcal{P} (\mathfrak{J})$ with $p_1 \circ p_2 =0$ and $\unit = p_1+p_2$. By construction, $\overline{\mu}$ is linear when restricted to the JBW$^*$-subalgebra generated by $p_1$ and $p_2$, thus $$\overline{\mu} (s) = \mu (p_1)- \mu (p_2)\leq 2 \max\{|\mu (p_1)|, |\mu (p_2)| \} \leq 2 V_{\mu} (\unit),$$ the rest is a clear consequence of $(c)$.  \smallskip

$(e)$ Given any $0 \leq x \leq \unit,$ in $\mathfrak{J}$ with $\overline{\mu} (x) \neq 0$, and $\varepsilon>0$, we can find projections $p_1,\ldots,p_m$ inside the associative JBW$^*$-subalgebra, $W^*(\unit,x),$ of $\mathfrak{J}$ generated by $\unit$ and $x,$ and positive real numbers $\lambda_1,\ldots, \lambda_m\in (0,\|x\|]$ such that $\displaystyle \sum_{j=1}^m \lambda_j=\|x\|\leq 1,$ and $\displaystyle \left\| x - \sum_{j=1}^m \lambda_j p_j\right\| <\frac{\varepsilon}{\|\overline{\mu}|_{W^*(\unit,x)}\|} $. Therefore, by $(a)$, we get $$ \left| \overline{\mu} (x) - \sum_{j=1}^m \lambda_j \overline{\mu}(p_j)\right| <{\varepsilon}, \hbox{ and thus }  \overline{\mu} (x) - \varepsilon  < \sum_{j=1}^m \lambda_j \overline{\mu}(p_j),$$ which assures the existence of $j$ with $\overline{\mu}(p_j) > \overline{\mu} (x) - \varepsilon$. The desired statement follows from the arbitrariness of $\varepsilon>0$. 
\end{proof}
    
Let $\mathfrak{J}$ be a JBW$^*$-algebra. Henceforth, for each bounded finitely additive measure $\mu: \PP(\mathfrak{J}) \rightarrow \mathbb{R}$, the unique function $\overline{\mu}: \mathfrak{J}  \rightarrow \mathbb{C}$ whose existence is guaranteed by \Cref{prop_quasilin} will be called the \emph{quasi-linear functional} associated with $\mu$, or the \emph{quasi-linear extension of $\mu$}.\smallskip

Let us finish this section with an observation. The \emph{quasi-linear functional} $\overline{\mu}$ associated with a bounded finitely additive measure $\mu$ is linear if, and only if, it is additive on positive elements.\label{label additivity on positive elements is enough} Namely, the ``only if'' implication is clear. The ``if'' implication will follow if we show that $\overline{\mu}|_{\mathfrak{J}_{sa}} : \mathfrak{J}_{sa}\to \mathbb{R}$ is additive. We prove first that $\overline{\mu} (a-b) = \overline{\mu} (a) - \overline{\mu}(b)$ for every pair of positive elements $a,b$ in $\mathfrak{J}$. Namely, take a positive real number $\alpha$ such that $\alpha\unit -b \geq 0$.  Since $\overline{\mu}$ is additive on positive elements and linear on $W^*(\unit, a-b)$ and on $W^*(\unit, b)$ it follows that 
$$\alpha \overline{\mu} (\unit) + \overline{\mu} (-b + a) = \overline{\mu} (\alpha\unit -b + a)= \overline{\mu} (\alpha\unit -b ) + \overline{\mu} ( a) $$
$$\hspace{1.1cm} =\alpha \overline{\mu} (\unit)- \overline{\mu} (b ) + \overline{\mu} ( a),$$ which shows that $ \overline{\mu} (a-b ) =  \overline{\mu} ( a) -  \overline{\mu} (b )$, for all positive elements $a,b\in \mathfrak{J}$. Finally, given $h,k\in \mathfrak{J}_{sa}$, there exists positive elements $h^+, h^-, k^+,$ and $k^-$ such that $h = h^+ - h^-$ and  $k= k^+ - k^-$. Therefore $$\overline{\mu} (h+k) = \overline{\mu} (h^+ + k^+) - \overline{\mu} (h^- + k^-) \hspace{3cm}$$ $$\hspace{2cm} = \overline{\mu} (h^+ ) + \overline{\mu} (k^+) - \overline{\mu} (h^- ) - \overline{\mu} ( k^-) = \overline{\mu} (h) + \overline{\mu} (k),
$$ witnessing that $\overline{\mu}|_{_{\mathfrak{J}_{sa}}}$ is additive. 

\section{Uniform continuity and \texorpdfstring{JBW$^*$}{JBW*}-algebras of type \texorpdfstring{$I_n$}{I}}\label{sec: modular type In}

Let $\mathfrak{J}$ be a JBW$^*$-algebra of type ${I_{mod}}$ with no type $I_2$ summands. The goal of this section is to prove that every bounded finitely additive signed measure $\mu  : \PP(\mathfrak{J}) \rightarrow \mathbb{R}$ extends to a linear functional on $\mathfrak{J}$. The strategy will consist in proving that the {quasi-linear functional} associated with $\mu$ is linear on every JBW$^*$-subfactor of type $I_n$ ($n\geq 3$). We shall show that $\mu$ is uniformly continuous on $\PP(\mathfrak{J})$ in order to apply Christensen's method for locally finite functions (see \cite[Theorem 8.6]{Maeda1989}).\smallskip

Our first lemma has been borrowed from \cite{BunceWright1989}.

\begin{lemma}\label{l BunceWright lemma 2.1}{\rm(Bunce, Wright \cite[Lemma 2.1]{BunceWright1985})}
Let $\mathfrak{J}$ be a JBW$^*$-algebra containing distinct projections $f,g,e$ and a symmetry $s$ such that $f\circ g \neq 0,$  $e\circ f = e\circ g = 0$, $U_s(f) = e,$ $U_f(g) \in \mathbb{R}f,$ and $U_g(f) \in \mathbb{R}g.$ Then the self-adjoint part of the JBW$^*$-subalgebra of $\mathfrak{J}$ generated by $f,g,$ and $f\circ g$ is Jordan isomorphic to $M_3(\mathbb{R})_{sa}$, that is, $W^*(f,g, f\circ s)_{sa} \cong M_3(\mathbb{R})_{sa}$. Consequently, the JBW$^*$-subalgebra of $\mathfrak{J}$ generated by $f ,g,$ and $f \circ g$, $W^*(f,g,f\circ g)$, is Jordan $^*$-isomorphic to the JB$^*$-algebra $S_3(\mathbb{C})$ of all symmetric $3\times 3$ complex matrices. 
\end{lemma}

We can now deal with factor JBW$^*$-algebras of type $I_n$ with $n\geq 3$. In this case, the result proved by Bunce and Wright for positive measures on the lattice of projections of a JBW-algebra in \cite{BunceWright1985} plays a central role in the argument here.

\begin{theorem}\label{theo: mu_lin_fac} Let $\mathfrak{J}$ be a factor JBW$^*$-algebra of type $I_n$ with $3 \leq n < \infty$. Then every bounded finitely additive measure $\mu$ on $\PP(\mathfrak{J})$ extends to a bounded linear functional on $\mathfrak{J}$. 
\end{theorem}

\begin{proof} It is known from the structure theory that $\mathfrak{J}$ must be a finite dimensional JBW$^*$-algebra factor, actually $\mathfrak{J}_{sa}$ is Jordan isomorphic to $H_3(\mathbb{O})$ or to $M_n(\mathbb{F}),$ where $\mathbb{F} = \mathbb{R}, \mathbb{C}$ or $\mathbb{H}$ (cf. \cite[Theorem 5.3.8]{HOS}).\smallskip
	
It follows from \Cref{remark:Iochum} that $\mathfrak{J}$ admits a bounded linear faithful normal trace $\tau: \mathfrak{J}\to \mathbb{C}$. Observe that, by the finite-dimensionality of $\mathfrak{J}$ the set $\PP(\mathfrak{J})\backslash\{0\}$ is a norm compact subset of $\mathfrak{J}$, and $\tau$ never vanishes on $\PP(\mathfrak{J})\backslash\{0\}$. Therefore, the number $\kappa = \min\{\tau (p) : p\in \PP(\mathfrak{J})\backslash\{0\} \}$ is strictly positive. Consider now the following mapping $\nu : \PP(\mathfrak{J}) \to \mathbb{R}$, $\nu (p) := \frac{\alpha_{_\mu} (\unit)}{\kappa} \tau (p) -\mu (p)$. It is easy to check that $\nu$ is a positive, bounded, and finitely additive measure on $\PP(\mathfrak{J})$. We deduce from \cite[Theorem 2.2 or Theorem 3.8]{BunceWright1985} that $\nu$ extends to a positive linear functional on $\mathfrak{J}_{sa}$.  Since, by \Cref{prop_quasilin}$(b)$, $\overline{\mu}$ satisfies $\overline{\mu}(x + i y) = \overline{\mu}(x) + i \, \overline{\mu}(y)$ for all $x, y \in \mathfrak{J}_{sa}$, the mapping $\nu$ (and consequently, $\mu$) admits a complex linear extension to a bounded linear functional on the whole $\mathfrak{J}$.
\end{proof}

Henceforth, for each natural number $n\geq 2,$ the symbol $S_n(\mathbb{C})$ will stand for the JB$^*$-algebra of all complex $n\times n$ symmetric matrices.\smallskip

We turn now our attention to bounded finitely additive measures on modular type $I$ JBW$^*$-algebras without type $I_2$ summands.  Our next goal is to prove that any such a  measure is uniformly continuous on $\PP(\mathfrak{J})$. We shall need an additional technical result from \cite{BunceWright1989}.\smallskip

Let us first recall that, by the Shirshov-Cohn theorem (\cite[Theorem 2.4.14]{HOS} or \cite[Corollary 2.2]{Wright1977}), the JB$^*$-subalgebra, $J^*(\unit, a,b)$, of a unital JB$^*$-algebra $\mathfrak{J}$ generated by two self-adjoint elements $a$ and $b$ and $\mathbf{1}$ is isometrically JB$^*$-isomorphic to a JB$^*$-subalgebra subalgebra of a unital C$^*$-algebra $A$, and $\mathbf{1}$ is the unit of $A$. In case that $a$ and $b$ are two projections, say $p$ and $q,$ in $\mathfrak{J}$, it is not hard to check that \begin{equation}\label{eq norm Up(q)}\begin{aligned}
 \left\| U_p (p -q) \right\|_{\mathfrak{J}} &= \left\| U_p (\unit -q) \right\|_{\mathfrak{J}}	= \left\| U_p (\unit -q) \right\|_{A}	= \left\| p (\unit -q) p \right\|_{A} \\
 &= \left\| (\unit -q) p \right\|^2_{A}  = \left\| (p -q) p \right\|^2_{A} \leq \left\| p -q \right\|^2_{A} = \left\| p -q \right\|^2_{\mathfrak{J}}.
	\end{aligned}
\end{equation} \smallskip

Let $e$ and $f$ be projections in a JBW$^*$-algebra $\mathfrak{J}$. We say that $e,f$ are \emph{isoclinic} with angle $\theta \in [0, \frac{\pi}{2})$ if $U_e(f) = \cos^2 (\theta) e,$ and  $U_f(e) = \cos^2(\theta) f$.

\begin{lemma}\label{lemma: BunceWright_isoclinic}{\rm\cite[Lemma 3.6]{BunceWright1985}}
 Let $f$,$g$ be projections in a JBW$^*$-algebra $\mathfrak{J}$. Suppose that $\|f-g\|< 1,$ and that there is a projection $e$ in $\mathfrak{J}$ such that $e \sim f,$ $e\circ f = e\circ g = 0$. Suppose also that there exists an angle $\theta$ in the interval $[0, \pi/4)$  such that $sin^{-1}(U_f(\unit - g)^{1/2}) \leq 2\theta f$. Then there exists a projection $h$ in $\mathfrak{J}$ isoclinic with angle $\theta$ to both $f$ and $g$. In particular, this conclusion  holds when we take $\theta = \frac{1}{2} \sin^{-1}(\| U_f(\unit -g)^{1/2}\|) = \frac{1}{2} \sin^{-1}(\| U_f(f -g)^{1/2}\|),$ and in such a case we also have $\| f-h\| , \| g-h\| \leq \| f-g\|$.
\end{lemma}

\begin{proof} As commented above, the first part of the result is explicitly proved in \cite[Lemma 3.6]{BunceWright1985}. To see the last  statement, let us observe that having in mind that $f$ and $h$ are isoclinic with angle $\theta$, the JBW$^*$-subalgebra of $\mathfrak{J}$ generated by $f$ and $h$ can be identified with the algebra $S_2(\mathbb{C})$ fo all complex $2\times 2$ symmetric matrices in such a way that $f$ and $h$ are identified with the matrices $\left(\begin{matrix}
		1 & 0 \\
		0 & 0
	\end{matrix}\right)$ and $\left(\begin{matrix}
	\cos^2(\theta) & \cos(\theta) \sin(\theta) \\
	\cos(\theta) \sin(\theta) &  \sin^2(\theta)
	\end{matrix}\right)$, respectively (see, for example, \cite[Proposition 3.5]{PeSaa2025}). We therefore have 
$$\begin{aligned}
\|f-h\| &= \left\| \left(\begin{matrix}
	1- \cos^2(\theta) & -\cos(\theta) \sin(\theta) \\
	-\cos(\theta) \sin(\theta) &  - \sin^2(\theta)
\end{matrix}\right) \right\| \\
&= \sin(\theta ) \left\| \left(\begin{matrix}
	\sin(\theta) & -\cos(\theta)  \\
	-\cos(\theta) &  - \sin(\theta)
\end{matrix}\right) \right\| \\
&= \sin(\theta ) = \sin\left(\frac{1}{2} \sin^{-1}(\| U_f(\unit -g)^{1/2}\|) \right)\\
 &\leq  \| U_f(\unit -g)\|^{1/2}=\hbox{(by \eqref{eq norm Up(q)})} \leq \|f -g\|^{1/2}. 
\end{aligned}$$  The conclusion concerning the norm of $g-h$ follows by similar arguments.  
\end{proof}

We complete our list of technical tools with the next lemma. 

\begin{lemma}\label{lemma: symproj}
   Let $f, g$ be projections in the JBW-algebra $\mathfrak{J}$  such that $\|f - g \| < 1$. Then there is a symmetry $s$ in $\mathfrak{J}$ such that $U_s(f) = g$ and $\| p - U_s(p)\| \leq \sqrt{2} \|f-g\|^{\frac12}$ for any projection $p \leq f$.  
\end{lemma}

\begin{proof} By \cite[Remark 3.4(2)]{BunceWright1985} the symmetry $s = c^{-1/2}\circ (f + g -\unit)$ exchanges $f $ and $g$ (i.e. $U_s (f) =g$), where the element $c =\unit- (f - g )^2 = (f + g -\unit)^2 = \unit -f -g + 2 f\circ g \geq 0$ is invertible in $\mathfrak{J}$. According to the fixed notation, let $W^*(\unit, f,g)$ be the JBW$^*$-subalgebra of $\mathfrak{J}$ generated by $f$ and $g$. It is known that $W^*(\unit, f,g)$ is a JW$^*$-algebra (cf. \cite[Proposition 3.3]{PeSaa2025}). Working in any von Neumann algebra $W$ containing $W^*(\unit, f,g)$ as a JBW$^*$-subalgebra and sharing the same unit, it can be easily checked that $c$ and $f$ commute in $W$, and hence operator commute in $W^*(\unit, f,g)$, and in $\mathfrak{J}$  (cf. \cite[Proposition 1.2]{EscoPeVi2025} or \cite[Proposition 1]{Topping65}). Similarly, $c$ and $g$ operator commute in $\mathfrak{J}$. Furthermore, by the same arguments above we have  \begin{equation}\label{eq Uf(s)}
		U_f(s) = c^{-1/2}\circ U_f(g) \hbox{ in $W$ and in $\mathfrak{J}$.}
	\end{equation} The fundamental identity of Jordan algebras (see \cite[Proposition 3.4.15]{CabGarPalVol1}) assures that \begin{equation}\label{eq Uf(s)2} U_f(s)^2  = U_{U_f(s)} (\unit) = U_f U_s U_f (\unit) = U_f U_s (f) = U_f (g).
	\end{equation}
	     
Now, let $p$ be a projection in $\mathfrak{J}$ with $p \leq f$ and set $q = U_s(p)$. Clearly, $U_{s} (q) = p$ by just considering a von Neumann algebra containing the JBW$^*$-subalgebra generated by $\unit, s$ and $p$. Observe that, by orthogonality, we have
     $$\| p - q \|^2 = \| U_p(\unit -q) + U_{\unit-p}(q)\| = \max \{\| U_p(\unit-q)\|, \| U_{\unit-p}(q) \|\},$$
     where $\|U_{\unit-p}(q)\| = \| U_q(\unit-p)\|$ (cf. \cite[Lemma 3.5.2]{HOS}). Thus $\|p -q\|^2= \max\{\|U_p(\unit-q)\|,\|U_q(\unit-p)\|\}$. 
     However, since $U_s$ is an isometry \cite[Theorem 4.2.28$(vii)$]{CabGarPalVol1}, the fundamental identity of Jordan algebras also assures that $$\|U_q(\unit-p)\| = \| U_s U_q(\unit-p)\| = \| U_s U_q U_s (\unit-q)\| = \| U_{U_s (q)}  (\unit-q)\| = \| U_p(\unit-q)\|.$$ We finally compute the norm of $p-q$;
     $$ 
     \begin{aligned}
       \| p -q \|^2 & = \|U_p(\unit-q) \|= \|p-U_p(q) \| = \left\|p-U_{p} U_s (p) \right\| = \|p-U_{p} U_s U_f (p) \|
       \\ &= \| p - U_pU_fU_sU_f(p) \|= \| p - U_p U_{U_f(s)} (p) \| = \| p - U_p U_{U_f(s)} U_p (\unit) \| 
     \end{aligned}$$\vspace*{-4mm}
     $$
     \begin{aligned}
     	\hspace*{4mm} &= \| p - U_{U_p (U_f(s))} (\unit) \| = \left\| p - \Big( U_p \left(U_f(s)\right) \Big)^2 \right\| \stackrel{\hbox{ \eqref{eq Uf(s)} }}{=} \\
     	&= \left\|p - \left(U_p\left(c^{-1/2}\circ U_f(g)\right)\right)^2 \right\|  \leq 2 \left\| p - U_p\left(c^{-1/2}\circ U_f(g)\right) \right\| 
     \end{aligned}$$ \vspace*{-4mm}
     $$
     \begin{aligned}
     \hspace*{15mm} 	&= 2 \left\| U_p\left( f - c^{-1/2}\circ U_f(g)\right) \right\|  \leq 2 \left\| f- c^{-1/2}\circ U_f(g) \right\|\stackrel{\hbox{ \eqref{eq Uf(s)} }}{=} \\
     	&= 2 \|f- U_f (s) \| \stackrel{\hbox{ \eqref{eq Uf(s)2} }}{=}  2 \|f- U_f (g)^{\frac12} \| \leq 2 \| f - U_f(g)\| = 2 \|f-g\|.
     \end{aligned}
     $$
\end{proof}

We are now in a position to establish the uniform continuity of every bounded finitely additive real measure on the lattice of projections of a JBW$^*$-algebra. We are inspired by arguments originally developed by Christensen \cite{Christensen82} for positive measures on the lattice of projections in a von Neumann algebra, which were later extended by Bunce and Wright for signed measures (see \cite[Proposition 2.3]{BunceWright1992}). Here we show that positivity is not necessary in the arguments. The technique actually relies on the possibility of halving projections, the technical results gathered along Lemmata \ref{lemma: BunceWright_isoclinic}, \ref{lemma: symproj}, \Cref{theo: mu_lin_fac}, the properties of the quasi-linear extension of the measure in \Cref{prop_quasilin}, and an appropriate Jordan version of Filmore's theorem borrowed from \cite[Porposition 1.6]{BunceWright1989}.

\begin{proposition}\label{prop: cont_mu}
	Let $\mathfrak{J}$ be a JBW$^*$-algebra without type $I_2$ summand. Then every bounded finitely additive measure  $\mu : \PP(\mathfrak{J}) \rightarrow \mathbb{R}$ is uniformly continuous on $\PP(\mathfrak{J})$. 
\end{proposition}

\begin{proof} By structure theory and the assumptions here, $\mathfrak{J}$ decomposes as the orthogonal $\ell_{\infty}$-sum $\mathfrak{J} = \mathfrak{J}_{I_{mod}} \oplus \mathfrak{J}_{II_{1}} \oplus \mathfrak{J}_{\infty}$, where $\mathfrak{J}_{I_{mod}}$ is a (possibly zero) type $I_{mod}$ JBW$^*$-algebra, $\mathfrak{J}_{II_{1}}$ is a (possibly zero) modular type $II$ JBW$^*$-algebra, and $\mathfrak{J}_{\infty}$ is a (possibly zero) properly non-modular JBW$^*$-algebra (see subsection~\ref{subsec: structure}). Since $\mu$ is additive on finite families of mutually orthogonal projections, and the distance between vectors in the above orthogonal sum is computed as the maximum among the distances between the three components, it suffices to prove that $\mu$ is uniformly continuous when restricted to the lattice of projections in each summand.\smallskip
	
In all cases, we fix two arbitrary projections $p,q \in \PP(\mathfrak{J})$ with $\|p-q\| < \delta < \frac{1}{4}$. There is no loss of generality in assuming that $p \wedge q = 0$. Set $\mathfrak{J}_{mod}= \mathfrak{J}_{I_{mod}} \oplus \mathfrak{J}_{II_{1}},$ and denote by $p_{_{mod}}$ and $q_{_{mod}}$ the components of $p$ and $q$ in $\mathfrak{J}_{mod}$, respectively.\smallskip	

\noindent$(i)$ \emph{Modular part}. Let $\tau$ be the centre-valued trace on $\mathfrak{J}_{mod}$ whose existence is guaranteed by \Cref{theo: topping_trace} and \Cref{prop: trace_exceptional}. By \Cref{remark: decomp_finite_alge}, every projection $e$ in $\mathfrak{J}_{mod}$ can be written as an orthogonal sum  $e = e_1 + e_2$, where $e_1$ can be halved, and $e_2$ is abelian. We can thus write $ p_{_{mod}} = p_{_{mod,1}} + p_{_{mod,2}}$ and $q_{_{mod}} = q_{_{mod,1}} + q_{_{mod,2}}$, where each summand satisfies the conditions commented above. We will prove the continuity of $\mu$ at both types of projections.\smallskip

Assume that $p_{_{mod,1}}\neq 0$. Since this projection can be halved, there exist two orthogonal projections $p_{_{mod,1,1}}, p_{_{mod,1,2}} \in \PP(\mathfrak{J}_{mod})$ such that $p_{_{mod,1}} = p_{_{mod,1,1}} + p_{_{mod,1,2}}$ and $p_{_{mod,1,1}} \sim_n p_{_{mod,1,2}}$. As $p_{_{mod,1,j}} \leq p_{_{mod,1}} $ for $j =1,2$ and $\|p_{_{mod}} -q_{_{mod}}\| < \delta < 1$, by \Cref{lemma: symproj}  there exist a symmetry $s_{mod},$ depending on  $p_{_{mod}}$ and $q_{_{mod}},$ such that $U_{s_{mod}}(p_{_{mod}}) =q_{_{mod}}$ and $$\| p_{_{mod,1,j}} -  U_{s_{mod}}(p_{_{mod,1,j}}) \| \leq \sqrt{2} \ \|p_{_{mod}}-q_{_{mod}}\|^{\frac12} \leq \left(2 \delta\right)^{\frac12} $$ for $j=1,2$. To simplify the notation,  for $j =1,2$, we denote $q_{_{mod,1,j}} = U_{s_{mod}}(p_{_{mod,1,j}}) \leq$ $ U_{s_{mod}}(p_{_{mod}}) $ $=q_{_{mod}}$.\smallskip 

By known properties of the lattice $\PP(\mathfrak{J})$ (see \cite[Proposition 5.2.3$(iii)$]{HOS}), combined with the assumption that $p \wedge q =0$ (and hence $ p_{_{mod,1}}\wedge q_{_{mod,1}} = 0$), we obtain  $(p_{_{mod,1}} \vee q_{_{mod,1}})-q_{_{mod,1}} \sim p_{_{mod,1}} -p_{_{mod,1}} \wedge q_{_{mod,1}} = p_{_{mod,1}}.$ Therefore $$\tau(p_{_{mod,1}} \vee q_{_{mod,1}}) = \tau(p_{_{mod,1}}) + \tau(q_{_{mod,1}}) = 2 \tau(p_{_{mod,1}})$$ as $q_{_{mod,1}} \sim_1 p_{_{mod,1}}$. Similarly, since  $p_{_{mod}}\wedge q_{_{mod}} = 0$ and $q_{_{mod}} = U_{s_{mod}}(p_{_{mod}})$ we have $p_{_{mod,1,j}} \wedge U_{s_{mod}}(p_{_{mod,1,j}}) = 0 = p_{_{mod,1,j}} \wedge q_{_{mod,1,j}}$, and thus $$\tau (p_{_{mod,1,j}} \vee q_{_{mod,1,j}}) = 2 \tau(p_{_{mod,1,j}}) = \tau(p_{_{mod,1}}),$$ for $j = 1,2$ as $p_{_{mod,1,1}} \sim_n p_{_{mod,1,2}}$. Note now that 
$$
	(\unit-p_{_{mod,1,j}}) \wedge (\unit-q_{_{mod,1,j}}) = \unit - (p_{_{mod,1,j}} \vee q_{_{mod,1,j}}) \hspace{2.6cm}$$
	$$\hspace{3.6cm} \geq (p_{_{mod,1}}\vee q_{_{mod,1}}) - (p_{_{mod,1,j}} \vee q_{_{mod,1,j}}),$$ 
which implies that 
$$
\begin{aligned}
	\tau((\unit-p_{_{mod,1,j}})\wedge (\unit-q_{_{mod,1,j}})) & \geq 2 \tau(p_{_{mod,1}}) - \tau(p_{_{mod,1}}) = \tau(p_{_{mod,1}}) = 2 \tau(p_{_{mod,1,j}}),
\end{aligned}
$$ for $j = 1,2$. 
Thus, by \Cref{lemma: centre_trace}  $$p_{_{mod,1,j}} \lesssim (\unit-p_{_{mod,1,j}}) \wedge (\unit-q_{_{mod,1,j}}) \ \ (j=1,2).$$ 

Pick a projection $e \leq (\unit-p_{_{mod,1,j}}) \wedge (\unit-q_{_{mod,1,j}})$ with $e\sim p_{_{mod,1,1}}$.
By \Cref{lemma: BunceWright_isoclinic} (applied to $p_{_{mod,1,1}}$, $p_{_{mod,1,j}}$ and $e$),  there exist projections $h_{_{mod,j}}$ with $j = 1,2$ such that $h_{_{mod,j}}$ is isoclinic with certain angle $\theta$ to both $p_{_{mod,1,1}}$ and $p_{_{mod,1,j}}$ and 
$$ \| h_{_{mod,j}} -p_{_{mod,1,j}}\| , \|h_{_{mod,j}} - q_{_{mod,1,j}}\| \leq \| p_{_{mod,1,j}} - q_{_{mod,1,j}}\| < (2 \delta)^{\frac12}.$$
Since, $h_{_{mod,j}}$ is isoclinic to both $p_{_{mod,1,1}}$ and $p_{_{mod,1,j}}$, \Cref{l BunceWright lemma 2.1} and \Cref{theo: mu_lin_fac} assure that the quasi-linear extension, $\overline{\mu}$, of $\mu$ is linear on $W^*(\unit, h_{_{mod,j}}, p_{_{mod,1,j}}) $ and on $W^*(\unit, h_{_{mod,j}}, q_{_{mod,1,j}})$. Consequently, by \Cref{prop_quasilin}$(c)$, for $j = 1,2,$
we have \begin{equation}\label{eq one june3} \begin{aligned}
		|\mu(p_{_{mod,1,j}}) - \mu(q_{_{mod,1,j}}&)| \leq |\mu(p_{_{mod,1,j}}) - \mu(h_{_{mod,j}})| + |\mu(q_{_{mod,1,j}})- \mu(h_{_{mod,j}})| \\
		 = |&\overline{\mu}(p_{_{mod,1,j}} -h_{_{mod,j}})| + |\overline{\mu}(q_{_{mod,1,j}} - h_{_{mod,j}})|  \\ 
		\leq (&2 \alpha_{_{\mu}}(\unit) - \mu(\unit)) (\| h_{_{mod,j}} -p_{_{mod,1,j}}\| + \|h_{_{mod,j}} - q_{_{mod,1,j}}\|)\\
		 < 2& (2 \alpha_{_{\mu}}(\unit) - \mu(\unit)) (2\delta)^{\frac12},   \ \ (j=1,2).
	\end{aligned}
\end{equation} Having in mind that $\mu$ is finitely additive we derive from \eqref{eq one june3} that 
\begin{equation}\label{eq pmod1} | \mu(p_{_{mod,1}}) - \mu(q_{_{mod,1}}) | \leq 4 (2 \alpha_{_{\mu}}(\unit) - \mu(\unit)) (2\delta)^{\frac12}.
\end{equation}
	
We deal next with the abelian parts of $p$ and $q$. The projections $p_{_{mod,2}}$ and $q_{_{mod,2}}$ are abelian by construction.  We claim that we can reduce to the case $c(p_{_{mod,2}}) =c(q_{_{mod,2}}) = \unit_{_{mod}}$ and $\mathfrak{J}_{mod}= \mathfrak{J}_{I_{mod}}$. Namely, consider the decomposition of $\mathfrak{J}_{mod}$ (and of $\mathfrak{J}$) induced by $c(p_{_{mod,2}})$ in the form  $$U_{c(p_{_{mod,2}})}(\mathfrak{J}_{mod})\oplus U_{\unit-c(p_{_{mod,2}})}(\mathfrak{J}_{mod}) = c(p_{_{mod,2}})\circ \mathfrak{J}_{mod}\oplus (\unit-c(p_{_{mod,2}})) \circ \mathfrak{J}_{mod},$$ where clearly  $p_{_{mod,2}} \in U_{c(p_{_{mod,2}})}(\mathfrak{J}_{mod})$. Note that $q_{_{mod,2}} \in U_{c(p_{_{mod,2}})}(\mathfrak{J}_{mod})$ since $\| p_{_{mod,2}} -q_{_{mod,2}}\| < \delta$. Therefore $p_{_{mod,2}}, q_{_{mod,2}}$ belong to $U_{c(p_{_{mod,2}})}(\mathfrak{J}_{mod}),$ and moreover, $c(q_{_{mod,2}}) \leq c(p_{_{mod,2}}) $. The statement in the claim follows by symmetry.\smallskip 

By assumptions $\displaystyle \mathfrak{J}_{mod} = \bigoplus_{\mathbb{N}\ni n\geq 3} z_{n}\circ \mathfrak{J}_{mod}$, for a suitable family of (possibly zero) central projections $(z_{n})_{n}$ in $\mathfrak{J}$ such that $z_{n}\circ \mathfrak{J}_{mod} = z_{n}\circ \mathfrak{J}$ is of type $I_{n}$. We write $\displaystyle p_{_{mod,2}} = \sum_{n\geq 3}^{\infty} p\circ z_n$ and $\displaystyle q_{_{mod,2}} = \sum_{n\geq 3}^{\infty} q\circ z_n$, with respect to the previous decomposition. We consider just those natural numbers $n$ for which $p\circ z_n, q\circ z_n\neq 0$.\smallskip

JBW$^*$-algebras of type I are completely known through the classification of their self-adjoint parts, the JBW-algebras of type $I$ which can be found in \cite{Shultz79, Stacey81, Stacey82, Topping65}. Actually, for each natural $n$, the JBW$^*$-algebra $z_n\circ \mathfrak{J}_{mod}$ is isometrically isomorphic to $C(\Omega_n,\mathfrak{F}_n)$, where $\Omega_n$ is a hyperStonean compact Hausdorff space and $\mathfrak{F}_n$ is a factor JBW$^*$-algebra of rank $n\geq 1$. Note that, by hypotheses, $n\neq 2$. It is known that every type JBW$^*$-algebra of type $I_1$ is isometrically isomorphic to a commutative von Neumann algebra. So, on the type $I_1$ part of $\mathfrak{J}$ the conclusion is a of the proposition is a trivial consequence of the Mackey-Gleason theorem. We can therefore assume that $n\geq 3$.\smallskip  

Since $p_n = z_n\circ p_{_{mod,2}}$ and  $q_n = z_n\circ q_{_{mod,2}}$ are abelian projections in $z_n\circ \mathfrak{J}_{mod}$ and $n\geq 3$, it can be easily deduced that $p_n\sim_1 q_n$ with $p_n,q_n \lesssim (z_n-p_n)\wedge (z_n-q_n)$ in $z_n\circ \mathfrak{J}_{mod}$. By glueing all symmetries in the corresponding summands $z_n\circ \mathfrak{J}_{mod}$, we derive that $p_{_{mod,2}}\sim q_{_{mod,2}}$ and $p_{_{mod,2}},q_{_{mod,2}} \lesssim (\unit_{_{mod}}- p_{_{mod,2}})\wedge (\unit_{_{mod}}- q_{_{mod,2}})$ in $\mathfrak{J}_{mod}$. Taking $e\leq (\unit_{_{mod}}- p_{_{mod,2}})\wedge (\unit_{_{mod}}- q_{_{mod,2}}),$ with $p_{_{mod,2}}\sim e$ we are in a position to apply \Cref{lemma: BunceWright_isoclinic}. We therefore find a new projection $z_{_{mod,2}}\in \mathfrak{J}_{mod}$ which is  isoclinic to $p_{_{mod,2}}$ and $q_{_{mod,2}}$ and 
$$ \| z_{_{mod,2}} - p_{_{mod,2}} \| , \|z_{_{mod,2}}  - q_{_{mod,2}}\| \leq \| p_{_{mod,2}} - q_{_{mod,2}} \| < \delta.$$
Having in mind that $z_{_{mod,2}}$ is isoclinic to both $p_{_{mod,2}}$ and $q_{_{mod,2}}$, \Cref{l BunceWright lemma 2.1} and \Cref{theo: mu_lin_fac} guarantee that the quasi-linear extension, $\overline{\mu}$, of $\mu$ is linear on the JBW$^*$-subalgebras $W^*(\unit, z_{_{mod,2}}, p_{_{mod,2}}) $ and $W^*(\unit, z_{_{mod,2}}, q_{_{mod,2}})$. \Cref{prop_quasilin}$(c)$ now gives \begin{equation}\label{eq two June 3} \begin{aligned}
		|\mu(p_{_{mod,2}}) - \mu(q_{_{mod,2}})| &\leq |\mu(p_{_{mod,2}}) - \mu(z_{_{mod,2}})| + |\mu(z_{_{mod,2}})- \mu(q_{_{mod,2}})| \\
		 =& |\overline{\mu}(p_{_{mod,2}} - z_{_{mod,2}})| + |\overline{\mu}(q_{_{mod,2}} - z_{_{mod,2}})|  \\ 
		\leq& (2 \alpha_{_{\mu}}(\unit) - \mu(\unit)) (\| z_{_{mod,2}} - p_{_{mod,2}} \| + \|z_{_{mod,2}}  - q_{_{mod,2}}\|)\\
		 <& 2 (2 \alpha_{_{\mu}}(\unit) - \mu(\unit)) \delta.
	\end{aligned}
\end{equation}

Finally, by combining \eqref{eq one june3}, \eqref{eq two June 3} and the orthogonal additivity of $\mu$, we arrive to \begin{equation}\label{eq boundedness on the modular part}  | \mu(p_{_{mod}}) - \mu(q_{_{mod}})| \leq 2 (2 \alpha_{_{\mu}}(\unit) - \mu(\unit)) (\delta + (2\delta)^{\frac12}).
\end{equation}

\noindent$(ii)$ \emph{Properly non-modular part}. We deal next with projections in the properly non-modular part of $\mathfrak{J}$. Let $p_{_\infty},q_{_\infty}$ denote the components of $p$ and $q$ in  $\mathfrak{J}_{\infty}$.  The arguments are close  to those employed in the type $II_1$ part. We shall briefly comment the steps to follow. Observe that $\| p_{_\infty}- q_{_\infty}\|<\delta <\frac14$.\smallskip 

Let $s\in \mathfrak{J}_{\infty}$ be the symmetry given by \Cref{lemma: symproj} for $p_{\infty}$ and $q_{\infty}$. Set $x = U_{p_{_\infty}}(s)\in \mathfrak{J}_{\infty}$. We shall apply now a Jordan version of a classical result by Fillmore obtained by Bunce and Wright in \cite{BunceWright1989}. Note that by structure theory, $\mathfrak{J}_{\infty}$ is a JW$^*$-algebra (cf. \cite[Theorem 5.3.9]{HOS}). Concretely, Proposition 1.6 in \cite{BunceWright1989}, applied to $U_{p_{\infty}} (\mathfrak{J}_{\infty})$ and the element $x$, assures the existence of two orthogonal projections $p_1,p_2$ in $U_{p_{\infty}} (\mathfrak{J}_{\infty})$ satisfying $p_{\infty}= p_1 + p_2$,  $p_1\sim_1 p_2$, and $p_j$ and $x$ operator commute in $U_{p_{\infty}} (\mathfrak{J}_{\infty})$ for all $j=1,2$. By applying the latter property we get $U_{p_1,p_2} (x) =  2 (p_1\circ x)\circ p_2 - (p_1\circ p_2)\circ x=0$, and hence $$x= U_{p_{\infty}} (x) = U_{p_1} (x) + U_{p_2} (x) = U_{p_1} ( U_{p_{_\infty}}(s)) + U_{p_2} ( U_{p_{_\infty}}(s))= U_{p_1} (s) + U_{p_2} (s).$$ Take now $q_j = U_s(p_j)$ ($j=1,2$). The properties of the symmetry $s$ given in \Cref{lemma: symproj} imply that $\|p_j-q_j\| \leq (2 \|p_{\infty}-q_{\infty}\|)^{\frac12}<   (2 \delta)^{\frac12},$ for all $j=1,2$. By construction, we also have $p_i\circ q_j = 0$ for all $i \neq j$ in $\{1,2\}$. Therefore, $p_i \sim_1 p_j \leq (1-p_i) \wedge (1-q_i)$ in $\mathfrak{J}_{\infty}$ with $p_j \circ q_i =0$ for all $i \neq j$ in $\{1,2\}$. A new application of \Cref{lemma: BunceWright_isoclinic} assures the existence of projections $\tilde{h}_j\in \mathfrak{J}_{\infty}$ ($j=1,2$) satisfying that $h_j$ is isoclinic to both $p_j$ and $q_j$ and $\|p_j-h_j\|,\|q_j-h_j\| < \|p_j-q_j\|<(2 \delta)^{\frac12}$ for all $j=1,2$. Now by chaining \Cref{l BunceWright lemma 2.1}, \Cref{theo: mu_lin_fac}, and \Cref{prop_quasilin}$(c)$, the same trick used to obtain \eqref{eq pmod1} allows us to conclude that $$| \mu(p_{\infty}) - \mu(q_{\infty}) | \leq 4 (2 \alpha_{_{\mu}}(\unit) - \mu(\unit)) (2\delta)^{\frac12},$$ which concludes the proof.\end{proof}

The Mackey-Gleason theorem for signed measures on factors JBW$^*$-algebras of type $I_n$ with $n\neq 2$ in  \Cref{theo: mu_lin_fac}, together with the conclusion in the previous \Cref{prop: cont_mu}, allow us to apply a classic method developed by Christensen and Maeda (see \cite[Theorem 8.6]{Maeda1989}). As commented by Bunce and Wright in the case of von Neumann algebras (see \cite[comments prior to Proposition 2.4]{BunceWright1992}), the argument does not depend on the positivity of the measure, and as shown below, nor on the associativity.

\begin{theorem}{\label{theo: lin_type_In}}
Let $\mathfrak{J}$ be a JBW$^*$-algebra of type $I_n$ with $n \neq  2$ and $n < \infty$. Then every bounded finitely additive signed measure $\mu: \PP(\mathfrak{J}) \rightarrow \mathbb{R}$ can be extended to a bounded linear functional on $\mathfrak{J}$. Moreover, the quasi-linear extension, $\overline{\mu}$, of $\mu$ is linear.
\end{theorem}

\begin{proof} Observe that, by a new application of the structure theory of JBW$^*$-algebras, we can assume $\mathfrak{J}_{sa} = C(\Omega,\mathfrak{F}_n)$, where $\Omega$ is a hyper-Stonean space and $\mathfrak{F}_n$ is a factor JBW$^*$-algebra of rank $n\geq 3$. Recall that every type $I_1$ JBW$^*$-algebra is a commutative von Neumann algebra, and thus the desired result is a consequence of the Mackey-Gleason-Bunce-Wright theorem proved by Bunce and Wright in \cite[Proposition 1.1]{BunceWright1992}. If we discard the case $n=1$, it follows from the hypotheses that $n\geq 3$. We further know that in such a case $(\mathfrak{F}_{n})_{sa}$ is Jordan isomorphic to $H_3(\mathbb{O})$ or to $M_n(\mathbb{F}),$ where $\mathbb{F} = \mathbb{R}, \mathbb{C}$ or $\mathbb{H}$ (cf. \cite[Theorem 5.3.8, Remark 6.4.3]{HOS} and \cite[Theorem 3.9]{Shultz79}).\smallskip
	 
We denote by $C^0(\Omega,\mathfrak{F}_{n})$ the subspace of all locally constant mappings in the sense introduced in \cite[Definition 8.2]{Maeda1989}, that is, the functions $a\in C(\Omega,\mathfrak{F}_{n})$ whose image is finite. Clearly $C^0(\Omega,\mathfrak{F}_{n})$ is a Jordan $^*$-subalgebra of $\mathfrak{J}$. \smallskip

We claim that $\overline{\mu}$ is linear on $C^0(\Omega,\mathfrak{F}_{n})$. To prove the claim, let us take $a, b \in C^0(\Omega,\mathfrak{F}_{n})$, and a partition $\{K_1, \dots, K_m\} $ of $\Omega$ such that each $K_i$ is clopen and both $a$ and $b$ are constant mappings on each $K_i$ for all $i:1,\ldots, m$. Let $\chi_{_{K_i}}$ denote the characteristic function of the set $K_i$. The finite dimensional JB$^*$-subalgebra 
$\mathfrak{B}:= \bigoplus_{i=1}^m \mathfrak{F}_{n} \chi_{_{K_i}}\cong \bigoplus_{i=1,\ldots,m}^{\infty} \mathfrak{F}_{n}.$ Since $\mathcal{P} (\mathfrak{B}) = \bigoplus_{i=1}^m \mathcal{P}(\mathfrak{F}_{n}) \chi_{_{K_i}},$ it is not hard to check via \Cref{theo: mu_lin_fac} and the orthogonal additivity of $\mu$ that $\mu|_{\mathcal{P}(\mathfrak{B})}$ admits a unique extension to a bounded linear functional on $\mathfrak{B},$ which, by uniqueness, must coincide with $\overline{\mu}|_{\mathfrak{B}}$ (cf. \Cref{prop_quasilin}). Since $a, b \in \mathfrak{B}$, we get $\overline{\mu}(\gamma a + \beta b) = \gamma \overline{\mu}(a) + \beta \overline{\mu}(b)$ ($\gamma, \beta \in \mathbb{C})$. Therefore $\overline{\mu}$ is linear on $C^0(\Omega,\mathfrak{F}_{n})$.\smallskip

\Cref{prop_quasilin} now assures that $\overline{\mu}|_{C^0(\Omega,\mathfrak{F}_{n})}$ is a bounded linear functional with $\|\overline{\mu}|_{C^0(\Omega,\mathfrak{F}_{n})}\|\leq 2 (2 \alpha_{_{\mu}}(\unit) - \mu(\unit))$. Having in mind the norm density of $C^0(\Omega,\mathfrak{F}_{n})$ in $\mathfrak{J}$, we can find a unique extension of $\overline{\mu}|_{C^0(\Omega,\mathfrak{F}_{n})}$ to a bounded linear functional $\varphi: \mathfrak{J}\to \mathbb{C}$. By construction, $\varphi = \overline{\mu} = \mu$ on $\PP(C^0(\Omega,\mathfrak{F}_{n}))$, and the latter is dense in $\PP(\mathfrak{J})$. We can finally apply the uniform continuity of $\mu$ on $\PP(\mathfrak{J})$ (cf. \Cref{prop: cont_mu}) to deduce that $\varphi = \mu$ on $\PP(\mathfrak{J})$, and thus $\varphi$ is a bounded linear extension of $\mu$, and $\varphi = \overline{\mu}$.
\end{proof}

\begin{remark}\label{r after theorem type In} Let $\mu: \PP(\mathfrak{J}) \rightarrow \mathbb{R}$ be a bounded finitely additive signed measure on the lattice of projections of a JBW$^*$-algebra $\mathfrak{J}$. By combining \Cref{prop_quasilin} and the previous \Cref{theo: lin_type_In}, we can conclude that the quasi-linear extension of $\mu$ is linear when restricted to each JBW$^*$-subalgebra of type $I_n$ with $n\neq 2$ of $\mathfrak{J}$. 
\end{remark}

Another straightforward consequence of \Cref{theo: lin_type_In} assure that if $\mathfrak{J}$ is a finite sum of JBW$^*$-algebras of type $I_n$ with $2\neq n\in \mathbb{N}$, then every bounded finitely additive measure $\mu : \mathcal{P} (\mathfrak{J})\to \mathbb{R}$ can be extended to a bounded linear functional on $\mathfrak{J}$. The main result of this section concludes the case of type I modular JBW$^*$-algebras without type $I_2$ part, that is, arbitrary direct sums of JBW$^*$-algebras of type $I_n$ with $2\neq n\in \mathbb{N}$. We recall that a JBW$^*$-algebra of type $I$ has bounded dimension of irreducible representations if it can be written as a finite direct sum of JBW$^*$-algebras of type $I_n$ with $n$ finite. 

\begin{theorem}{\label{theo: lin_type_I}} Let $\mathfrak{J}$ be a type $I$ modular JBW$^*$-algebra with no direct summands of type $I_2$.  Let $\mu: \PP(\mathfrak{J}) \rightarrow \mathbb{R}$ be a bounded finitely additive signed measure. Let $\mathfrak{J}_0$ be a JBW$^*$-subalgebra of $\mathfrak{J}$ which is of type $I$ (modular) with no direct summands of type $I_2$ and has bounded dimension of irreducible representations. Then the restriction of $\mu$ to $\PP(\mathfrak{J}_0)$ can be extended to a bounded linear functional on $\mathfrak{J}_0$. Moreover, there exists a bounded linear functional $\overline{\mu}_0$ which coincides with $\mu$ on every projection belonging to a finite direct sum of JBW$^*$-algebra summands of type $I_n$ of $\mathfrak{J}$.  
\end{theorem}

\begin{proof} By assumptions, there exist (possibly zero) pairwise orthogonal central projections $z_n$  ($n \in \mathbb{N}$) satisfying that $z_n\circ\mathfrak{J}$ is of type $I_n,$ $\displaystyle \unit = \sum_{n=1}^{\infty} z_n$, and $z_2 = 0$. \Cref{theo: lin_type_In} (see also \Cref{r after theorem type In}) proves that $\overline{\mu}$ is linear on $z_n\circ \mathfrak{J}$ for each $n \in \mathbb{N}$, equivalently, the assignment $a\mapsto \overline{\mu}(z_n\circ a)$ defines a bounded linear functional $\overline{\mu}_n$ on $\mathfrak{J}$ for all natural $n$.\smallskip
	
Now, given a natural $n$, it follows from \Cref{prop_quasilin}$(c)$ that the inequality  
$ | \overline{\mu}(z_n\circ x)| \leq  \alpha_{_\mu} (z_n) + (\alpha_{_\mu} (z_n) -\mu (z_n))\leq \alpha_{_\mu} (z_n),$ for every $x = x^*$ and $\|x\| \leq 1$, and consequently, $\|\overline{\mu}_n\|\leq 2 \alpha_{_\mu} (z_n)$, for all natural $n$. By applying \Cref{r properties of alpha and V} we have  
$$\begin{aligned}
	\sum_{n= 1}^{k} \|\overline{\mu}_n \| &\leq \sum_{n= 1}^{k} 2 \alpha_{_\mu} (z_n) = 2 \alpha_{_\mu}\left(\sum_{n=1}^{k}z_n\right)  \leq 2 \alpha_{_\mu}(\unit),
\end{aligned}$$ for all natural $k$, which implies that the series $\displaystyle \sum_{n=1}^{\infty} \overline{\mu}_n$ is absolutely convergent in $\mathfrak{J}^*$. Take $\overline{\mu}_0 = \displaystyle \sum_{n=1}^{\infty} \overline{\mu}_n \in \mathfrak{J}$. Given a projection $\displaystyle p$ in $\displaystyle \bigoplus_{k=n_1,\ldots, n_m} z_k \circ \mathfrak{J}$, we have $$\mu (p) =\mu \left( \sum_{k=n_1,\ldots, n_m} p\circ z_k\right) = \sum_{k=n_1,\ldots, n_m} \mu (p\circ z_k) = \sum_{k=1}^{k} \overline{\mu}_k (p) = \overline{\mu}_0 (p).$$

If $\mathfrak{J}_0$ is a JBW$^*$-subalgebra of $\mathfrak{J}$ which is of type $I$ with no direct summands of type $I_2$ and has bounded dimension of irreducible representations, there exists a finite set of central projections $z_{n_1},\ldots, z_{n_k}$ in $\mathfrak{J}_0$ such that $\mathfrak{J}_0 = \displaystyle \bigoplus_{k=n_1,\ldots, n_m} z_k \circ \mathfrak{J}$, where $z_k \circ \mathfrak{J}$ is of type $I_{n_k}$. So the previous arguments give the desired statement. 
\end{proof}

\section{Properly non-modular \texorpdfstring{JW$^*$}{JW*}-algebras}\label{sec: properly non-modular}

In this section we study bounded finitely additive measures on the lattice of projections of properly non-modular JW$^*$-algebras. Let us recall that a Jordan subalgebra $\mathfrak{J}$ of an associative algebra $A$ is called \emph{reversible} if for all natural $n$  we have $$a_1,a_2,\ldots, a_n\in \mathfrak{J} \Rightarrow a_1 a_2 \ldots a_n + a_n \ldots a_{2} a_1 \in \mathfrak{J},$$ where the juxtaposition of two or several elements stands for their associative product in $A$ (cf. \cite[Definition 4.24]{AlfsenShultz2003}).  Reversible JW$^*$-algebras are well determined (see \cite[Theorem 4.29]{AlfsenShultz2003}). We simply observe that $3$-dimensional spin factors are examples of non-reversible JW$^*$-algebras \cite[Theorem 4.31]{AlfsenShultz2003}).\smallskip 

We begin this section with a technical result needed in our arguments, which is a novelty by itself in the Jordan setting. 

\begin{proposition}\label{prop a la Christensen} Let $\mathfrak{J}$ be a JW$^*$-algebra. Let $p, q_1, q_2, q_3$ be mutually orthogonal projections in $\mathfrak{J}$ such that $p\sim_1 q_i$ for all $i\in \{1,2,3\}.$ Let us additionally assume that $s_i$ is a symmetry in $\mathfrak{J}$ such that $U_{s_i} (p) = q_i.$ Then, for every $a,b$ in $\mathfrak{J}$ with $0\leq a,b\leq \frac12 p$,  the elements \begin{equation}\label{eq r} \begin{aligned}
r =& a + 2 a\circ s_1 + U_{s_1} (a) + U_{s_2} (p-2 a) + 2 U_{s_2,a^{\frac12}} \left( \left( p-2 a \right)^{\frac12} \right) \\
&+  2 U_{s_1,s_2} \left( \left( a\circ (p-2 a) \right)^{\frac12} \right),
\end{aligned}
\end{equation} and \begin{equation}\label{eq q}\begin{aligned}
q =& b - 2 b\circ s_1 + U_{s_1} (b) + U_{s_3} (p-2 b) + 2 U_{s_3,b^{\frac12}} \left( \left( p-2 b \right)^{\frac12} \right) \\
&-  2 U_{s_1,s_3} \left( \left( b\circ (p-2 b) \right)^{\frac12} \right),    
\end{aligned}
\end{equation} are two orthogonal projections in $\mathfrak{J}$ satisfying $U_{p} (r) = a$ and $U_p (q) = b.$\smallskip

\noindent Furthermore, if $\mathfrak{J}$ is a reversible JW$^*$-algebra and $p, q_1, q_2, q_3,q_4,$ and $q_5$ are non-trivial mutually orthogonal projections in $\mathfrak{J}$ such that $p\sim_1 q_i$ for all $i\in \{1,2,3,4,5\},$ then for every $c,d$ in $\mathfrak{J}$ with $0\leq c,d$, $ c+d\leq p$, there exist orthogonal projections $\tilde{r},\tilde{q}\in \mathfrak{J}$ satisfying $U_{p} (\tilde{r}) = c$ and $U_p (\tilde{q}) = d.$
\end{proposition}

\begin{proof} We can assume, by hypotheses, that $\mathfrak{J}$ is a JBW$^*$-subalgebra of some von Neumann algebra $A$. Henceforth, the associative product of $A$ will be denoted by mere juxtaposition. We can clearly assume that $\mathfrak{J}$ and $A$ share the same unit, and thus each $s_i$ is a symmetry in $A$, and each $q_i$ is a projection in $A$.\smallskip    

Let us deal with the first statement. The elements $$x = a^{\frac12} + a^{\frac12} s_1 + (p-2 a)^{\frac12} s_2, \hbox{ and } y = b^{\frac12} + b^{\frac12} s_1 + (p-2 b)^{\frac12} s_3,$$ are not necessarily in $\mathfrak{J}$, but both lie in the von Neumann algebra $A$. Working in $A$, it is not hard to check that $xx^* = yy^* =p$, $x^* x = r,$ $y^* y = q,$ and $ x y^*= y^* x =0.$ Therefore $x$ and $y$ are partial isometries in $A$, and thus $r$ and $q$ are orthogonal projections in $A$. Since, clearly, $r$ and $q$ are defined in terms of Jordan products of elements in $\mathfrak{J}$ (cf. \eqref{eq q} and \eqref{eq q}), and they are thus projections in $\mathfrak{J}$, which concludes the proof of the first statement.\smallskip

For the second statement, we rely on the extra hypotheses. We keep the notation above, and we assume that $\mathfrak{J}$ is a reversible JW$^*$-algebra inside the von Neumann algebra $A$. Let $p, q_1,\ldots, q_5$ and $s_1,\ldots, s_5$ be the mutually orthogonal projections and the symmetries in $\mathfrak{J}$ given by the hypotheses (that is, $U_{s_i} (p) = q_i$ for all $i=1,\ldots, 5$). As in the the first part of the result, we consider the elements 
$$z = \left(\frac{c}{2}\right)^{\frac12} + \left(\frac{d}{2}\right)^{\frac12} s_1 + \left(\frac{c}{2}\right)^{\frac12} s_2+ \left(\frac{d}{2}\right)^{\frac12} s_3+ \left(p- c-d\right)^{\frac12} s_4, \hbox{ and  }$$ $$t = \left(\frac{c}{2}\right)^{\frac12} - \left(\frac{d}{2}\right)^{\frac12} s_1 - \left(\frac{c}{2}\right)^{\frac12} s_2+ \left(\frac{d}{2}\right)^{\frac12} s_3+ \left(p- c-d\right)^{\frac12} s_5,$$
which are in $A$ but not necessarily in $\mathfrak{A}$. Operating in the von Neumann algebra $A$ we can check that
$z z^* =p,$ $t t^* =p,$ $ z t^* = \frac{c}{2} -\frac{d}{2}-\frac{c}{2} + \frac{d}{2}=0$, which shows that $\tilde{r} = z^* z$ and $\tilde{q} = t^* t$ are two orthogonal projections in the von Neumann algebra $A$. Moreover $U_p (\tilde{r}) = p \tilde{r} p = c$ and $U_p (\tilde{q}) = p \tilde{q} p = d$. Consider now the identity 
$$\tilde{r} = z^*z = \frac{c}{2} + 2 \left\{ \left(\frac{c}{2}\right)^{\frac12}, \left(\frac{d}{2}\right)^{\frac12}, s_1 \right\} + 2 \left(\frac{c}{2}\right) \circ s_2  +2 \left\{\left(\frac{c}{2}\right)^{\frac12}, \left(\frac{d}{2}\right)^{\frac12},s_3 \right\}$$
$$+ 2\left\{\left(\frac{c}{2}\right)^{\frac12}, \left(p-c-d\right)^{\frac12}, s_4  \right\} + \left\{s_1, \frac{d}{2}, s_1\right\} + 2 \left\{s_1, \frac{d}{2}, s_3\right\}
$$ 
$$+ \left\{s_3, \frac{d}{2}, s_3\right\} + \left\{s_2, \frac{c}{2}, s_2\right\}  +  \left\{s_4, (p -c-d) , s_4\right\} \hspace*{1.75cm} $$ $$ +\! s_1 \left(\frac{d}{2}\right)^{\frac12} \left(\frac{c}{2}\right)^{\frac12}\! s_2 + s_2 \left(\frac{c}{2}\right)^{\frac12} \left(\frac{d}{2}\right)^{\frac12}\! s_1 
+ s_3 \left(\frac{d}{2}\right)^{\frac12} \left(\frac{c}{2}\right)^{\frac12}\! s_2 + s_2 \left(\frac{c}{2}\right)^{\frac12} \left(\frac{d}{2}\right)^{\frac12}\! s_3 $$ 
$$ +s_1 \left(\frac{d}{2}\right)^{\frac12} \left(p-c-d\right)^{\frac12} s_4 + s_4 \left(p-c-d\right)^{\frac12} \left(\frac{d}{2}\right)^{\frac12} s_1 \hspace{1.3cm} $$
$$+s_2 \left(\frac{c}{2}\right)^{\frac12} \left(p-c-d\right)^{\frac12} s_4 + s_4 \left(p-c-d\right)^{\frac12} \left(\frac{c}{2}\right)^{\frac12} s_2 \hspace{1.49cm} $$ $$ +s_3 \left(\frac{d}{2}\right)^{\frac12} \left(p-c-d\right)^{\frac12} s_4 + s_4 \left(p-c-d\right)^{\frac12} \left(\frac{d}{2}\right)^{\frac12} s_3.\hspace{1.2cm}$$ Clearly the summands in the first three lines of the previous idenity lie in $\mathfrak{J}$ because they are given by Jordan products of elements in $\mathfrak{J}$, while the reversibility of $\mathfrak{J}$ in $A$ implies that the summands in the last four lines are in $\mathfrak{J}$. Therefore $\tilde{r} = z^* z\in \mathfrak{J}.$\smallskip   
 
Similarly, the identity 
$$\tilde{q} = t^*t = \frac{c}{2} - 2 \left\{ \left(\frac{c}{2}\right)^{\frac12}, \left(\frac{d}{2}\right)^{\frac12}, s_1 \right\} - 2 \left(\frac{c}{2}\right) \circ s_2  +2 \left\{\left(\frac{c}{2}\right)^{\frac12}, \left(\frac{d}{2}\right)^{\frac12},s_3 \right\} $$
$$+ 2\left\{\left(\frac{c}{2}\right)^{\frac12}, \left(p-c-d\right)^{\frac12}, s_5  \right\} + \left\{s_1, \frac{d}{2}, s_1\right\} - 2 \left\{s_1, \frac{d}{2}, s_3\right\} + \left\{s_3, \frac{d}{2}, s_3\right\} \hspace{0.3cm}$$  
 $$+ \left\{s_2, \frac{c}{2}, s_2\right\}  +  \left\{s_5, (p -c-d) , s_5\right\} +s_1 \left(\frac{d}{2}\right)^{\frac12} \left(\frac{c}{2}\right)^{\frac12} s_2 + s_2 \left(\frac{c}{2}\right)^{\frac12} \left(\frac{d}{2}\right)^{\frac12} s_1 $$
$$- s_3 \left(\frac{d}{2}\right)^{\frac12} \left(\frac{c}{2}\right)^{\frac12} s_2 - s_2 \left(\frac{c}{2}\right)^{\frac12} \left(\frac{d}{2}\right)^{\frac12} s_3  -s_1 \left(\frac{d}{2}\right)^{\frac12} \left(p-c-d\right)^{\frac12} s_5 \hspace{1.6cm} $$ $$- s_5 \left(p-c-d\right)^{\frac12} \left(\frac{d}{2}\right)^{\frac12} s_1 -s_2 \left(\frac{c}{2}\right)^{\frac12} \left(p-c-d\right)^{\frac12} s_5 - s_5 \left(p-c-d\right)^{\frac12} \left(\frac{c}{2}\right)^{\frac12} s_2$$ $$ +s_3 \left(\frac{d}{2}\right)^{\frac12} \left(p-c-d\right)^{\frac12} s_5 + s_5 \left(p-c-d\right)^{\frac12} \left(\frac{d}{2}\right)^{\frac12} s_3, \hspace{3.6cm}$$ implies that $\tilde{q}$ is a projection in $\mathfrak{J}$.
\end{proof}

\subsection{Bunce-Wright equivalence} \  \smallskip

If $\mathfrak{J}\subseteq B(H)$ is a JW$^*$-algebra, we write $R(\mathfrak{J}_{sa})^{-}$ for the weak$^*$-closure of the real norm closed algebra generated by $\mathfrak{J}_{sa}$. It is known that if $\mathfrak{J}$ is reversible, the equality $R(\mathfrak{J}_{sa})^{-}_{sa} = \mathfrak{J}_{sa}$ holds (\cite[Remark 2.5]{Stormer65}, \cite[1.1]{Ayupov82}).\smallskip

The statement in the previous \Cref{prop a la Christensen} has been given in terms of the usual Jordan equivalence of projections, $\sim$, recalled in page~\pageref{def Jordan equivalence}. For the goals in this section the usual Jordan equivalence will not be the appropriate tool. We shall consider a weaker relation introduced by Bunce and Wright in \cite[\S 4]{BunceWright1985}. Let $\mathfrak{J}$ be a reversible JW$^*$-algebra and let $e,f$ be projections in $\mathfrak{J}$. We shall say that $e$ and $f$ are \emph{Bunce-Wright equivalent} ($e\approx f$) if there exists a partial isometry $u$ in $R(\mathfrak{J}_{sa})^{-}$ such that $e = u^*u$ and $f = uu^*$. The symbol $e\lessapprox f$ will stand to denote the existence of a partial isometry $u$ in $R(\mathfrak{J}_{sa})^{-}$ such that $e = u^*u$ and $u u^* \leq f$. It is important to note that $u u^* \in R(\mathfrak{J}_{sa})^{-}_{sa} = \mathfrak{J}_{sa}.$\smallskip

The Bunce-Wright equivalence has been employed in texts like \cite[pages 23--24]{AyupovRakhimovUsmanovBook1997} without an explicit name.\smallskip 

The next lemma gathers some properties of the Bunce-Wright equivalence established in \cite{BunceWright1985}.
 
\begin{lemma}\label{l properties of the BWsim}{\rm\cite[Lemma 4.4 and Proposition 4.5]{BunceWright1985}} Let $\mathfrak{J}$ be a reversible JW$^*$-algebra. Then,
\begin{enumerate}[$(i)$]\item For $e,f\in \PP(\mathfrak{J})$ we have $e\sim f$ implies $e\approx f$, and $\approx$ is an equivalence relation on $\PP(\mathfrak{J})$.
\item If $(e_i)_i,$ $(f_i)_i$ are both families of mutually orthogonal projections in $\mathfrak{J}$
such that $e_i \approx f_i$, for every $i$, then $\displaystyle \sum_i e_i \approx \sum_i f_i$.
\item  If $e,f\in \PP(\mathfrak{J})$ with $e\lessapprox f$ and $f\lessapprox e$, then $e\approx f$.	
\item If $e,f\in \PP(\mathfrak{J})$, then there exists a central projection $z$ in $\mathfrak{J}$ such that $e\circ z \lessapprox f\circ z$ and $f\circ (\unit-z) \lessapprox e\circ (\unit-z)$.
\item A projection $e$ in $\mathfrak{J}$ is properly non-modular if and only if any of the equivalent statements hold:
\begin{enumerate}[$(a)$]
	\item There exists an infinite sequence $(p_n)_n$ of mutually orthogonal projections in $\mathfrak{J}$ such that $\displaystyle\sum_{n=1}^{\infty} p_n = e \approx p_m$ for every $m$.
	\item There exists a projection $f$ in $\mathfrak{J}$ with $e\gtrapprox f\approx e-f.$
\end{enumerate}
\end{enumerate}	
\end{lemma}

The next corollary is a consequence of the above properties. 

\begin{corollary}\label{c properly modular is preserved by BW equivalence} Let $p$ and $q$ be projections in a reversible JW$^*$-algebra $\mathfrak{J}$ with $p\approx q$. Then $p$ is properly non-modular if and only if $q$ is properly non-modular. 
\end{corollary}

\begin{proof} Suppose that $q$ is properly non-modular. Assume that $\mathfrak{J}$ is reversible in $B(H)$. The associative product on $B(H)$ is denoted by juxtaposition. By \Cref{l properties of the BWsim} $(v)$, there exists a projection $q_1$ in $\mathfrak{J}$ with $q\gtrapprox q_1\approx q-q_1.$ Take a partial isometry $u\in R(\mathfrak{J}_{sa})^{-}$ such that $u u^* = p$ and $u^* u =q$. Set $u_1 = u q_1, u_2 = u (q-q_1) \in R(\mathfrak{J}_{sa})^{-}$. It is easy to check that $u_1^* u_1 = q_1 u^* u q_1 = q_1$, $u_2^* u_2 = (q-q_1) u^* u (q-q_1) = q-q_1$, and hence $u_1$ and $u_2$ are partial isometries, and $p_1=u_1 u_1^*, p_2 = u_2 u_2^*\in R(\mathfrak{J}_{sa})^{-}_{sa} = \mathfrak{J}_{sa}$ are two projections with $p\approx q\approx q_1\approx p_1$, $p_2 + p_1 =p$ and $p_2\approx q-q_1\approx q\approx p$ (see \Cref{l properties of the BWsim}).
\end{proof}

\begin{proposition}\label{p four BW equivalent projections} Let $p_1,p_2,p_3,$ and $p_4$ be mutually orthogonal non-zero projections in a reversible JW$^*$-algebra $\mathfrak{J}$. Then the following statements hold:
\begin{enumerate}[$(a)$]
\item Suppose there exists a partial isometry $v_{_{12}}\in R(\mathfrak{J}_{sa})^{-}$ such that $v_{_{12}}^* v_{_{12}} = p_1$ and $v_{_{12}} v_{_{12}}^* = p_2$. Assume additionally that $p_1\approx p_j,$ for all $j$ in $\{3,4\}$. Then there exists a JBW$^*$-subalgebra $\tilde{\mathfrak{B}}$ of $\mathfrak{J}$ which is Jordan $^*$-isomorphic to $S_4(\mathbb{C})$ and contains $p_1,p_2,p_3,$  $p_4$ and $w_{_{12}}=v_{_{12}}+ v_{_{12}}^*$. In particular $\tilde{\mathfrak{B}}$ contains the JB$^*$-subalgebra generated by the elements $p_1,p_2,$ and $w_{_{12}}$ which is Jordan isomorphic to $S_2(\mathbb{C})$. Moreover, if $\mathfrak{D}$ is another JBW$^*$-subalgebra of $\mathfrak{J}$ which is orthogonal to $p_1,p_2,p_3$ and $p_4$, we can also assume that the JBW$^*$-subalgebra $\tilde{\mathfrak{B}}$ is orthogonal to $\mathfrak{D}$.
\item Suppose there exist partial isometries $v_{_{12}},v\in R(\mathfrak{J}_{sa})^{-}$ such that $v_{_{12}}^* v_{_{12}} = p_1,$ $v_{_{12}} v_{_{12}}^* = p_2$, $v p_1 v^* = p_3$, $v p_2 v^* = p_4$, $v^* p_3 v = p_1$, and $v^* p_4 v = p_4$. Then $p_1\approx p_j,$ for all $j$ in $\{3,4\},$ and there exists a JBW$^*$-subalgebra $\tilde{\mathfrak{B}}$ of $\mathfrak{J}$ which is Jordan $^*$-isomorphic to $S_4(\mathbb{C})$ and contains $p_1,p_2,p_3,$  $p_4$ and $v+v^*$. In particular $\tilde{\mathfrak{B}}$ contains the JB$^*$-subalgebra generated by the elements  $p_1+p_2,p_3+p_4,$ and $v+v^*,$ and the latter is Jordan $^*$-isomorphic to $S_2(\mathbb{C})$. Moreover, if $\mathfrak{D}$ is another JBW$^*$-subalgebra of $\mathfrak{J}$ which is orthogonal to $p_1,p_2,p_3$ and $p_4$, we can also assume that the JBW$^*$-subalgebra $\tilde{\mathfrak{B}}$ is orthogonal to $\mathfrak{D}$.
\end{enumerate}

\end{proposition}

\begin{proof} $(a)$ For $3\leq j \leq 4$ let us take a partial isometry $v_{_{1j}}\in R(\mathfrak{J}_{sa})^{-}$ such that $ v_{_{1j}}^* v_{_{1j}} = p_1$ and $v_{_{1j}} v_{_{1j}}^* = p_j$. We define the following elements: $$v_{_{23}} := v_{_{13}} v_{_{12}}^*,\  v_{_{24}} := v_{_{14}} v_{_{12}}^*,\  v_{_{34}} = v_{_{14}} v_{_{13}}^*,\ w_{_{12}} = v_{_{12}} + v_{_{12}}^*,\ w_{_{13}} = v_{_{13}} + v_{_{13}}^*,$$ 
$$	w_{_{14}} = v_{_{14}} + v_{_{14}}^*,\ w_{_{23}} = v_{_{23}} + v_{_{23}}^*,\ w_{_{24}} = v_{_{24}} + v_{_{24}}^*,$$ and $w_{_{34}} = v_{_{34}} + v_{_{34}}^*$. We first observe that $w_{_{12}}, w_{_{13}}, w_{_{14}}, w_{_{23}}, w_{_{24}},$ and $w_{_{34}}$ all lie in $R(\mathfrak{J}_{sa})^{-}_{sa} = \mathfrak{J}_{sa}$. Relying on the weak$^*$-closed real algebra $R(\mathfrak{J}_{sa})^{-},$ whose product is denoted by mere juxtaposition, it is not hard to check that the elements $v_{_{ij}}$ ($2\leq i < j \leq 4$) and $w_{_{ij}}$ ($1\leq i < j \leq 4$) are partial isometries in $R(\mathfrak{J}_{sa})^{-}$, with $v_{_{ij}}^* v_{_{ij}} = p_i$ and  $v_{_{ij}} v_{_{ij}}^* = p_j$ for all $2\leq i < j \leq 4$. Furthermore, the following identities hold:
$$w_{_{ij}}^2 = p_i + p_j, \ p_i \circ w_{_{ij}} = \frac12 w_{_{ij}},\ p_k\perp w_{_{ij}}\  (1\leq i < j \leq 4, \ 1\leq k\leq 4, \ k\neq i,j), $$
$$w_{_{12}}\circ w_{_{13}} = \frac12 w_{_{23}}, \ w_{_{12}}\circ w_{_{14}} = \frac12 w_{_{24}},\ w_{_{13}}\circ w_{_{14}} = \frac12 w_{_{34}}, w_{_{12}}\circ w_{_{23}} = \frac12 w_{_{13}},$$
$$w_{_{12}}\circ w_{_{24}} = \frac12 w_{_{14}},\ w_{_{12}}\circ w_{_{34}} = 0,\ w_{_{13}}\circ w_{_{23}} = \frac12 w_{_{12}},\ w_{_{13}}\circ w_{_{24}} = 0,$$ 
$$ w_{_{13}}\circ w_{_{34}} = \frac12 w_{_{14}}, \ w_{_{14}}\circ w_{_{23}} = 0, \ w_{_{14}}\circ w_{_{24}} = \frac12 w_{_{12}}, \ w_{_{14}}\circ w_{_{34}} = \frac12 w_{_{13}},$$
$$ w_{_{23}}\circ w_{_{24}} = \frac12 w_{_{34}},\ w_{_{23}} \circ w_{_{34}} = \frac12 w_{_{24}},\  w_{_{24}}\circ w_{_{34}} = \frac12 w_{_{23}}.$$ 

It follows from the previous identities that the subspace $$\tilde{\mathfrak{B}}:= \hbox{span}\{p_1,p_2,p_3,p_4, w_{i,j}: 1\leq i <j \leq 4\},$$ is a JB$^*$-subalgebra of $\mathfrak{J}$, contains $p_1,p_2,p_3,$ and $p_4$, and is Jordan $^*$-isomorphic to $S_4(\mathbb{C})$ via the Jordan $^*$-isomorphism given by the following assignment: $$p_1 \mapsto \left(\begin{matrix} 1& 0&0&0 \\
 0& 0&0&0 \\  
 0& 0&0&0 \\ 
 0& 0&0&0 \\
\end{matrix}\right),\ p_2 \mapsto \left(\begin{matrix} 0& 0&0&0 \\
0& 1&0&0 \\  
0& 0&0&0 \\ 
0& 0&0&0 \\
\end{matrix}\right),\ p_3 \mapsto \left(\begin{matrix} 0& 0&0&0 \\
0& 0&0&0 \\  
0& 0&1&0 \\ 
0& 0&0&0 \\
\end{matrix}\right),$$ 

$$p_4 \mapsto \left(\begin{matrix} 0& 0&0&0 \\
0& 0&0&0 \\  
0& 0&0&0 \\ 
0& 0&0&1 \\
\end{matrix}\right),\ w_{_{12}} \mapsto \left(\begin{matrix} 0& 1&0&0 \\
1& 0&0&0 \\  
0& 0&0&0 \\ 
0& 0&0&0 \\
\end{matrix}\right),\ w_{_{13}} \mapsto \left(\begin{matrix} 0& 0&1&0 \\
0& 0&0&0 \\  
1& 0&0&0 \\ 
0& 0&0&0 \\
\end{matrix}\right),$$

$$w_{_{14}} \mapsto \left(\begin{matrix} 0& 0&0&1 \\
	0& 0&0&0 \\  
	0& 0&0&0 \\ 
	1& 0&0&0 \\
\end{matrix}\right),\ w_{_{23}} \mapsto \left(\begin{matrix} 0& 0&0&0 \\
0& 0&1&0 \\  
0& 1&0&0 \\ 
0& 0&0&0 \\
\end{matrix}\right), \ w_{_{24}} \mapsto \left(\begin{matrix} 0& 0&0&0 \\
0& 0&0&1 \\  
0& 0&0&0 \\ 
0& 1&0&0 \\
\end{matrix}\right),$$

$$\hbox{ and }\ \  w_{_{34}} \mapsto \left(\begin{matrix} 0& 0&0&0 \\
	0& 0&0&0 \\  
	0& 0&0&1 \\ 
	0& 0&1&0 \\
\end{matrix}\right).$$ The rest of the statement is clear.\smallskip

$(b)$ The idea is to reduce our argument to the previous case. Let us first observe that $v v^* = v (p_1+p_2) v^* = v p_1 v^* + v p_2 v^* = p_3+ p_4,$ and similarly, $v^* v = p_1+ p_2$. Set $v_{_{13}}:= v p_1,$ $v_{_{24}}:= v p_2$ and $v_{_{14}} := v_{_{24}} v_{_{12}}$. Note that $v_{_{13}}$ and $v_{_{14}}$ are partial isometries in $R(\mathfrak{J}_{sa})^{-}$ with $v_{_{13}} v_{_{13}}^* = p_3$, $v_{_{13}}^* v_{_{13}} = p_1$, $v_{_{24}} v_{_{24}}^* = p_4$, and $v_{_{24}}^* v_{_{24}} = p_2$. These identities show that $p_1\approx p_3$ and $p_1\approx p_2 \approx p_4$.\smallskip

By observing that, according to the previous definitions, the identity $v_{_{24}} = v_{_{14}} v_{_{12}}^*$ holds, we are in a position to repeat the arguments in case $(a)$ with $w_{_{ij}} = v_{_{ij}}+v_{_{ij}}^*$ ($1\leq i < j \leq 4$) to conclude that the JB$^*$-subalgebra $$\tilde{\mathfrak{B}}:= \hbox{span}\{p_1,p_2,p_3,p_4, w_{i,j}: 1\leq i <j \leq 4\}$$ contains $p_1,p_2,p_3,$ and $p_4$, and is Jordan $^*$-isomorphic to $S_4(\mathbb{C})$. We finally observe that $v+v^* = v_{_{13}}+ v_{_{24}}+ v_{_{13}}^*+ v_{_{24}}^* = w_{_{13}} + w_{_{24}}\in \tilde{\mathfrak{B}}$, $(v+v^*)^2 = p_1+p_2+p_3+p_4$, $(p_1+p_2)\circ (v+v^*) = \frac12 (v+v^*) = (p_3+p_4)\circ (v+v^*)$, and thus $\hbox{span}\{p_1+p_2, p_3+p_4, v+v^*\}$ is a JB$^*$-subalgebra of $\tilde{\mathfrak{B}}$, which is Jordan $^*$-isomorphic to $S_2(\mathbb{C})$.
\end{proof}

An appropriate version of \Cref{prop a la Christensen} for the Bunce-Wright equivalence is stated next. 

\begin{proposition}\label{prop a la Christensen with the BunceWright equivalence} Let $\mathfrak{J}$ be a  reversible JW$^*$-algebra. Let $p, q_1, q_2, q_3$ be mutually orthogonal projections in $\mathfrak{J}$ such that $p\approx q_i$ for all $i\in \{1,2,3\}.$ Then, for every $a,b$ in $\mathfrak{J}$ with $0\leq a,b\leq \frac12 p$, there are orthogonal projections $r$, $q$ in  $\mathfrak{J}$ satisfying $r,q\leq p + q_1 + q_2+q_3,$ $U_{p} (r) = a,$ and $U_p (q) = b.$
\end{proposition}

\begin{proof} The proof outlined by Christensen in \cite[Proof of Theorem 4.1]{Christensen82} is also valid here up to a small change. By assumptions, there are partial isometries $u_1,$ $u_2,$ and $u_3\in R(\mathfrak{J}_{sa})^{-}$ such that $u_i u_i^* = q_i$ and $u_i^* u_i =p$ for all $i=1,2,3$. Define 
$$x = a^{\frac12} + u_1 a^{\frac12} + u_2 (p-2 a)^{\frac12} \in R(\mathfrak{J}_{sa})^{-}\!\!\!\!\!, \ \& \  y = b^{\frac12} - u_1 b^{\frac12} + u_3 (p-2 b)^{\frac12}\in R(\mathfrak{J}_{sa})^{-}\!\!\!\!\!\!.$$ It is not hard to check that $x$ and $y$ are partial isometries in $R(\mathfrak{J}_{sa})^{-}$, whose range projections $x x^* = r\in R(\mathfrak{J}_{sa})^{-}_{sa} = \mathfrak{J}_{sa}$, and $y y^*= q\in R(\mathfrak{J}_{sa})^{-}_{sa} = \mathfrak{J}_{sa}$ are orthogonal with $r,q\leq p + q_1 + q_2+q_3,$, $U_{p} (r) =a,$ and $U_p (q) = b$.
\end{proof}


In this subsection we need to deal with JBW$^*$-subalgebras generated by a couple of projections (see  \cite[\S 3]{PeSaa2025} for a recent discussion on this topic).\smallskip
 
Having in mind that the Mackey-Gleason theorem fails in the case of $M_2(\mathbb{C})$, and hence it also fails for $S_2(\mathbb{C})$, the next lemma might result a bit surprising at first look. From now on, if $\mathfrak{A}$ and $\mathfrak{J}$ are two JB$^*$-algebras, we shall write $\mathfrak{A}\cong \mathfrak{J}$ if $\mathfrak{A}$ and $\mathfrak{J}$ are Jordan $^*$-isomorphic. 

\begin{lemma}\label{l subalgebra isomorphic to S2} Let $\mu: \PP(\mathfrak{J}) \rightarrow \mathbb{R}$ be a bounded finitely additive signed measure, where $\mathfrak{J}$ is a properly non-modular JBW$^*$-algebra. Let $e$ and $p$ be two projections in $\mathfrak{J}$, where $e$ is properly non-modular, and  let $W^*(\unit, e,p)$ denote the JBW$^*$-subalgebra of $\mathfrak{J}$ generated by $\{\unit, e,p\}$. Then the restriction of the quasi-linear extension of $\mu$ to $W^*(\unit, e,p)$ is linear. 
\end{lemma}

\begin{proof} We begin by observing that $\mathfrak{J}$ is a reversible JW$^*$-algebra (cf. \cite[Theorem 5.3.10]{HOS}). We consider the algebra $R(\mathfrak{J}_{sa})^{-}$, whose product will be denoted by juxtaposition. By \cite[Proposition 3.5]{PeSaa2025} there exists a hyper-Stonean compact space $\Omega$, $m\in \mathbb{N}\cup \{0\},$ and a Jordan $^*$-isomorphism $\Psi: W^*(\unit, e,p) \to   C(\Omega, S_2(\mathbb{C})) \oplus^{\ell_{\infty}} \mathbb{C}^m $ such that $\Psi (e)= \left(\begin{matrix} 1& 0 \\ 0&0 \end{matrix}\right) + p_0$, where $p_0$ is a projection in $\mathbb{C}^m$. We can clearly reduce to the case $m=0$. Let us define $r := \Psi^{-1} \left( \left(\begin{matrix} 0& 0 \\ 0&1 \end{matrix}\right) \right)$ and $u := \Psi^{-1} \left( \left(\begin{matrix} 0& 1 \\ 1&0 \end{matrix}\right) \right)$. Note that the element $v = r u $ is a partial isometry in $R(\mathfrak{J}_{sa})^{-}$ satisfying $vv^* = e$, $v^* v = r$ and $v+v^* =u$. In particular $e\approx r$.\smallskip

As in the proof of \Cref{theo: lin_type_In}, elements in $W^*(\unit, e,p)\cong C(\Omega, S_2(\mathbb{C}))$ can be approximated in norm by elements in JBW$^*$-subalgebras which are finite $\ell_{\infty}$-sums of copies of $W^*(e,r, u)\cong S_2(\mathbb{C}).$ We focus on one of these summands. Since $e$ is properly non-modular, we can find two orthogonal subprojections $p_1$ and $p_2$ such that $p_1\approx p_2\approx e$ and $e = p_1+p_2$ (cf. \Cref{l properties of the BWsim}$(v)$). Let $v_{_{12}}$ denote a partial isometry in $R(\mathfrak{J}_{sa})^{-}$ satisfying $v_{_{12}}^* v_{_{12}} = p_1$ and $v_{_{12}} v_{_{12}}^* =p_2$. The elements $p_3=U_{u} (p_1)$ and $p_4 = U_{u} (p_2)$ are orthogonal projections in $\mathfrak{J}$ satisfying $r = U_{u} (e) = U_{u} (p_1+p_2) = p_3 +p_4$, $U_{u} (p_3)= U_{u} U_{u} (p_1) = U_{u^2} (p_1) = U_{e+r} (p_1)= p_1$, and similarly $U_{u} (p_4) = p_2$. It is not hard to check that the identities $v p_1 v^* = p_3$, $v p_2 v^* = p_4$, $v^* p_3 v = p_1$, and $v^* p_4 v = p_4$ hold. For example, $$v p_1 v^* = r u p_1 u r = U_r U_u (p_1) = U_r (p_3) = p_3,$$
$$ \hbox{and } v^* p_3 v = u r p_3  r u = U_{u} U_{r} (p_3) = U_u (p_3) =  p_1.$$ We are therefore in a position to apply \Cref{p four BW equivalent projections}$(b)$ to deduce the existence of a JBW$^*$-subalgebra $\tilde{\mathfrak{B}}$ of $\mathfrak{J}$ which is Jordan $^*$-isomorphic to $S_4(\mathbb{C}),$ contains the elements $p_1,p_2,p_3,$ $p_4$ and $v+v^*$, and $W^*(\unit, e, p)$ as a JB$^*$-subalgebra. \Cref{theo: lin_type_In} (see also \Cref{r after theorem type In}) prove that the restriction of the quasi-linear extension of $\mu$ to $\tilde{\mathfrak{B}}$ is linear. Consequently, $\overline{\mu}|_{_{W^*(\unit, e, p)}}$ is linear.
\end{proof}

\subsection{Linearity of the quasi-liner extension in the case of properly non-modular \texorpdfstring{JW$^*$}{JW*}-algebras}\ \smallskip

We can now establish a version of the Mackey-Gleason-Bunce-Wright theorem in the case of properly non-modular JBW$^*$-algebras. Let us fix some notation and conventions valid in the remaining part of this section. 

\begin{remark}\label{remark: alpha}
Let $\mathfrak{J}$ be JW$^*$-algebra without type $I_2$ part, $\mu: \PP(\mathfrak{J}) \rightarrow \mathbb{R}$ a bounded finitely additive signed measure. As in \cite[Section 3]{BunceWright1992}, we can (and will) assume that $\mu(\unit) = 0$. Namely, let $\phi$ be a self-adjoint functional in $\mathfrak{J}_{sa}^*$ with $\phi(\unit) = \mu(\unit)$. The mapping $\mu - \phi|_{\PP(\mathfrak{J})} : \PP(\mathfrak{J}) \rightarrow \mathbb{R}$, is a bounded finitely additive signed measure mapping $\unit$ to zero. Obviously, $\mu$ admits an extension to a bounded linear function on $\mathfrak{J}$ if and only if $\mu - \phi|_{\PP(\mathfrak{J})} $ does.\smallskip 

We can also assume, without loss of generality, that $$ \sup\Big\{|\overline{\mu}(x)| \, : \, x^* = x, \, \|x\| \leq 1 \Big\} = 1,$$ an assumption, which combined with $\mu(\unit) =0$ and  \Cref{prop_quasilin}$(c),$ gives 
$$ \alpha_{_\mu}(\unit) = \sup \Big\{ \mu(p) \, : \, p \in \PP(\mathfrak{J})\Big\} = -\inf \Big\{ \mu(p) \, : \, p \in \PP(\mathfrak{J})\Big\} = \frac{1}{2}$$
\end{remark}

We state next some technical results, which are appropriate Jordan versions of results in \cite[Section 3]{BunceWright1992}.

\begin{proposition}\label{prop: p_decomp_acot}
Let $\mu: \PP(\mathfrak{J}) \rightarrow \mathbb{R}$ be a bounded finitely additive signed measure satisfying the assumptions in \Cref{remark: alpha}, where $\mathfrak{J}$ is a properly non-modular JBW$^*$-algebra. Let $ 0 < \varepsilon < \frac{1}{2}$ and let $q$ be a properly infinite projection in $ \mathfrak{J}$ satisfying $\mu(q)> \frac{1}{2} - \varepsilon$. 
Then the inequality
$$ \left|\mu(p) - \overline{\mu}(U_q(p)) -\overline{\mu}(U_{\unit-q}(p)) \right|< 4 \varepsilon^{\frac{1}{2}},$$ holds for every projection $p \in \mathfrak{J}$, where $\overline{\mu}$ denotes the quasi-linear extension of $\mu$.
\end{proposition}

\begin{proof} As in the proof of \Cref{prop a la Christensen}, we can assume from the hypothesis that $\mathfrak{J}$ is a JBW$^*$-subalgebra of some von Neumann algebra $A$ sharing the same unit element, where the associative product on $A$ will be denoted by mere juxtaposition. Clearly, each symmetry and each projection in $\mathfrak{J}$ is a symmetry and a projection in $A$, respectively.\smallskip

Having in mind that $0\leq U_{\unit-q}(p) \leq \unit$, \Cref{prop_quasilin}$(c)$ implies that \begin{equation}\label{eq 1 0506} \left| \overline{\mu} (U_{\unit-q}(p)) \right|\leq 2 \alpha_{_\mu} (\unit) = 1.
\end{equation}

Let us define the elements $x_{\pm} = (1- \varepsilon)^{\frac{1}{2}}q \pm \varepsilon^{\frac{1}{2}}(\unit-q)p \in A$. Then $$x_{\pm}^*x_{\pm} = (1-\varepsilon)q + \varepsilon U_p(\unit-q) \in \mathfrak{J},$$ and moreover $x_{\pm}^*x_{\pm} \leq \unit$. This implies that 
    $$ 0 \leq x_{\pm}x_{\pm}^* = (1-\varepsilon) q + \varepsilon U_{\unit-q}(p) \pm 2 \varepsilon^{\frac{1}{2}}(1-\varepsilon)^{\frac{1}{2}}\{q,p,\unit-q\} \leq \unit,$$
    which also assures that $x_{\pm}x_{\pm}^* \in \mathfrak{J}$. Having in mind that $q$ is properly non-modular, it follows from \Cref{l subalgebra isomorphic to S2} that the restriction of $\overline{\mu}$ to $W^*(\unit, p, q)$ is linear, thus by the assumptions in  \Cref{remark: alpha} together with \Cref{prop_quasilin}$(e)$ we get
    $$\frac{1}{2}= \alpha_{_\mu} (\unit) \geq (1-\varepsilon)\overline{\mu}(q)  + \varepsilon \overline{\mu}(U_{\unit-q}(p)) \pm 2 \varepsilon^{\frac{1}{2}}(1-\varepsilon)^{\frac{1}{2}}\overline{\mu}(\{q,p,\unit-q\})$$
    $$\hspace{0.3cm}\geq\hbox{ (by \eqref{eq 1 0506}) } \geq (1-\varepsilon)\left(\frac{1}{2}-\varepsilon\right) -\varepsilon \pm 2\varepsilon^{\frac{1}{2}}(1-\varepsilon)^{\frac{1}{2}}\overline{\mu}(\{q,p,\unit -q\}),$$ and consenquently $$ 2 \left|\overline{\mu}(\{q,p,\unit -q\}) \right|\leq \frac{(5/2-\varepsilon)\varepsilon}{\varepsilon^{\frac{1}{2}}(1-\varepsilon)^{\frac{1}{2}}} = \varepsilon^{\frac{1}{2}} \frac{(5/2-\varepsilon)}{(1-\varepsilon)^{\frac{1}{2}}}  < 4\varepsilon^{\frac{1}{2}}. $$
    Finally, we compute
    $$| \mu(p) - \overline{\mu}(U_q(p)) - \overline{\mu}(U_{\unit -q}(p)) |  = | \overline{\mu}(p) - \overline{\mu}(U_q(p)) - \overline{\mu}(U_{\unit -q}(p)|\hspace{1.293cm} $$ $$\hspace{2.8cm} = | 2 \overline{\mu}(\{q, p, \unit -q\}) |  < 4\varepsilon^{\frac{1}{2}}. $$\end{proof}

We continue our discussion with a couple of technical lemmas. 

\begin{lemma}\label{lemma: bound_V_mu}
    Let $\mathfrak{J}$ be a properly non-modular JW$^*$-algebra, and let $\mu: \PP(\mathfrak{J}) \rightarrow \mathbb{R}$ be a bounded finitely additive signed measure satisfying the assumptions in \Cref{remark: alpha}. Suppose $e \in \PP(\mathfrak{J})$ is properly non-modular. Then for every  $ \varepsilon > 0$, there is a projection $p \leq e$ satisfying $p \approx e \approx e-p$ in $\mathfrak{J}_{p} = U_p (\mathfrak{J})$ and 
    $$V_{_\mu}(p) = \sup \{ | \mu(q)| \, : \,  q \in \PP(\mathfrak{J}), \ q \leq p\} \leq \varepsilon.$$
\end{lemma}

\begin{proof} Note that $\mathfrak{J}$ is reversible (cf. \cite[Theorem 5.3.10]{HOS}).  Let us take $n \in \mathbb{N}$ such that $2^n \varepsilon > 1$. Since $e$ is properly non-modular \Cref{l properties of the BWsim}$(v)$ assures the existence of a set $\{ e_i, i = 1, \dots, 2^n \}$ of mutually orthogonal projections  such that $e = e_1 + \dots + e_{_{2^n}}$ and $e_j \approx e_i \approx e$ in $\mathfrak{J}$ for every $ i ,j \in \{1, \dots , 2^n\}$. If $V_{_\mu}(e_j) > \varepsilon$ for every $j \in \{1, \dots , 2^n\}$, by definition, there exist projections $r_j \leq e_j$ in $\mathfrak{J}$ with $| \mu(r_j)| > \varepsilon$. Define  $$\Gamma^+ = \{j \in \{1,\ldots,2^n\}:  \mu(r_j) >0 \}, \hbox{ and } \Gamma^- = \{j \in \{1,\ldots,2^n\}:  \mu(r_j) <0\}.$$ We shall distinguish two possible cases: \smallskip
    
\noindent $(a)$ Suppose first that  $ \#\Gamma^{+} \geq \#\Gamma^{-}$.  Define $f^+ = \sum_{j\in \Gamma^+} r_j$. It is clear that $f^+$ is a projection in $\PP(\mathfrak{J})$ with  $\mu(f^+) > 0$, and moreover $ \mu(f^+) = \sum_{j\in \Gamma^+} \mu(r_j) \geq \frac{2^n}{2}\varepsilon = 2^{n-1}\varepsilon > \frac{1}{2},$ which is impossible since, by our assumptions, $\alpha_{_\mu} (\unit) = \frac{1}{2}$ (cf.  \Cref{remark: alpha}).\smallskip
    
\noindent $(b)$ Similarly, if $ \#\Gamma^{-} \geq\#\Gamma^{+}$, we set $f^{-} = \sum_{j\in \Gamma^-} r_j\in \PP(\mathfrak{J}),$ and we get  $ - | \mu (f^-)| = \mu(f^-) = \sum_{j\in \Gamma^-} \mu(r_j) \leq -\frac{2^n}{2}\varepsilon = -2^{n-1}\varepsilon < -\frac{1}{2},$
    which also contradicts that $\alpha_{_\mu}(\unit) = \frac{1}{2}$ because $\mu (\unit-f^-) = -\mu (f^-) >\frac12$.
\end{proof}

\begin{lemma} {\label{lemma: sup_eq_half}}
    Let  $\mu: \PP(\mathfrak{J}) \rightarrow \mathbb{R}$ be a bounded finitely additive signed measure satisfying the assumptions in \Cref{remark: alpha}, where $\mathfrak{J}$ is a properly non-modular JW$^*$-algebra. Then 
    $$ \sup \Big\{ \mu(e) \, : \, e \in \PP(\mathfrak{J}), e \approx \unit \approx \unit -e \Big\} = \frac{1}{2}. $$
\end{lemma}

\begin{proof} Let us take $0 < \varepsilon < \frac{1}{8}$ and $e \in \PP(\mathfrak{J})$ such that $\mu(e) > \frac{1}{2} - \varepsilon$. By \Cref{remark: decomp_central-proj} there exist a central projection $z \in \PP(\mathfrak{J})$ satisfying 
$ e = z \circ e + (\unit-z)\circ e$, where $z\circ e$ is modular and $(\unit-z) \circ e$ is properly non-modular.\smallskip

Observe that $p = z\circ (\unit-e)$ and $e + p = z + (\unit -z)\circ e $ are properly non-modular projections. Otherwise, if for example $z\circ (\unit-e)$ is not properly non-modular, there exists a central projection $\tilde{z}$ in $\mathfrak{J}$ such that $\tilde{z}\circ (z\circ (\unit-e)) = (\tilde{z}\circ z)\circ (\unit-e)$ $= \tilde{z}\circ z - (\tilde{z}\circ z)\circ e$ is a non-zero modular projection.\smallskip 

Since $(\tilde{z}\circ z)\circ e$ is a modular projection and the sum of two orthogonal modular projections is a modular projection \cite[Theorem 7.6.4]{HOS}, the projection $\tilde{z}\circ z$ must be modular too, which contradicts that $\mathfrak{J}$ is properly non-modular. If on the other hand, $e + p = z + (\unit -z)\circ e $ is not properly non-modular, there exists a central projection $\tilde{z}$ in $\mathfrak{J}$ such that $\tilde{z}\circ (e+p)=\tilde{z}\circ ( z + (\unit -z)\circ e ) = \tilde{z}\circ z + \tilde{z}\circ ((\unit -z)\circ e ) $ is a non-zero modular projection. However, $(\unit -z)\circ e$ being properly non-modular implies that $\tilde{z}\circ ((\unit -z)\circ e )=0$, and thus $\tilde{z}\circ z = \tilde{z}\circ (e+p)$ must be a non-zero modular projection, which is also impossible.\smallskip  

Now, having in mind that that $p$ and $e+p$ are properly non-modular projections, by \Cref{lemma: bound_V_mu}, there exists a projection $f\leq p$ such that $f \approx p \approx p-f$ and $V_{_\mu}(f) < \varepsilon$. Then $q := e +f \approx e +p$ (cf. \Cref{l properties of the BWsim}$(ii)$) is properly non-modular (see \Cref{l properties of the BWsim}$(v)$) and $\mu(q) > \frac{1}{2} - 2 \varepsilon$. Since $q$ is properly non-modular, we apply \Cref{lemma: bound_V_mu} to find a projection $r \leq q$ with $r \approx q \approx q-r$ and $V_{\mu}(q-r) < \varepsilon$. Then $\unit - r = (\unit - q) + (q -r) \approx (\unit -q) + q = \unit$ (cf. \Cref{l properties of the BWsim}$(ii)$). Therefore 
$$ \mu(\unit - r) = \mu(\unit - q) + \mu(q -r) = - \mu(q) + \mu(q-r) < -\frac{1}{2} + 3 \varepsilon.$$
Finally, since $\unit- r$ is a properly non-modular projection (cf. \Cref{c properly modular is preserved by BW equivalence}), we apply \Cref{lemma: bound_V_mu} once again to find a projection $h \leq \unit -r$ such that $h \approx \unit - r \approx \unit -r-h$ and $V_{\mu}(\unit -r-h) < \varepsilon$. It follows from the above that $h \approx \unit-r \approx \unit$, $\unit - h =  (\unit -r-h) + r\approx \unit-r +r = \unit,$ and $\mu(\unit - h) > \frac{1}{2} - 4 \varepsilon$, which completes the proof.
\end{proof}

The next proposition completes the technical tools required for the main result in this section.

\begin{proposition}\label{prop: lin_aprox}
Let $\mu: \PP(\mathfrak{J}) \rightarrow \mathbb{R}$ be a bounded finitely additive signed measure satisfying the assumptions in \Cref{remark: alpha}, where $\mathfrak{J}$ is a properly non-modular JW$^*$-algebra. Then for each $\delta>0$ there exists a properly non-modular projection $e$ in $\PP(\mathfrak{J})$ satisfying the following statements: 
\begin{enumerate}[$(a)$]
\item $\displaystyle\Big|\mu(p) -\overline{\mu}(U_e(p)) - \overline{\mu}\left(U_{\unit - e}(p)\right) \Big| < \delta$, for all $p \in \PP(\mathfrak{J})$.
\item $\displaystyle \left| \sum_{i=1}^n\overline{\mu}(x_i) - \overline{\mu} \left(\sum_{i=1}^n x_i\right) \right|< (n-1)\delta$,  for every $n\in \mathbb{N}$, $n\geq 2$ and every   $x_1,\ldots, x_n$ in $\mathfrak{J}$ such that $\displaystyle 0 \leq x_i\leq \sum_{i=1}^n x_i \leq e $.
\item $\displaystyle \left| \sum_{i=1}^n\overline{\mu}(x_i) - \overline{\mu} \left(\sum_{i=1}^n x_i\right) \right|< (n-1)\delta,$ for every $n\in \mathbb{N}$, $n\geq 2$ and every   $x_1,\ldots, x_n$ in $\mathfrak{J}$ such that $\displaystyle 0 \leq x_i\leq \sum_{i=1}^n x_i \leq \unit-e $.
\end{enumerate}
\end{proposition}

\begin{proof}
\noindent$(a)$ Given $\delta>0$, let us take $0 < \varepsilon < \frac{1}{4}$ with $24 \varepsilon^{\frac12}<\delta$. By \Cref{lemma: sup_eq_half}, there exists a properly non-modular projection $e$ in $\PP(\mathfrak{J})$ such that $e \approx \unit \approx \unit -e$ and $\mu(e) > \frac{1}{2} - \varepsilon$. Under these conditions \Cref{prop: p_decomp_acot} implies that 
\begin{equation}\label{eq (a) improved in the proof of p4.9} \left|\mu(p) -\overline{\mu}\left(U_e(p)\right) - \overline{\mu}\left(U_{\unit - e}(p)\right)\right| < 4 \varepsilon^{\frac{1}{2}} <\delta,
\end{equation} for all $p\in \PP(\mathfrak{J})$, which completes the first assertion. \\

\noindent$(b),(c)$ We prove the statements in the case $n = 2$. The general case, which is left to the reader, follows by a straightforward induction on $n \geq 2$.\smallskip

Since $\mathfrak{J}$ is properly non-modular, every central projection in $\mathfrak{J}$ is properly non-modular. We go back to the projection $e$ obtained in the first paragraph of this proof. Since $\unit - e\approx \unit\approx e$, the projections $e$ and $\unit-e$ are both properly non-modular (cf. \Cref{c properly modular is preserved by BW equivalence}). Thus, by \Cref{lemma: bound_V_mu}, there is a projection $r \leq \unit - e$, satisfying $ r, \unit-e-r \approx \unit - e \approx \unit$ and $V_{_\mu}(r) < \varepsilon$. \Cref{c properly modular is preserved by BW equivalence} implies that $r$ is also properly non-modular. By \Cref{l properties of the BWsim}$(v)$ there exist orthogonal projections $r_1,r_2,r_3\in \PP(\mathfrak{J})$ satisfying $r = r_1+r_2+r_3$ and $r_i \approx r$ for all $i = 1,2,3$. In particular, $r_i \approx \unit-e \approx \unit$.\smallskip 

Given positive elements $x,y$ in $\mathfrak{J}$ with $x,y \leq e$, we apply \Cref{prop a la Christensen with the BunceWright equivalence} to $e,r_1, r_2, r_3$, $\frac{1}{2}x,$ and $\frac{1}{2}y$, to obtain two orthogonal projections $\tilde{p},\tilde{q}$ in $\mathfrak{J}$ satisfying $\tilde{p},\tilde{q}\leq  h := e +r\approx e + \unit -e = \unit \approx r,$ $\frac{1}{2} x = U_e(\tilde{p}),$ and $\frac{1}{2} y = U_e(\tilde{q})$. Since $x \in W^*(\unit, e, \tilde{p},$ and by \Cref{l subalgebra isomorphic to S2}, $\overline{\mu}$ is linear when restricted to $W^*(\unit, e, \tilde{p})$, we deduce that
$$| 2 \mu(\tilde{p}) - \overline{\mu}(x)| = | \overline{\mu}(2 \tilde{p} -x)| = \left|\overline{\mu}\left(2 \tilde{p} - 2 U_e(\tilde{p})\right)\right| = 2 \left|\overline{\mu}\left(\tilde{p} - U_e(\tilde{p})\right)\right|.$$
Observe now that the projections $ e, r $ and $\tilde{p} \leq h$ lie in the JBW$^*$-subalgebra $\mathfrak{J}_h = U_h(\mathfrak{J})$ (see \cite[Proposition 2.9]{AlfsenShultz2003}), therefore $e$ lies in the $^*$-subalgebra $ W^*(h, r, \tilde{p})$ of $\mathfrak{J}_h$ generated by $r$, $\tilde{p}$ and the unit of
$\mathfrak{J}_h$. On the other hand, $W^*(\unit,h, r, \tilde{p}) \cong \mathbb{C} (\unit-h)\oplus^{\infty} W^*(h, r, \tilde{p})$, where $r\approx \unit$ is properly non-modular, $W^*(h, r, \tilde{p})\subseteq \mathfrak{J}_{h}$, and the latter is properly non-modular. \Cref{prop_quasilin}, \Cref{l subalgebra isomorphic to S2}, and the properties of $\mu$ assure that $\overline{\mu}$ is linear on $W^*(\unit, h, r, \tilde{p})$. We therefore have 
$$ \left| \overline{\mu}(\tilde{p}  - U_e(\tilde{p})) \right|  = \left| \overline{\mu}(U_{h} (\tilde{p})  - U_e(\tilde{p})) \right| = \left| \overline{\mu}(U_{e+r} (\tilde{p})  - U_e(\tilde{p})) \right| \hspace{1cm}$$ 
 $$\hspace{2.5cm} = \left| \overline{\mu}\left(U_r(\tilde{p}) + 2\{r,\tilde{p},e\}\right) \right| = \left| \overline{\mu}\left(U_r(\tilde{p})\right) + 2\overline{\mu}\left( \{r,\tilde{p},e\}\right) \right| $$ $$\leq \left|\overline{\mu}\left(U_r(\tilde{p})\right) \right| + \left| \overline{\mu}\left(2 \{r,\tilde{p},e\}\right)\right|. \hspace{1.1cm}$$
 
  By applying \Cref{prop_quasilin}$(c)$ and $(d)$ to $\mu|_{\PP(\mathfrak{J}_r)}$ we derive that $$|\mu(U_r(\tilde{p}))| \leq 2 V_{_{\mu|_{\PP(\mathfrak{J}_r)}}}(r) = 2 V_{_\mu}(r)  < 2 \varepsilon.$$ We consequently have 
$$|\overline{\mu}(\tilde{p}) - \overline{\mu}(x)| < 4 \varepsilon  + 2 \left| \overline{\mu}\left(2 \{r,\tilde{p},e\}\right)\right|.$$ Note that $\{r,\tilde{p},e\} = \{\unit-e,\tilde{p},e\}$. A new application of \Cref{l subalgebra isomorphic to S2} assures that $\overline{\mu}|_{W^*(\unit, e,\tilde{p})}$ is linear. It then follows from $(a)$ (actually from \eqref{eq (a) improved in the proof of p4.9}) that
$$\left| \overline{\mu}\left(2 \{r,\tilde{p},e\}\right)\right| = \left| \overline{\mu}\left(2 \{\unit-e,\tilde{p},e\}\right)\right| = \left| \overline{\mu}\left(U_{\unit} (\tilde{p}) -U_e (\tilde{p}) -U_{\unit-e} (\tilde{p})  \right) \right| $$
$$\hspace{1.3cm}= \left| \overline{\mu}\left(U_{\unit} (\tilde{p}) \right) - \overline{\mu}\left(U_e (\tilde{p})\right) - \overline{\mu}\left(U_{\unit-e} (\tilde{p})  \right) \right| < 4 \varepsilon^{\frac12},$$ which in turn gives $$ \left|\mu(p) - \overline{\mu}(x)\right| < 4 \varepsilon + 8 \varepsilon^{\frac{1}{2}} < 12\varepsilon^{\frac{1}{2}}.$$

We can similarly obtain $$ |\mu(\tilde{q}) - \overline{\mu}(y)|< 12 \varepsilon^{\frac{1}{2}}, \hbox{ and } |\mu(\tilde{p}+\tilde{q}) - \overline{\mu}(x+ y)| < 12\varepsilon^{\frac{1}{2}}. $$
By combining the previous conclusions we get
$$| \overline{\mu}(x+y) - \overline{\mu}(x) - \overline{\mu}(y)| < 24 \varepsilon^{\frac{1}{2}}<\delta,$$
for every $0 \leq x,y \leq  e$.  Since $(-\mu)(\unit - e) = \mu(e) > \frac{1}{2} - \varepsilon$, we get the same inequality for every $0\leq x,y \leq \unit -e$.
\end{proof}

\begin{theorem}\label{theo: mu_lin_prop_inf} Let $\mu : \PP(\mathfrak{J}) \rightarrow \mathbb{R}$ be a bounded finitely additive signed measure, where $\mathfrak{J}$ is a properly non-modular JW$^*$-algebra. Then $\mu$ extends to a linear functional on $\mathfrak{J}$. 
\end{theorem}

\begin{proof} We can clearly assume that $\mu$ satisfies the assumptions in \Cref{remark: alpha}.\smallskip

According to the notation fixed in this note, $\overline{\mu}$ will stand for the quasi-linear extension of $\mu$. Let us take $x, y \in \mathfrak{J}$ with $0 \leq x,y \leq \frac{1}{2}\unit$, which implies that $0 \leq x,y \leq x+y \leq \unit$. It is well-known that the JBW$^*$-algebras $W^*(\unit, x), W^*(\unit,y)$ and $ W^*(\unit, x+y)$ are associative, and hence commutative von Neumann algebras \cite[Proposition 2.11]{AlfsenShultz2003}. \smallskip

It is also part of the folklore in von Neumann algebra theory that for each $\varepsilon>0$ we can find finite families of projections $\{p_n\}_{1\leq n\leq m_1} \subset W^*(\unit, x)$, $\{q_n\}_{1\leq n\leq m_2} \subset W^*(\unit, y)$ and $\{r_n\}_{1\leq n\leq m_3} \subset W^*(\unit, x+y),$ and finite collections of positive numbers $(\lambda_n)_{1\leq n\leq m_1},$ $(\mu_n)_{1\leq n\leq m_2},$ and $(\delta_n)_{1\leq n\leq m_3}$ such that $\displaystyle\sum_{n=1}^{m_1} \lambda_n = \|x\|\leq \frac12,$ $\displaystyle \displaystyle\sum_{n=1}^{m_2} \mu_n = \|y\|\leq \frac12,$  $\displaystyle\sum_{n=1}^{m_3} \delta_n =\|x+y\|\leq 1,$ $\displaystyle x - \sum_{n=1}^{m_1} \lambda_n p_n\geq 0,$ $   \displaystyle y - \sum_{n=1}^{m_2} \mu_n q_n\geq 0,$   $\displaystyle x + y  - \sum_{n=1}^{m_3} \delta_n r_n\geq 0,$ and	$$\left \| x - \sum_{n=1}^{m_1} \lambda_n p_n \right \|, \  \left \|y - \sum_{n=1}^{m_2} \mu_n q_n  \right \|, \ \left \| x + y  - \sum_{n=1}^{m_3} \delta_n r_n \right \|  < \varepsilon,$$ (see, for example, \cite[Proposition 4.2.3]{HOS}).\smallskip
	
To simplify the notation, set $\displaystyle 0\leq a := x - \sum_{n=1}^{m_1} \lambda_n p_n \in W^*(\unit , x)$. By \Cref{prop_quasilin}, the restriction of $\overline{\mu}$ to $W^*(\unit, x)$, is linear and hence  
    $$\overline{\mu}(x) = \overline{\mu} \left( a + \sum_{n=1}^{m_1} \lambda_n p_n \right) = \overline{\mu} (a) + \overline{\mu}\left(\sum_{n=1}^{m_1} \lambda_n p_n \right) =\overline{\mu}(a) + \sum_{n=1}^{m_1} \lambda_n \overline{\mu}(p_n).$$

Set $m_0 =\max\{m_1,m_2,m_3\}$. Let $e$ be the properly non-modular projection in $\PP(\mathfrak{J})$ given by \Cref{prop: lin_aprox} for $\delta = \frac{\varepsilon}{m_0}$. \Cref{l subalgebra isomorphic to S2} proves that $\overline{\mu}|_{W^*(\unit, p_j,e)}$ is linear for each $j \in \mathbb{N}$. We therefore deduce that  
    $$ \overline{\mu}(p_j) = \overline{\mu}\Big(p_j \pm U_e(p_j) \pm U_{\unit -e}(p_j)\Big) \hspace{4.3cm}$$ 
    $$\hspace{1.3cm}= \overline{\mu} \Big(p_j - U_e(p_j) - U_{\unit -e}(p_j)\Big) + \overline{\mu}\Big(U_e(p_j)\Big) + \overline{\mu}\Big(U_{\unit -e}(p_j)\Big), $$
  which implies that 
  \begin{equation}\label{eq: 1}
      \left\{\begin{aligned}
          \overline{\mu}(x) = \overline{\mu}(a) &+ \sum_{n=1}^{m_1} \lambda_n \overline{\mu}\Big(p_n - U_e(p_n) - U_{\unit -e}(p_n)\Big) \\
          & + \sum_{n=1}^{m_1} \lambda_n \overline{\mu}\Big(U_e(p_n)\Big) + \sum_{n=1}^{m_1} \lambda_n \overline{\mu}\Big(U_{\unit - e}(p_n)\Big).
      \end{aligned}\right.
  \end{equation}

\Cref{prop: lin_aprox} proves that the following inequalities hold:  
\begin{equation}\label{eq first application of 4.11 for x+y} \left\{
\begin{aligned}
		&| \overline{\mu}(U_e(x+y)) - \overline{\mu}(U_e(x)) - \overline{\mu}(U_e(y)) | <\delta \leq \varepsilon, \\
	&|\overline{\mu}(U_{\unit -e}(x+y)) - \overline{\mu}(U_{\unit -e}(x)) - \overline{\mu}(U_{\unit -e}(y))| <\delta \leq \varepsilon,
\end{aligned}\right.
\end{equation}
\begin{equation}\label{eq: 2207}
\big|\overline{\mu}\Big(p_n - U_e(p_n) - U_{\unit -e}(p_n)\Big)\Big| < \delta ,
\end{equation}
\begin{equation}\label{eq: 2}
    \left| \overline{\mu}\Big(U_e(x)\Big) - \overline{\mu}\Big(U_e(a)\Big) - \sum_{n=1}^{m_1}  \lambda_n \overline{\mu}\Big(U_e(p_n)\Big)\right| < m_1 \delta <\varepsilon ,
\end{equation}
and 
\begin{equation}\label{eq: 3}
    \left| \overline{\mu}\Big( U_{\unit-e}(x)\Big) - \overline{\mu}\Big(U_{\unit -e}(a)\Big) - \sum_{n=1}^{m_1} \lambda_n \overline{\mu}\Big(U_{\unit -e}(p_n)\Big) \right| < m_1 \cdot \delta<\varepsilon. 
\end{equation}

By combining \eqref{eq: 1}, \eqref{eq: 2207}, \eqref{eq: 2} and \eqref{eq: 3} we deduce that 
\begin{equation}\label{eq inequ for x in 2207} \begin{aligned}
		\Big|\overline{\mu}(x) - \overline{\mu}\Big(U_e(x)\Big) - \overline{\mu}\Big(U_{\unit - e}(x)\Big) \Big| & \leq \Big| \overline{\mu}(a) - \overline{\mu}(U_e(a)) - \overline{\mu}(U_{\unit -e}(a))\Big| \\
		 + 2m_1 \delta +&  \left | \sum_{n=1}^{m_1} \lambda_n \, \overline{\mu}\Big(p_n - U_e(p_n) - U_{\unit -e}(p_n)\Big)\right |\\
		&< 3 \| a\| + 2\varepsilon + \delta < 6 \varepsilon.
	\end{aligned}
\end{equation} 

Similar arguments show that the previous inequality also holds when $x$ is replaced by $y$ and $x+y$, respectively, that is, \begin{equation}\label{eq inequ for y in 2207} \begin{aligned}
		\Big|\overline{\mu}(y) - \overline{\mu}\Big(U_e(y)\Big) - \overline{\mu}\Big(U_{\unit - e}(y)\Big) \Big|< 6 \varepsilon,
	\end{aligned}
\end{equation} and 
\begin{equation}\label{eq inequ for x+y in 2207} \begin{aligned}
		\Big|\overline{\mu}(x+y) - \overline{\mu}\Big(U_e(x+y)\Big) - \overline{\mu}\Big(U_{\unit - e}(x+y)\Big) \Big|< 6 \varepsilon.
	\end{aligned}
\end{equation}

Finally, by combining \eqref{eq inequ for x in 2207}, \eqref{eq inequ for y in 2207}, \eqref{eq inequ for x+y in 2207}, and \eqref{eq first application of 4.11 for x+y} we derive that 
$$|\overline{\mu}(x+y) - \overline{\mu}(x) - \overline{\mu}(y)| < | \overline{\mu}(U_e(x+y)) - \overline{\mu}(U_e(x)) - \overline{\mu}(U_e(y)) | \hspace{2.59cm} $$
    $$\hspace{2.19cm} + |\overline{\mu}(U_{\unit -e}(x+y)) - \overline{\mu}(U_{\unit -e}(x)) - \overline{\mu}(U_{\unit -e}(y))| $$
     $$\hspace{0.61cm} + \Big|\overline{\mu}(x) - \overline{\mu}\Big(U_e(x)\Big) - \overline{\mu}\Big(U_{\unit - e}(x)\Big) \Big| $$
     $$\hspace{0.73cm} + \Big|\overline{\mu}(y) - \overline{\mu}\Big(U_e(y)\Big) - \overline{\mu}\Big(U_{\unit - e}(y)\Big) \Big| $$
     $$\hspace{2.590083cm} + \Big|\overline{\mu}(x+y) - \overline{\mu}\Big(U_e(x+y)\Big) - \overline{\mu}\Big(U_{\unit - e}(x+y)\Big) \Big|
     < 20\varepsilon.$$ It follows from the arbitrariness of $\varepsilon>0$ that $\overline{\mu}$ is actually linear. 
\end{proof}

\section{Intermediate value property for centre-valued traces on modular \texorpdfstring{JW$^*$}{JW-star}-algebras}\label{sec: modular JBW-algebras}

This section is devoted to establish a Mackey-Gleason-Bunce-Wright theorem for quantum measures on the lattice of projections of an arbitrary  modular JW$^*$-algebra, by proving that the quasi-linear extension of any such a measure ${\mu},$ as defined in \Cref{sec: prelim}, is linear. We shall try to adapt the methods introduced by F.W. Yeadon \cite{Yeadon1984} and L.J. Bunce and J.D.M. Wright \cite{BunceWright1994} from the setting of finite von Neumann algebras to modular JBW$^*$-algebras. However, the process is not totally trivial and will required the development of new non-trivial tools.  One of the main novelties required in our arguments --namely the \emph{intermediate value property for centre-valued traces on JW$^*$-algebras of type $II_1$}-- names the whole section. In the case of a type $II_1$ von Neumann the intermediate value property forms part of the tools employed by F.W. Yeadon \cite[Lemma 1]{Yeadon1984} and L.J. Bunce and J.D.M. Wright \cite[Lemma 2.3]{BunceWright1994}. The reader interested on a detailed proof of the result for type $II_1$ von Neumann algebras can consult the recent lecture notes by de Santiago and Nelson in \cite[Theorem 5.4.13.]{deSantiagoNelson2024}.\smallskip
 
Before stating the result in the case of JW$^*$-algebras, recall that every modular JW$^*$-algebra admits a unique faithful normal centre-valued trace $\tau$  (see \cite[Section 18]{Topping65} or subsection~\ref{subsec: traces}). 

\begin{proposition}[Intermediate value property for centre-valued traces]\label{prop: trace_surject} Let $\mathfrak{J}$ be a JW$^*$-algebra of type $II_1$, and let $\tau: \mathfrak{J} \rightarrow Z(\mathfrak{J})$ denote the normal centre-valued unital faithful trace. Let $p \in \PP(\mathfrak{J})$. Then, for each $w \in Z(\mathfrak{J})$ with $0\leq w \leq \tau(p),$ there exist $q \in \PP(\mathfrak{J})$ satisfying $q \leq p$ and $\tau(q) = w$.
\end{proposition}

\begin{proof} We can assume, without loss of generality, that $0\neq p$ and $0 < w < \tau(p)$. Fix a natural $n$. By \cite[Theorem 17]{Topping65} there exists a family of mutually orthogonal projections $\{p_1, \dots, p_{2^n} \}$ in $\mathfrak{J}$ such that $p_i \sim_1 p_j$ for all $i,j \in \{1, \dots, 2^n\}$ and $p = p_1 + p_2 \dots+ p_{2^n}$. Thus, by the linearity of the trace we have $\displaystyle \sum_{i=1}^{2^n} \tau(p_i) = \tau(p) $. Furthermore, since $p_i \sim_1 p_j$ we have $\tau(p_i) = \tau(p_j),$ for all $i,j$ (cf. \Cref{lemma: centre_trace}). Thus, $\tau (p_j) = \frac{1}{2^n}\tau(p),$ for all $j \in \{1, \dots, 2^n\}$, where $n$ is a fixed but arbitrary natural number.\smallskip
	
Consider now the set 
$$ \PP(\mathfrak{J})^-_{\tau, z, p} :=\{ \tilde{q} \in \PP(\mathfrak{J}) \, | \, 0< \tilde{q} \leq p ,\, \tau(\tilde{q}) \leq z \} .$$
We claim that $\PP(\mathfrak{J})^+_{\tau, z, p}$ is not empty. Namely, by structure theory $Z(\mathfrak{J}) \cong C(K)$, where $K$ is a hyper-Stonean compact Hausdorff space. Since $0 < w$, there exists a non-zero central projection $z_0\in \mathfrak{J}$ satisfying $w(t) >0$ whenever $t\in K$ with $z_0(t) \neq 0$ (which, in turn, implies that $w(t)\geq \rho_0\in \mathbb{R}^+$ for all $t\in \{t\in K: z_0(t) =1\}$). Observe that $z_0\circ w \leq z_0\circ  \tau(p) = \tau(z_0\circ p).$	Take a natural number $n$ satisfying $\frac{1}{2^n}<\rho_0$. By the conclusions in the first part of this proof, there exists a non-zero projection $\tilde{q}$ such that $\tilde{q} \leq p$ and $\tau(\tilde{q}) = \frac{1}{2^n}\tau(p)$. It follows from this that $z_0\circ  w \geq \rho_0 z_0 > \frac{1}{2^n} z_0\circ \tau(p) = z_0 \circ \tau(\tilde{q}) = \tau (z_0\circ \tilde{q}),$ and hence the element $z_0\circ \tilde{q}$ lies in $\PP(\mathfrak{J})^-_{\tau, z, p} $.\smallskip 
	
Consider next the set  
    $$ \mathcal{C}_z := \left \{ (q_i)_{i \in \Gamma} \hbox{ totally ordered }:  q_i\in \PP(\mathfrak{J})^-_{\tau, z, p}  \, \, \forall i \in \Gamma\right \},$$ which is inductive with respect to the partial order given by inclusion. Hence, by Zorn's Lemma, there exists a maximal element $\{q_i\}_{i \in \Gamma_0}\in \mathcal{C}_z$. Setting 
    $ q:= \bigvee_{i \in \Gamma_0} q_i,$ we get a projection in $\mathfrak{J}$. Moreover, since $\tau$ is normal we have 
    $$ \tau(q) = \bigvee_{i \in \Gamma_0}\tau(q_i) \leq z.$$
    Note that $0 < q \leq p $ since $0< q_i \leq p$ for all $i \in \Gamma_0$. We finally proof that $\tau(q) = z$. Otherwise, $z - \tau(q) > 0$ and $ p - q \in \PP(\mathfrak{J})\backslash\{0\}$ (since $z < \tau(p)$). Then 
    $ 0 < z-\tau(q) \leq \tau(p -q)$. As before, the set $\PP(\mathfrak{J})^-_{\tau, z-\tau(q), \tau(p-q)}$ is also non-empty, there exists a non-zero projection $r \in \PP(\mathfrak{J})$ with $r \leq p -q$,  $\tau(r) \leq z- \tau(q)$, and $r \perp q$ since $p-q \perp q$. Then the set $\{q_i\}_{i \in \Gamma_0} \cup \{q + r\}$ is totally ordered and strictly bigger than $\{q_i\}_{i \in \Gamma_0}$, a contradiction with the fact that $\{q_i\}_{i \in \Gamma_0}$ is maximal. 
\end{proof}

We shall discuss next some preliminary assumptions assumed along this section. 

\begin{remark}\label{remark assumptions for section 5}
Let $\mu : \PP(\mathfrak{J}) \rightarrow \mathbb{R}$ be a bounded finitely additive signed measure, where $\mathfrak{J}$ is a modular JBW$^*$-algebra without type $I_2$ part. By structure theory there exists a sequence of central projections $(z_n)_{n\in \mathbb{N}\cup \{0\}}$ such that $z_n \circ \mathfrak{J}$ is a (possibly zero) type $I_n$ JBW$^*$-algebra, $z_2=0$, and $z_0\circ\mathfrak{J}$ is a  (possibly zero) type $II_1$ JBW$^*$-algebra (see \cite[\S 5]{HOS} or subsection~\ref{subsec: structure}). We can always assume that $\mathfrak{J}$ contains no summands of type $I_n$ with $1\leq n\leq m_0$ for a fixed natural number $m_0$ --consequently $\mathfrak{J}$ is a modular JW$^*$-algebra. Namely, if we decompose $\mathfrak{J} = \left( \bigoplus_{1\leq n\leq m_0}^{\ell_{\infty}} z_n \circ \mathfrak{J}\right) \bigoplus^{\ell_{\infty}} \left( \bigoplus_{n\geq m_0+1}^{\ell_{\infty}} z_n \circ \mathfrak{J}\right) \bigoplus^{\ell_{\infty}} z_0\circ \mathfrak{J}.$ The restriction of $\mu$ to the lattice $\PP\left( \bigoplus_{1\leq n\leq m_0}^{\ell_{\infty}} z_n \circ \mathfrak{J} \right)$ admits an extension to a bounded linear functional on $\bigoplus_{1\leq n\leq m_0}^{\ell_{\infty}} z_n \circ  \mathfrak{J}$ (cf. \Cref{theo: lin_type_I}). So, it suffices to prove that the restriction of $\mu$ to the lattice of projections of $\mathfrak{J}_{m_0+1}:=\left( \bigoplus_{n\geq m_0+1}^{\ell_{\infty}} z_n \mathfrak{J}\right) \bigoplus^{\ell_{\infty}} z_0\circ \mathfrak{J}$ admits an extension to a bounded linear functional on $\mathfrak{J}_{m_0+1}$.\smallskip

Furthermore, by \Cref{theo: lin_type_I} we can find a bounded linear functional $\overline{\mu}_0\in \mathfrak{J}_{m_0+1}^*$ such that $(\mu-\overline{\mu}_0)|_{\PP\left(\mathfrak{J}_{m_0+1} \right)}$ vanishes on every projection belonging to a finite sum of factors of type $I_n$. Along the rest of this section we shall assume that $\mathfrak{J}$ contains no summands of type $I_n$ with $1\leq n\leq m_0$ for a fixed $3\leq m_0\in \mathbb{N}$, and $\mu$ vanishes on every projection belonging to a finite sum of factors of type $I_n$.
\end{remark}

Let us return to the mappings $\alpha_{\mu}$ and $V_{\mu}$ introduced in \Cref{def: alpha and variation}. The following lemma is essentially in \cite[Lemma 2.2]{BunceWright1994}.

\begin{lemma}\label{lemma: tracial_posit}
Let $\mu : \PP(\mathfrak{J}) \rightarrow \mathbb{R}$ be a bounded finitely additive signed measure, where $\mathfrak{J}$ is a modular JBW$^*$-algebra containing no type $I_n$ summands for all $1\leq n\leq 3$. Let $\tau$ denote the normal unital centre-valued faithful trace on $\mathfrak{J}$. 
Assume additionally that $\mu$ vanishes on every projection belonging to a finite sum of factors of type $I_n$. Then the mappings $\alpha_{\mu}|_{\PP(Z(\mathfrak{J}))}, \alpha_{-\mu}|_{\PP(Z(\mathfrak{J}))}$ extend to  bounded linear functionals $\overline{\alpha}_{\mu}$ and $\overline{\alpha}_{-\mu}$ on $Z(\mathfrak{J})$ and the inequality $$V_{\mu}(p) \leq (\overline{\alpha}_{\mu}+\overline{\alpha}_{-\mu}) \tau (p) \leq 2 V_{\mu}(p)$$ holds for all $p\in \PP(Z(\mathfrak{J})).$ Moreover, the functional $\phi_{_\mu}= (\overline{\alpha}_{\mu}+\overline{\alpha}_{-\mu}) \tau\in \mathfrak{J}^*$ vanishes on every central projection belonging to a finite sum of factors of type $I_n$. Up to replacing $\mu$ with an appropriate multiple of it, we can assume that $\phi_{_\mu} (\unit) =1$.
\end{lemma}

The proof given in \cite[Lemma 2.2]{BunceWright1994} in the case of von Neumann algebras actually works in our setting, so details are omitted.\smallskip

We are now in a position to apply the intermediate value property in \Cref{prop: trace_surject} in the next Jordan version of \cite[Lemma 1]{Yeadon1984}. 

\begin{lemma}\label{lemma: Yeadon IVT for functionals vanishing on bounded type projections} 
Let $\mathfrak{J}$ be a modular JBW$^*$-algebra containing no type $I_n$ summands for all $1\leq n\leq 3$, and let $\tau$ denote the normal unital centre-valued faithful trace on $\mathfrak{J}$. Suppose that for a positive $\varphi\in Z(\mathfrak{J})_{sa}^*$, the functional $\phi = \varphi \tau$ vanishes on every projection belonging to a finite sum of factors of type $I_n$. Then, for each $p\in \PP(\mathfrak{J})$ and each real $\alpha$ with $0<\alpha < \phi(p)$, there exists $q \in \PP(\mathfrak{J})$ satisfying $q \leq p$ and $\phi(q) = \alpha$.
\end{lemma}

\begin{proof} The proof is very similar to that in \cite[Lemma 1]{Yeadon1984}, we include it here for completeness. We can clearly assume that $\mathfrak{J}$ is type $I$ modular or type $II_1$.\smallskip
	
Assume first that $\mathfrak{J}$ is type $II_1$.  Since $\frac{\alpha}{\phi (p)}<1$, \Cref{prop: trace_surject} assures the existence of a projection $q\leq p$ such that $\tau (q) = \frac{\alpha}{\phi (p)} \tau (p)$, and thus $\phi (q) = \alpha$. \smallskip

We assume next that $\mathfrak{J}$ is type $I$ modular, and thus there exists a sequence of central projections $(z_n)_{n\geq 4}$ such that $z_n\circ \mathfrak{J}$ is a type $I_n$ JBW$^*$-algebra and $\mathfrak{J} = \bigoplus_{n\geq 4}  z_n\circ \mathfrak{J}$. It is well known that $z_n\circ \mathfrak{J}\cong C(\Omega_n, \mathfrak{F}_n)$, where $\Omega_n$ is a hyper-Stonean space and $\mathfrak{F}_n$ is a JBW$^*$-algebra factor of type $I_n$ ($n\geq 4$). Since the unique tracial state on $\mathfrak{F}_n$ only takes the values $\{0,\frac1n, \frac2n, \ldots, 1\}$ on $\PP (\mathfrak{F}_n)$ (see \cite[Proposition 5.22]{AlfsenShultz2003}), the restriction of $\tau$ to $\PP (z_n\circ \mathfrak{J})$ only contains functions in $C(\Omega_n)$ whose image is inside $\{0,\frac1n, \frac2n, \ldots, 1\}$. Working on each summand, we can find a projection $q\in \PP(\mathfrak{J})$ satisfying $z_n\circ \tau(q) \leq \frac{\alpha}{\phi (p)} z_n\circ \tau (p) \leq z_n\circ \tau(q) + \frac1n z_n$, for all $n\geq 4$, that is, \begin{equation}\label{eq 2407} \tau(q) \leq \frac{\alpha}{\phi (p)} \tau (p) \leq \tau(q) + \sum_{n\geq 4}^{\infty} \frac1n z_n.
\end{equation} If we apply $\varphi$ on \eqref{eq 2407} we get $$\begin{aligned}
\phi (q) =  \varphi (\tau(q)) &\leq \frac{\alpha}{\phi (p)} \varphi(\tau (p)) = \alpha \\ &\leq \varphi(\tau(q)) + \varphi\left(\sum_{n\geq 4}^{\infty} \frac1n z_n\right) = \phi(q) + \phi \left(\sum_{n\geq 4}^{\infty} \frac1n z_n\right).
\end{aligned}$$ It follows from the hypothesis on $\phi$ that $\displaystyle \phi \left(\sum_{n\geq 4}^{\infty} \frac1n z_n\right) = \phi \left(\sum_{n\geq m}^{\infty} \frac1n z_n\right)$ for all $m\geq 4$, and since the sequence $\displaystyle \left(\sum_{n\geq m}^{\infty} \frac1n z_n\right)_{m\geq 4}$ converges to $0$ in norm, the continuity of $\phi$ gives $\displaystyle \phi \left(\sum_{n\geq 4}^{\infty} \frac1n z_n\right) =0$.  
\end{proof}

\subsection{Technical arguments}\label{subsec: technical arguments}\ \smallskip

Along the rest of this section, let $\mu : \PP(\mathfrak{J}) \rightarrow \mathbb{R}$ be a bounded finitely additive signed measure, where $\mathfrak{J}$ is a modular JBW$^*$-algebra containing no type $I_n$ summands for all $1\leq n\leq 3$. Let $\tau$ denote the normal unital centre-valued faithful trace on $\mathfrak{J}$. Assume additionally that $\mu$ vanishes on every projection belonging to a finite sum of factors of type $I_n$. Let $\phi_{_\mu}\in \mathfrak{J}^*$ stand for the positive functional given by \Cref{lemma: tracial_posit}, which is assumed to be unital. As in \cite{BunceWright1994}, for each each $e \in \PP(\mathfrak{J})$ and each real number $t\in (0,\phi_{_\mu}(e))$  we define
$$ \mathcal{P}(e,t):=\{ p \in \PP(\mathfrak{J}) \, : \, p \leq e \, , \, \phi_{_\mu}(p) = t\}.\label{def: P(e,t)} 
$$ A combination of \Cref{lemma: tracial_posit} and \Cref{lemma: Yeadon IVT for functionals vanishing on bounded type projections} shows that $ \mathcal{P}(e,t)$ is a non-empty set. We also define 
    $$ 
    \begin{aligned}
        M(e,t) &= \sup\{ \mu(p) \, : \, p \in \mathcal{P}(e,t)\}, \hbox{ and }
        m(e,t) = \inf\{ \mu(p) \, : \, p \in \mathcal{P}(e,t)\}.
    \end{aligned}
    $$
    
Most of the technical results in this subsection are inspired from those in \cite{BunceWright1994}, the proofs are actually valid by just replacing Yeadon's results in \cite{Yeadon1984} by their Jordan versions in \Cref{lemma: tracial_posit} and \Cref{lemma: Yeadon IVT for functionals vanishing on bounded type projections}. We opted for including some full arguments for completeness reasons.\smallskip

As in \cite[Lemma 3.1]{BunceWright1994}, the following result is a direct consequence of the definitions.      
    
\begin{lemma}\label{l properties of M and m} 
Under the hypotheses assumed in this subsection, the following statements hold for all $e \in \PP(\mathfrak{J})$ and all $t \in (0,\phi_{_\mu}(e))$:
     \begin{enumerate}[$(a)$]
         \item $-V_{_\mu}(e) \leq m(e,t) \leq M(e,t) \leq V_{_\mu}(e),$
         \item $M(e,t) = \mu(e) - m(e,\phi_{_\mu}(e) -t).$
     \end{enumerate} Consequently, $M(e, \cdot) $ and $m(e,\cdot)$ are bounded functions on $[0,\phi_{_\mu}(e)]$. \hfill$\Box$
\end{lemma} 

The next lemma relies on \Cref{lemma: Yeadon IVT for functionals vanishing on bounded type projections}.

\begin{lemma}\label{lemma: 3.3} Under the hypotheses assumed in this subsection, let $e \in \PP(\mathfrak{J})$  with $\phi_{_\mu} (e) > 0$. Then the following statements hold for all $0 < t_1 \leq t_2 < \phi_{_\mu} (e)$:
\begin{enumerate}[$(a)$]
\item Given $\varepsilon > 0$ and a projection $ p \in \mathcal{P} (e, t_1)$ such that $\mu(p) > M(e, t_1) - \varepsilon,$ then there exists a projection $q \geq p$ such that $q \in \mathcal{P}(e, t_2)$ and $\mu(q) > M(e, t_2) - \varepsilon$.
\item Given $\varepsilon > 0$ and a projection $ p \in \mathcal{P}(e, t_2)$ such that $\mu(p) > M(e, t_2) - \varepsilon,$ then there exists a projection $q \leq p$ such that $q \in \mathcal{P}(e, t_1)$ and $\mu(q) > M(e, t_1) - \varepsilon$.
\end{enumerate}
\end{lemma}

\begin{proof}$(a)$ Set $\delta = \mu(p) - M(e,t_1) + \varepsilon >0$ and take $h \in \mathcal{P}(e, t_2)$ such that $\mu(h) > M(e,t_2) - \delta$. By applying \cite[Proposition 5.2.3$(i)$]{HOS} we arrive to 
    $$ h - h \wedge (\unit-p) \sim p - p\wedge (\unit - h) \leq  p.  $$ So $t_2 - t_1 = \phi_{_\mu} (h) - \phi_{_\mu} (p) \leq \phi_{_\mu} ((\unit -p)\wedge h).$
    Therefore, by \Cref{lemma: Yeadon IVT for functionals vanishing on bounded type projections} we can choose a projection $h_0 \in \mathcal{P}(e, t_2-t_1)$ with $h_0 \leq (\unit - p)\wedge h\leq h\leq e$. Then $h-h_0 \in \mathcal{P}(e, t_1)$, and consequently
    $$ M(e,t_1) \geq \mu(h-h_0) = \mu(h) - \mu(h_0) > M(e, t_2) - \delta - \mu(h_0),$$
    and $ \mu(h_0) > M(e, t_2) - M(e, t_1)-\delta.$
By defining $q = p + h_0$, we have $q \geq p,$ $q \in \mathcal{P}(e, t_2)$, and 
$$ \begin{aligned}
        \mu(q) &= \mu(p) + \mu(h_0) > \mu(p) + M(e,t_2) - M(e,t_1) - \delta \\
        &=M(e,t_2)- (M(e,t_1) - \mu(p) + \delta) = M(e,t_2) - \varepsilon. 
    \end{aligned}
    $$

$(b)$ In this case, we set $\delta = \mu(p) - M(e,t_2) + \varepsilon> 0,$ and we take $h \in \mathcal{P}(e, t_1)$ such that $\mu(h) > M(e, t_1) - \delta$. By applying \cite[Proposition 5.2.3$(i)$]{HOS} we get $$ p - p \wedge (\unit-h) \sim h- h \wedge (\unit - p) \leq h,$$ and by the properties of the trace we arrive to  
$$ t_2-t_1 =\phi_{_\mu} (p) - \phi_{_\mu} (h) \leq \phi_{_\mu}((\unit -h)\wedge p).$$ \Cref{lemma: Yeadon IVT for functionals vanishing on bounded type projections} guarantees the existence of $h_0 \in \mathcal{P}(e, t_2-t_1)$ with $h_0 \leq p\wedge(\unit -h)\leq p$ (and $h\perp h_0$). Then $\phi_{_\mu} (p-h_0) = t_1$ and
    $$ 
    \begin{aligned}
        M(e,t_1) & \geq \mu(h-h_0) = \mu(h) - \mu(h_0) \\
        &> M(e, t_1) - \delta - \mu(h_0) = M(e,t_1) + M(e,t_2) - \mu(p) - \varepsilon - \mu(h_0).
    \end{aligned}
    $$
    Setting $q = p - h_0$ we get $q \leq p$, $q \in \mathcal{P}(e,t_1)$ and 
    $$ \mu(q) = \mu(p) - \mu(h_0) > M(e, t_1) - \varepsilon.$$
\end{proof}

When in the proof of \cite[Lemma 3.4]{BunceWright1994}, our previous \Cref{lemma: 3.3} and \Cref{lemma: Yeadon IVT for functionals vanishing on bounded type projections} replace \cite[Lemma 3.3 and Lemma 2.3]{BunceWright1994}, respectively, the arguments are literally valid to get the following lemmas.  

\begin{lemma}\label{l 5.8} Under the hypotheses assumed in this subsection, let $e$ be a projection in $\mathfrak{J}$ satisfying $\phi_{_\mu}(e) \neq 0$. Then the mapping $M_e: (0, \phi_{_\mu}(e)) \rightarrow \mathbb{R}$, $M_{e} (t) = M(e,t)$ (respectively, $m_e: (0, \phi_{_\mu}(e)) \rightarrow \mathbb{R}$, $m_{e} (t) = m(e,t)$) is continuous and concave (respectively, convex), and the limits $\displaystyle \lim_{t\rightarrow 0}M(t)$ and $\displaystyle \lim_{t\rightarrow \phi_{_\mu}(e)}M(t)$ (respectively, $\displaystyle \lim_{t\rightarrow 0}m(t)$ and $\displaystyle \lim_{t\rightarrow \phi_{_\mu}(e)}m(t)$) exist. \hfill$\Box$
\end{lemma}

The statement concerning $m(e,\cdot)$ follows from \Cref{l properties of M and m}$(b)$.\smallskip
 
As in \cite[Lemma 3.5]{BunceWright1994}, the next lemma is a corollary of the previous \Cref{l 5.8}. The proof is omitted. 

\begin{lemma}\label{lemma: 3.5} Under the hypotheses assumed in this subsection, let $e \in \PP(\mathfrak{J})$ with $\phi_{_\mu}(e) > 0$, and let $0 <t_1,t_2 \in \mathbb{R}$ such that $t_1 + t_2 \leq \phi_{_\mu}(e)$. Then the following statements hold:
\begin{enumerate}[$(a)$]
\item $m(e,t_1) + m(e, t_2) \leq m(e, t_1 + t_2).$
\item $M(e,t_1) + M(e, t_2) \geq M(e, t_1 + t_2).$
\item $\displaystyle\lim_{t\rightarrow 0} m_e (t) \leq 0$.
\item $\displaystyle \lim_{t\rightarrow 0}M_e(t) \geq 0$. \hfill$\Box$
\end{enumerate}
\end{lemma}  

The next lemma is essentially in \cite[Lemma 3.6]{BunceWright1994}, the change of nomenclature invites us include the proof. 

\begin{lemma}\label{lemma: Lambda_acot_mu}  Under the hypotheses assumed in this subsection, let $e,p$ be orthogonal projections in $\mathfrak{J}$ with $\phi_{_\mu}(e), \phi_{_\mu}(p) > 0$. Let $\delta = M(e+p, \phi_{_\mu}(e)) - \mu (e)$. Then there exists a real number $\lambda$ such that 
    $$ \overline{\mu} (U_e(a)) \geq \lambda \phi_{_\mu} (U_e(a)) -\delta, \hbox{ and }\  \overline{\mu}(U_p(a)) \leq \lambda \phi_{_\mu}(U_p(a)) +\delta,$$ for all $0\leq a \leq \unit$.
\end{lemma}

\begin{proof} Let us fix $0 < t \leq \min \{\phi_{_\mu}(e), \phi_{_\mu}(p)\}$. By \Cref{lemma: tracial_posit} and \Cref{lemma: Yeadon IVT for functionals vanishing on bounded type projections} (see page~\pageref{def: P(e,t)}), we can take $e_1\in \mathcal{P}(e,t)$ and $p_1\in \mathcal{P}(p,t)$. Note that 
    $$ \phi_{_\mu} (e -e_1 +p_1) = \phi_{_\mu} (e) -t +t = \phi_{_\mu} (e),$$
    which implies that 
    $$ \mu(e -e_1 + p_1) \leq M(e + p, \phi_{_\mu} (e)) = \mu(e) + \delta. $$
    It is clear that $e_1 \circ p_1 = 0$ ($e_1\perp p_1$), therefore $\overline{\mu}$ is linear on $W^*(\unit, e, e_1,p_1)$ by definition since the latter is a JW$^*$-subalgebra of $\mathfrak{J}$. Therefore,  
    $$ \mu(e -e_1 + p_1) = \overline{\mu}(e-e_1+p_1) = \mu(e) -\mu(e_1) + \mu(p_1),$$ which proves that  $\mu(p_1) \leq \mu(e_1) + \delta.$ The arbitrariness of $e_1\in  \mathcal{P}(e,t)$ and $p_1\in \mathcal{P}(p,t)$ imply that 
    \begin{equation}\label{eq: T_acot}
        \mu(p_1) \leq m(e,t) +\delta,\, \hbox{ and } \, M(p,t) \leq m(e,t) + \delta. 
    \end{equation} Consider now the following sets in $\mathbb{R}^2$
    $$ 
    \begin{aligned}
        U_1 &= \{(t,y) \, : \, 0< t< \phi_{_\mu}(e) \text{ and } y > m(e,t)\},\\
        U_2 &= \{(t,y) \, : \, 0< t< \phi_{_\mu}(p) \text{ and } y < M(p,t)\}.
    \end{aligned}
    $$ Clearly $U_1$ and $U_2$ are non-empty open and disjoint by the previous arguments. It follows from the convexity of $m(e,\cdot)$ and the concavity of $M(e,\cdot)$ (cf. \Cref{l 5.8}) that $U_1$ and $U_2$ are convex subsets. So, by the Hahn-Banach theorem we can separate $U_1$ and $U_2$ by an affine hyperplane, in particular, there exist $\lambda,\rho \in \mathbb{R}$ satisfying
    $$ M(e,t) + \varepsilon' \leq \lambda t + \rho \leq M(p,t) + \varepsilon' $$ for all $0 < t \leq \min \{\phi_{_\mu}(e),\phi_{_\mu}(p)\}$ and $\varepsilon'>0$, which in turn gives
    $$ M(e,t) \leq \lambda t + \rho \leq M(p,t)\leq  m(e,t) + \delta , $$ for all $0 < t < \min \{\phi_{_\mu}(e),\phi_{_\mu}(p)\}.$
    By letting $t \rightarrow 0$ and applying \Cref{lemma: 3.5}$(c)$ and $(d)$ we get 
    $ 0 \leq \rho \leq \delta.$\smallskip

    Finally, let $q \in \PP(\mathfrak{J})$ with $q \leq e$. Then $$\mu(q) \geq m(e, \phi_{_\mu}(q)) \geq \lambda \phi_{_\mu} (q) + \rho -\delta \geq \lambda \phi_{_\mu} (q) - \delta.$$
    Let $ a \in \mathfrak{J}$ such that $0 \leq a \leq \unit$, and take $U_e(a)$. By \cite[Proposition 4.2.3]{HOS} for each positive $\varepsilon''$ there are mutually orthogonal projections $\{q_1, \ldots, q_n\}$ in $W^*(\unit, U_e(a)) $ and non-negative real numbers $\{\alpha_1,\ldots,\alpha_n\}$ such that the inequalities $\displaystyle \sum_{k=1}^n \alpha_k\leq \|a\| \leq 1$ and 
    $\displaystyle \left\| U_e(a) - \sum_{k=1}^n \alpha_k q_k \right\|<\varepsilon''$ hold.
    It follows from the above conclusions that
    $$\overline{\mu}\left(\sum_{k=1}^n \alpha_k q_k \right) = \sum_{k=1}^n \alpha_k \overline{\mu} (q_k) = \sum_{k=1}^n \alpha_k {\mu} (q_k) \geq  \sum_{k=1}^n \alpha_k \left(\lambda \phi_{_\mu} (q) - \delta\right)$$
  $$ \hspace{1.75cm}= \lambda \phi_{_\mu} \left( \sum_{k=1}^n \alpha_k \right) - \delta \sum_{k=1}^n \alpha_k \geq \lambda \phi_{_\mu} \left( \sum_{k=1}^n \alpha_k \right) - \delta.$$ The arbitrariness of $\varepsilon''$ and the continuity of $\overline{\mu}|_{_{W^*(\unit, U_e(a))}}$ and $\phi_{_\mu}$ can be now applied to conclude that  $$ \overline{\mu}(U_e(a)) \geq \lambda  \phi_{_\mu}(U_e(a)) - \delta.$$
    We can similarly get 
    $ \overline{\mu}(U_p(a)) \leq \lambda \phi_{_\mu} (U_p(a)) +  \delta.$
\end{proof}


We continue with a Jordan version of \cite[Lemma 4.1]{BunceWright1994}.

\begin{lemma}\label{lemma: e_pm} Let $\mathfrak{J}$ be a JW$^*$-algebra without $I_2$ part. Let us take  $0 <\varepsilon< \frac{1}{2}$ and  $p, q , e \in \PP(\mathfrak{J})$ such that $p \circ q = 0$ and $e \leq p + q$. Define 
    $$ c = \frac{1}{2}p + \left(\frac{1}{4} p - \varepsilon^4 U_pU_e(q)\right)^{\frac12}, \quad d =\frac{1}{2}q - \left(\frac{1}{4} q - \varepsilon^4 U_qU_e(p)\right)^{\frac12},$$
    $$e_{-} = c +d+\varepsilon^2 (2\{p,e,q\}), \hbox{ and } e_{+} = c + d - \varepsilon^2 (2\{p,e,q\}).$$
    Then $e_{-}, e_+$ are projections dominated by $p +q$. We further know that $e_- \sim p \sim e_+$, $ 0 \leq p-c \leq \frac{1}{2}\varepsilon^4 p$ and $0 \leq d \leq \frac{1}{2} \varepsilon^4 q$.
\end{lemma}

\begin{proof} We can assume, by the assumptions in this subsection that $\mathfrak{J}$ is a JBW$^*$-subalgebra of some von Neumann algebra $A$, where the associative product of $A$ will be denoted by juxtaposition. We can also assume that $A$ and $\mathfrak{J}$ share the same unit. Clearly, the elements $p,q,$ and $e$ are projections in $A$ satisfying the same hypotheses. \smallskip
    
Set $b = peq\in A$. Note that $b$ need not be an element in $\mathfrak{J}$. However, the elements 
    $$ c =\frac{1}{2}p + \left(\frac{1}{4} p - \varepsilon^4 bb^*\right)^{\frac12} = \frac{1}{2}p + \left(\frac{1}{4} p - \varepsilon^4 U_p U_e U_q (\unit)\right)^{\frac12},$$ $$\hbox{ and } d =\frac{1}{2}q - \left(\frac{1}{4} q - \varepsilon^4 b^*b\right)^{\frac12} = \frac{1}{2}q - \left(\frac{1}{4} q - \varepsilon^4 U_q U_e U_p (\unit)\right)^{\frac12}$$ lie in $\mathfrak{J}$, and the same occur to  
    $$e_{-} = c +d+\varepsilon^2 (b+b^*) = c +d+2 \varepsilon^2  \{p,e,q\}, \hbox{ and } e_{+} = c + d - 2 \varepsilon^2  \{p,e,q\}.$$ By applying 
    Lemma 11.10 in \cite{Maeda1989} with $\lambda = \varepsilon^2$ we conclude that $e_{-}, e_+$ are projections dominated by $p +q$, $ 0 \leq p - c \leq \frac{1}{2} \varepsilon^4 p$, and $0 \leq d \leq \frac{1}{2} \varepsilon^4 q$ in $A$ (and hence in $\mathfrak{J}$), and $\| p - e_{-}\|, \ \|p - e_+\| <\frac14 + \frac{1}{16}<1$. It is well known that if $p$ and $q$ are projections in a unital C$^*$-algebra $B$ satisfying $\|p-q\|< 1$, then $p$ and $q$ are unitarily equivalent \cite{RoLarLaustBook2000}. Thus, $e_- $ and $ p$ and $e_+$ and $p$ are unitarily equivalent in $A$. An application of \cite[Corollary 1.3]{BunceWright1989} shows that $e_- \sim p \sim e_+$ in $\mathfrak{J}$. 
\end{proof}

Our next goal is a version of \Cref{l subalgebra isomorphic to S2} for modular JW$^*$-algebras without type $I_2$ summands.

\begin{lemma}\label{l subalgebra isomorphic to S2(C) in a modular} Under the hypotheses assumed in this subsection, let $\mathfrak{B}$ be a JBW$^*$-subalgebra of $\mathfrak{J}$ which is JB$^*$-isomorphic to $S_2 (\mathbb{C})$. Then the restriction of $\overline{\mu}$ to $\mathfrak{B}$ is linear.
\end{lemma}

\begin{proof} Note that $\mathfrak{J}$ is reversible (cf. \cite[Theorem 5.3.10]{HOS}). Let $e$ denote the unit element in $\mathfrak{B}\cong S_2 (\mathbb{C})$. We can then find a couple of orthogonal projections $q_1,q_2$ and a symmetry $w_{_{12}}$ in $\mathfrak{B}$ corresponding to the matrices $\left(\begin{matrix}
		1& 0 \\ 0 &0
	\end{matrix}\right),$ $\left(\begin{matrix}
		0& 0 \\ 0 &1
	\end{matrix}\right),$ and $\left(\begin{matrix}
		0& 1 \\ 1 &0
	\end{matrix}\right),$ respectively. Observe that $v_{_{12}} = w_{_{12}} q_1$ is a partial isometry in $R(\mathfrak{B}_{sa})^{-}\subseteq R(\mathfrak{J}_{sa})^{-}$ satisfying $v_{_{12}}^* v_{_{12}} = q_1$ and $v_{_{12}} v_{_{12}}^* = q_2$. \smallskip

We can obviously reduce our arguments to the following cases: $(1)$ $\mathfrak{J}$ is type $I$ modular; $(2)$ $\mathfrak{J}$ is type $II_1$.\smallskip

$(1)$ $\mathfrak{J}$ is type $I$ modular. By structure theory, there exists a sequence of central projections $(z_n)_{n\geq 4}$ such that $\mathfrak{J} = \bigoplus_{n\geq 4}^{\infty}  z_n\circ \mathfrak{J}$, and $z_n\circ \mathfrak{J}$ is a (possibly zero) type $I_n$ JBW$^*$-algebra of the form $C(\Omega_n, \mathfrak{F}_n)$, where $\Omega_n$ is a hyper-Stonean space and $\mathfrak{F}_n$ is a finite-dimensional JBW$^*$-algebra factor of type $I_n$ ($n\geq 4$). Fix $n\geq 4$ and consider the projection $z_n \circ q_1.$ We can clearly assume that $z_n\circ q_1\neq 0$ for all $n\geq 4$. Arguing in the type $I$ modular JBW$^*$-subalgebra $U_{(z_n\circ q_1)}  (\mathfrak{J})$, as in the proof of \cite[Proposition 2.3, Case II]{EscoPeVi2025}, we can find (possibly zero) pairwise orthogonal projections $p_{1n},p_{2n},$ and $q_{1n}$ such that $q_{1n}$ is abelian, $p_{1n}\sim p_{2n}$, and $z_n\circ q_1 = p_{1n}+p_{2n} + q_{1n}$. We define the following projections $p_1,p_2$ and $r_1$ determined, uniquely, by $p_1\circ z_n  = p_{1n},$ $p_2\circ z_n  = p_{2n},$ and $r_1\circ z_n = q_{1n}.$ It follows from these definitions that $q_1 = p_{1} +p_{2} + r_{1}$, where $p_{1},p_{2},$ and $r_{1}$ are mutually orthogonal, $p_{1}\sim p_{2},$ and $r_1$ is abelian.\smallskip

Since $q_1\approx q_2$ (via $v_{_{12}}$, as in the proof of \Cref{p four BW equivalent projections}$(b)$), the decomposition of $q_1$ can be transferred to an orthogonal decomposition of $q_2$ in the form $q_2 = p_{3} +p_{4} + r_{2},$ where $p_{3}\sim p_{4},$ and $r_2$ is abelian. 
\smallskip

Let $\mathfrak{B}_1$ (respectively, $\mathfrak{B}_2$) denote the JBW$^*$-subalgebra of $\mathfrak{J}$ generated by $p_1+p_2, p_3+p_4,$ and $v_{_{12}} (p_1+p_2) +  (p_1+p_2) v_{_{12}}^*$ (respectively, $r_1,$ $r_2$, and $v_{_{12}} r_1 +  r_1 v_{_{12}}^*$). According to this construction, $\mathfrak{B}$ is a JBW$^*$-subalgebra of $\mathfrak{B}_1\bigoplus^{\ell_{\infty}}\mathfrak{B}_2$. \smallskip

If $p_1\neq 0 $, by \Cref{p four BW equivalent projections}$(b)$, there exists a JBW$^*$-subalgebra $\tilde{\mathfrak{B}}_1$ of $\mathfrak{J}$ which is Jordan $^*$-isomorphic to $S_4(\mathbb{C}),$ contains the JB$^*$-subalgebra $\mathfrak{B}_1$ and is orthogonal to $\mathfrak{B}_2$. \Cref{theo: lin_type_In} (see also \Cref{r after theorem type In}) assures that $\overline{\mu}|_{_{\tilde{\mathfrak{B}}_1}}$ (and hence $\overline{\mu}|_{_{{\mathfrak{B}}_1}}$) is linear. \smallskip

We deal next with ${{\mathfrak{B}}_2}$, which is Jordan $^*$-isomorphic to $S_2(\mathbb{C})$, where the minimal projections in its diagonal (i.e., $r_1$ and $r_2$) are two orthogonal abelian projections in $\mathfrak{J}$. It follows that $c(r_1+r_2) -(r_1+r_2)$ is a non-zero projection in $\mathfrak{J}$ with $c(r_1) = c(r_2) =c(r_1+r_2)$. Lemma 5.3.2$(iii)$ in \cite{HOS} proves that $c(r_1+r_2) -(r_1+r_2)$ dominates an abelian projection $r_3$ with $c(r_1) = c(r_3)$, and thus $r_3\sim r_1\sim r_2$ (cf. \cite[Lemma 5.3.2$(ii)$]{HOS}). Arguing as in the proof of \Cref{p four BW equivalent projections}, it is not hard to see that there is another JBW$^*$-subalgebra $\tilde{\mathfrak{B}}_2$ of $\mathfrak{J}$ which is Jordan $^*$-isomorphic to $S_3(\mathbb{C}),$ contains the JB$^*$-subalgebra $\mathfrak{B}_2$ and is orthogonal to $\tilde{\mathfrak{B}}_1$. A new application of \Cref{theo: lin_type_In} (see also \Cref{r after theorem type In}) shows that $\overline{\mu}|_{_{\tilde{\mathfrak{B}}_2}}$ (and hence $\overline{\mu}|_{_{{\mathfrak{B}}_2}}$) is linear. This concludes the proof in this case.\smallskip

$(2)$ $\mathfrak{J}$ is type $II_1$. Let us observe that $q_1$ must be a modular projection, and there are no abelian projections under $q_1$. So $U_{q_1} (\mathfrak{J})$ is a type $II_1$ JBW$^*$-algebra. By the halving lemma (see \cite[5.2.14]{HOS}), there exists subprojections $p_1,p_2\leq q_1$ such that $p_1\sim p_2$ and $p_1+ p_2 = q_1$. As in the previous case (and in the proof of \Cref{p four BW equivalent projections}$(b)$), the decomposition of $q_1$ can be transferred to an orthogonal decomposition of $q_2$ in the form $q_2 = p_{3} +p_{4},$ with $p_3 \sim p_4$. We are in a position to apply \Cref{p four BW equivalent projections}$(b)$, which leads to the existence of a JBW$^*$-subalgebra $\tilde{\mathfrak{B}}\cong S_4 (\mathbb{C})$ containing $\mathfrak{B}$. The desired statement follows from \Cref{theo: lin_type_In} (see also \Cref{r after theorem type In}), as in the final argument of the previous case. 
\end{proof}

We can now mimic some of the ideas applied in the proof of \Cref{theo: lin_type_In}.

\begin{lemma}\label{l subalgebra isomorphic to S2 modular and properly non-modular} Under the hypotheses assumed in this subsection, let $e$ and $p$ be two projections in $\mathfrak{J}$, and let $W^*(\unit, e,p)$ denote the JBW$^*$-subalgebra of $\mathfrak{J}$ generated by $\{\unit, e,p\}$. Then the restriction of the quasi-linear extension of $\mu$ to $W^*(\unit, e,p)$ is linear. 
\end{lemma}

\begin{proof} \cite[Proposition 3.5]{PeSaa2025}  there exists a hyper-Stonean compact space $\Omega$, $m\in \mathbb{N}\cup \{0\},$ such that $W^*(\unit, e,p)$ is Jordan $^*$-isomorphic to $C(\Omega, S_2(\mathbb{C})) \oplus^{\ell_{\infty}} \mathbb{C}^m $.\smallskip
	
Let $C^0(K,S_{2}(\mathbb{C}))$ denote the subspace of all locally constant mappings, that is, the functions $a\in C(K,S_{2}(\mathbb{C}))$ whose image is finite. Clearly $C^0(K,S_{2}(\mathbb{C}))$ is a Jordan $^*$-subalgebra of $\mathfrak{J}$. \smallskip

Let us prove that $\overline{\mu}$ is linear on $C^0(K,S_{2}(\mathbb{C}))$. Take $a, b \in C^0(K,S_{2}(\mathbb{C}))$, and a partition $\{K_1, \dots, K_m\} $ of $K$ such that each $K_i$ is clopen and both $a$ and $b$ are constant mappings on each $K_i$ for all $i:1,\ldots, m$. Let $\chi_{_{K_i}}$ denote the characteristic function of the set $K_i$. The finite dimensional JB$^*$-subalgebra 
$\displaystyle \mathfrak{B}:= \bigoplus_{i=1,\ldots,m}^{\ell_{\infty}} S_{2}(\mathbb{C}) \chi_{_{K_i}}\cong \bigoplus_{i=1,\ldots,m}^{\infty} S_{2}(\mathbb{C}).$ Note that $S_{2}(\mathbb{C}) \chi_{_{K_i}}$ is a JBW$^*$-subalgebra of $\mathfrak{J}$ for all $i\in \{1,\ldots, m\}$. Since $\displaystyle \mathcal{P} (\mathfrak{B}) = \bigoplus_{i=1}^m \mathcal{P}(S_{2}(\mathbb{C})) \chi_{_{K_i}},$ it is not hard to check via \Cref{l subalgebra isomorphic to S2(C) in a modular} that $\mu|_{\mathcal{P}(\mathfrak{B})}$ admits a unique extension to a bounded linear functional on $\mathfrak{B},$ which, by uniqueness, must coincide with $\overline{\mu}|_{\mathfrak{B}}$ (cf. \Cref{prop_quasilin}). Since $a, b \in \mathfrak{B}$, we get $\overline{\mu}(\gamma a + \beta b) = \gamma \overline{\mu}(a) + \beta \overline{\mu}(b)$ ($\gamma, \beta \in \mathbb{C})$. Therefore $\overline{\mu}$ is linear on $C^0(K,S_{2}(\mathbb{C}))$.\smallskip

\Cref{prop_quasilin} now assures that $\overline{\mu}|_{C^0(K,S_{2}(\mathbb{C}))}$ is a bounded linear functional with $\|\overline{\mu}|_{C^0(K,S_{2}(\mathbb{C}))}\|\leq 2 (2 \alpha_{_{\mu}}(\unit) - \mu(\unit))$. Having in mind the norm density of $C^0(K,S_{2}(\mathbb{C}))$ in $C(K,S_{2}(\mathbb{C}))$, we can find a unique extension of $\overline{\mu}|_{C^0(K,S_{2}(\mathbb{C}))}$ to a bounded linear functional $\varphi: C(K,S_{2}(\mathbb{C}))\to \mathbb{C}$, and a posteriori, an extension to a bounded linear functional on $W^*(\unit, e,p )$, which is also denoted by $\varphi$. By construction, $\varphi = \overline{\mu} = \mu$ on $\PP(C^0(K,\mathfrak{F}_{n})\oplus^{\ell_{\infty}} \mathbb{C}^m)$, and the latter is dense in $\PP(W^*(\unit, e,p ))$. We can finally apply the uniform continuity of $\mu$ on $\PP(W^*(\unit, e,p ))$ (cf. \Cref{prop: cont_mu}) to deduce that $\varphi = \mu$ on $\PP(W^*(\unit, e,p ))$, and thus $\varphi$ is a bounded linear extension of $\mu$, and $\varphi = \overline{\mu}$.
\end{proof}

The next technical result, which seems to be a novelty by itself, simplifies part of the subsequent arguments.

\begin{lemma}\label{l new almost three projections} Let $\mathfrak{J}$ be a reversible JW$^*$-subalgebra of a von Neumann algebra $A$, and assume that both algebras share the same unit $\unit$. Let $p,q$ and $e$ be three projections in $\mathfrak{J}$ such that $p$ and $q$ are orthogonal and $e\leq p+q$. Then the elements $U_q (e) = \{q,e,q\}$ and $U_{p,q} (e) = \{p,e,q\}$ belong to the JW$^*$-algebra $W^*(\unit, p,e)$.\end{lemma}

\begin{proof} We agree to denote the product of elements in the von Neumann algebra $A$ by mere juxtaposition. Set $\overline{p} := p\vee e\in  W^*(\unit, e,p)\subseteq \mathfrak{J}\subseteq A$. By definition $\overline{p}\circ p =  \overline{p} p = p \overline{p}= p,$ and $\overline{p}\circ e =  \overline{p} e = e \overline{p}= e$. It is also clear, by construction, that $\overline{p}\leq p+q,$ It then follows from these identities that $\overline{p} = \overline{p} (p+q) = p +  \overline{p} q$, which implies that $\overline{p} q = \overline{p} -p = q \overline{p}$. By relying on the associative structure of $A$ we can easily check that $ q e = q \overline{p} e= (\overline{p}-p) e,$ $e q =e \overline{p} q= e (\overline{p}-p),$
	$$\begin{aligned}
		\{q,e,q\}=U_q (e) = U_q U_{\overline{p}} (e) &= q \overline{p} e \overline{p} q = (\overline{p}-p) e (\overline{p}-p) \\
		&= U_{\overline{p}-p} (e)=\{\overline{p}-p,e,\overline{p}-p\}\in W^*(\unit, e,p).
	\end{aligned}$$	 Finally, 
	$$\begin{aligned}  2 \{p, e, q\} &=  p e q + q e p = p e  \overline{p} q + q \overline{p} e p  = p e (\overline{p}-p) + (\overline{p}-p) e p  \\
		&= 2 \{p,e,\overline{p}-p\} = \{\overline{p},e,\overline{p}\} -\{{p},e,{p}\} - \{\overline{p}-p,e,\overline{p}-p\} \\
		&= e -\{{p},e,{p}\} - \{q,e,q\} \in W^*(\unit, p, e).
	\end{aligned}$$
\end{proof}

\begin{lemma}\label{lemma: pierce_acot_theta_lin}
Under the hypotheses assumed in this subsection, let $p,q$ be orthogonal projections in $\mathfrak{J}$ such that $\mu (p) > M(p +q, \phi_{_\mu}(p)) - \frac{1}{4} \varepsilon^4$, where $0 <\varepsilon < \frac{1}{3}$. Then $ | \overline{\mu}(e) - 
 \overline{\mu}(U_p(e)) - \overline{\mu}(U_q(e))| < \varepsilon,$ for all projection $e \leq p +q$. Furthermore, for each $j \in \{1,2,3,4\}$, let $\theta_j$ be a bounded linear functional satisfying 
 $$
 \begin{aligned}
     \theta_j(x) &= \theta_j(U_p(x)), \text{ for } j \in \{1,2\}, \hbox{ and } \theta_j (x) = \theta_j(U_q(x)), \text{ for } j \in \{3,4\}.
 \end{aligned}
 $$
Finally, let $\gamma_1,\gamma_2,\gamma_3,$ and $\gamma_4$ be positive real numbers satisfying  
$$ \theta_1 (a) - \gamma_1 \leq \overline{\mu}(U_p(a)) \leq \theta_2(a) + \gamma_2,$$ and 
 $$   \theta_3 (a) - \gamma_3 \leq \overline{\mu}(U_q(a)) \leq \theta_4(a) + \gamma_4,$$ for all $a \in \mathfrak{J}$ with $0 \leq a \leq \unit$.
Then the inequalities 
\begin{equation}\label{eq: RN_1}
    (\theta_1 + \theta_3)(a) - (\gamma_1 + \gamma_3 + \varepsilon) \leq \overline{\mu}(U_{p+q}(a)) \leq (\theta_2 + \theta_4)(a) + (\gamma_2 + \gamma_4 + \varepsilon),
\end{equation} hold for all $0 \leq a \leq \unit$ in $\mathfrak{J}$.
\end{lemma}

\begin{proof} As observed before, $\mathfrak{J}$ is a reversible JW$^*$-algebra (cf. \cite[Theorem 5.3.10]{HOS}). Under these conditions, there exists a von Neumann algebra $A$ satisfying that $\mathfrak{J}$ is a reversible JBW$^*$-subalgebra of $A$. The product of elements in the von Neumann algebra $A$ will be denoted by juxtaposition.\smallskip 
	
Let us take $e \in \PP(\mathfrak{J})$ such that $e \leq p +q$. By applying \Cref{l new almost three projections} and its proof, we deduce that $\{q,e,q\}=U_q (e)= U_{\overline{p}-p} (e)=\{\overline{p}-p,e,\overline{p}-p\},$  $\{p,e, \overline{p}-p\}=\{p,e,q\}= U_{p,q} (e) \in W^*(\unit, p, e),$ where $\overline{p} := p\vee e\in  W^*(\unit, e,p)$.\smallskip
  
We shall next show that $|\overline{\mu}(2 \{p,e,q\})|$ is less than $\varepsilon$. To get the desired conclusion, we apply \Cref{lemma: e_pm} to $p,\overline{p}-p$ and $e$ to deduce that, taking $$ c = \frac{1}{2}p + \left(\frac{1}{4} p - \varepsilon^4 U_p U_e(\overline{p}-p)\right)^{\frac12}\in W^*(\unit, p, e),$$ and $$d =\frac{1}{2}(\overline{p}-p) - \left(\frac{1}{4} (\overline{p}-p) - \varepsilon^4 U_{\overline{p}-p} U_e(p)\right)^{\frac12}\in W^*(\unit, p, e),$$ then there exist projections $e_+, e_-\in W^*(\unit, p, e)$ satisfying $$e_{\pm} \sim p, \ 0 \leq p -c \leq \frac{1}{2}\varepsilon^4 p,\ 0 \leq d \leq \frac{1}{2}\varepsilon^4 (\overline{p}-p),$$ and $$e_{\pm} -c-d = \mp \varepsilon^2 (2 \{p,e, \overline{p}-p\}) = \mp \varepsilon^2 (2 \{p,e, q\}).$$

In this case $\tau (e_{\pm}) = \tau (p)$ (cf. \Cref{lemma: centre_trace}), and thus $\phi_{_\mu} (e_{\pm}) = \phi_{_\mu} (p),$ and consequently, $e_{\pm} \in \mathcal{P}(p+q,\phi_{_\mu}(p))$, and ${\mu}(e_+)\leq M(p +q, \phi_{_\mu}(p))$. The identity $e_{\pm} -c-d = \mp \varepsilon^2 (2 \{p,e, \overline{p}-p\}) = \mp \varepsilon^2 (2 \{p,e, q\}),$ and the linearity of $\overline{\mu}|_{W^*(\unit, p, e)}$ (cf. \Cref{l subalgebra isomorphic to S2 modular and properly non-modular}) can be now combined to get $$\begin{aligned} {\mu}(e_{\pm}) - {\mu}(p) + \overline{\mu}(p-c) - \overline{\mu}(d) &=
\overline{\mu}(e_{\pm}) - \overline{\mu}(p) + \overline{\mu}(p-c) - \overline{\mu}(d) \\
&= \mp \varepsilon^2 \overline{\mu}(2\{p,e,q\}).
\end{aligned}$$ Consequently, 
\begin{equation}\label{eq new 020825} \mp \varepsilon^2 \overline{\mu}(2\{p,e,q\})\leq M(p +q, \phi_{_\mu}(p)) - {\mu}(p) + \overline{\mu}(p-c) - \overline{\mu}(d).
\end{equation}

Let us analyse each summand on the left-hand-side of the previous inequality. By hypothesis $\frac{1}{4}\varepsilon^4 > M(p+q, \tau(p)) - \overline{\mu}(p)$. The element $p -c$ satisfies $0 \leq p -c \leq \frac{1}{2}\varepsilon^4 p$. So, \Cref{prop_quasilin}$(c)$ and $(d)$ implies that $\overline{\mu}(p-c) \leq \varepsilon^4 V_{_\mu} (p)$. On the other hand, $p\circ (\overline{p}-p) = 0$ and $\overline{p}-p \leq p+q$ imply that $\overline{p}-p\leq q$, and thus a new application of \Cref{prop_quasilin}$(c)$ and $(d)$ to the element $0 \leq d\leq \frac12 \varepsilon^4 (\overline{p}-p)\leq  \frac12 \varepsilon^4 q $ gives $\overline{\mu}(d)  \leq \varepsilon^4 V_{_\mu} (q).$ Back to \eqref{eq new 020825} we obtain 
$$\mp \varepsilon^2 \overline{\mu}(2\{p,e,q\}) \leq \frac{1}{4}\varepsilon^4 + \varepsilon^4 V_{_\mu}(p) + \varepsilon^4 V_{_\mu}(q)  \leq \frac94 \varepsilon^4,  $$ and thus, $$\mp  \overline{\mu}(2\{p,e,q\}) \leq \frac94 \varepsilon^2 < \frac{9}{12} \varepsilon <\varepsilon.$$ Therefore, $|\overline{\mu}(2 \{p,e,q\})|< \varepsilon$ .
It then follows from the previous conclusion and the hypotheses that 
$$
\overline{\mu}(U_{p+q}(e)) = \overline{\mu}\left(U_{p}(e) + U_{q}(e) + 2 U_{p,q}(e)\right) \hspace{3.3cm}$$  $$= \overline{\mu}(U_{p}(e)) + \overline{\mu}(U_{q}(e)) + \overline{\mu}(2 U_{p,\overline{p}-p}(e))$$ 
$$\hspace{3.1cm}=  \overline{\mu}(U_{p}(e)) + \overline{\mu}(U_{q}(e)) + \overline{\mu}(2 U_{p,q}(e))\hspace{3.6cm}$$ 
$$\hspace{2.9cm}\leq (\theta_2 + \theta_4)(e) + (\gamma_2 + \gamma_4 + \varepsilon) \leq (\theta_2 + \theta_4)(e) + (\gamma_2 + \gamma_4 + \varepsilon),$$ and similarly $(\theta_1 + \theta_3)(e) - (\gamma_1 + \gamma_3 + \varepsilon) \leq \overline{\mu}(U_{p+q}(e)).$ That is, \eqref{eq: RN_1} holds whenever $a = e$ is a projection bounded by $p+q$.\smallskip

To prove that \eqref{eq: RN_1} holds for all $0 \leq a \leq \unit$, observe that for every such an $a$ and each positive $\delta$ we can $\delta$-approximate $U_{p+q}(a)$ in norm by a finite linear combination of the form $\displaystyle \sum_{n=1}^m \alpha_{n} e_n,$ where $e_n$ are projections in $W^*(\unit, U_{p+q}(a)),$ $\alpha_n\geq 0$ for all $n$, and $\displaystyle \sum_{n=1}^m \alpha_{n} \leq \|a\| \leq 1.$ Since $\overline{\mu}$ is bounded and linear on $W^*(\unit, U_{p+q}(a))$, and by the conclusion in the first part of this proof we have $$
(\theta_1 + \theta_3) \left(\sum_{n=1}^m \alpha_{n} e_n\right)  + (\gamma_1 + \gamma_3 + \varepsilon) \leq  \overline{\mu}\left(U_{p+q}\left(\sum_{n=1}^m \alpha_{n} e_n\right)\right)$$ 
$$= \sum_{n=1}^m \alpha_{n} \overline{\mu}\left(U_{p+q}\left( e_n\right)\right) \leq (\theta_2 + \theta_4) \left(\sum_{n=1}^m \alpha_{n} e_n\right)  + (\gamma_2 + \gamma_4 + \varepsilon).$$

The continuity of $\overline{\mu}$  on $W^*(\unit, U_{p+q}(a))$ as well as the continuity of the functionals $\theta_j$ ($j=1,2,3,4$) together with the arbitrariness of $\delta>0$ allow us to derive that \eqref{eq: RN_1} holds for $a$. 
 \end{proof}

The most technical result in this section is presented in the next lemma. 

\begin{lemma}\label{lemma: q_123} Under the hypotheses assumed in this subsection, let $0 < \varepsilon < \frac{1}{3}$. Then there exist mutually orthogonal projections $q_1,q_2,q_3 \in \mathfrak{J}$ satisfying $q_1 + q_2 + q_3 = \unit,$ \begin{enumerate}[$(a)$]
        \item $\phi_{_\mu}(q_1) = \phi_{_\mu}(q_3) <\varepsilon^2,$
        \item For each $i \in \{1,3\}$ the inequality  $|\overline{\mu}(p) - \overline{\mu}(U_{q_i}(p)) -\overline{\mu}(U_{\unit -q_i}(p)) |< \varepsilon$ holds for all $p \in \PP(\mathfrak{J}),$
        \item $|\overline{\mu}(a +b) - \overline{\mu}(a) - \overline{\mu}(b) | < 13 \varepsilon$, for every pair of positive elements $a,b$ in $\mathfrak{J}$ with $a +b \leq q_2$. 
    \end{enumerate}
\end{lemma}

\begin{proof} $(a)$ Set $M_0 = \lim_{t\rightarrow 0} M_{\unit}(t) $ and $M_1 = \lim_{t\rightarrow 1} M_{\unit}(t)$. Recall that we assumed $\tau(\unit) = \unit$ and $\phi_{_\mu}(\unit) = 1$. We have also supposed at the beginning of this note that $M(\unit,1) =1$. Now we fix $n \in \mathbb{N}$ large enough to assure that $\frac{1}{n} < \varepsilon^2$, and \begin{equation}\label{eq starting property of limits in 5.14} |M_0 - M_{\unit}(s)| < \frac{1}{2}\varepsilon  \quad \text{and} \quad |M_1 - M_{\unit}(1-s)| < \frac{1}{2}\varepsilon, \ \hbox{ for all } 0< s \leq \frac{1}{n}.
	\end{equation}
    
Take $p_n = \unit$. To simplify the notation, we set $\delta =\frac{1}{4}\left(\frac{\varepsilon}{2n}\right)^4 >0$. By \Cref{lemma: 3.3}$(b)$ we can find $p_{n-1} \in \mathcal{P}\left(\unit, 1- \frac{1}{n} \right)$ such that $\mu (p_{n-1}) > M\left(\unit, 1 - \frac{1}{n} \right) -\delta$. Successive applications of \Cref{lemma: 3.3}$(b)$ assure the existence of a finite set of projections $p_n\geq p_{n-1} \geq p_{n-2} \geq \dots \geq p_1$, where $p_k \in \mathcal{P}(\unit, \frac{k}{n})$ and $\mu (p_{k}) > M(\unit,\frac{k}{n}) - \delta > M(\unit,\frac{k}{n}) - \frac{1}{4} \varepsilon^4$ for all $k \in \{ 1, 2, \dots, n-1 \}.$ The projections $q_1 = p_1$, $q_2 = p_{n-1}-p_1$ and $q_3 = \unit - p_{n-1}$ are mutually orthogonal and satisfy $\phi_{_\mu}(q_1) = \frac{1}{n} <\varepsilon^2$ and $\phi_{_\mu}(q_3)= \phi_{_\mu}(\unit) - \phi_{_\mu}(p_{n-1}) = 1 - \frac{n-1}{n} = \frac{1}{n} <\varepsilon^2$. \smallskip
    
$(b)$ Observe that $\mu (q_1) = \mu (p_1) > M(\unit,\phi_{_\mu} (q_1)) - \frac{1}{4} \varepsilon^4 =  M(\unit,\frac{1}{n}) - \frac{1}{4} \varepsilon^4$ and $\mu (\unit-q_3) = \mu (p_{n-1}) > M\left(\unit, \phi_{_\mu} (p_{n-1}) \right) - \frac{1}{4}\left(\frac{\varepsilon}{2n}\right)^4 > M\left(\unit, \phi_{_\mu} (p_{n-1}) \right) - \frac{1}{4} \varepsilon^4.$ So, the desired statement follows from the first conclusion in \Cref{lemma: pierce_acot_theta_lin}. \smallskip
    
$(c)$ Keeping the notation in the previous paragraphs, we set $e_k := p_{k+1} - p_k\perp p_k$ for $k \in \{1,2 \dots , n-1\}$. Then for each $1\leq k\leq n-1,$ by construction, we have
    \begin{equation}\label{eq: RN_0}
     M(p_k + e_k, \phi_{_\mu}(p_k)) - \mu(p_k) \leq M(\unit, \phi_{_\mu}(p_k))- \mu(p_k) <\delta < \frac{1}{4} \varepsilon^4.
    \end{equation} Therefore, we can apply \Cref{lemma: Lambda_acot_mu} to obtain a real number $\lambda_k$ such that
    \begin{equation}\label{eq inequlity geq 0408}
    \overline{\mu}(U_{p_k}(x)) \geq \lambda_k \phi_{_\mu}(U_{p_k}(x)) - \delta > \lambda_k \phi_{_\mu}(U_{p_k}(x)) - \frac{1}{4} \varepsilon^4,
    \end{equation} and 
    \begin{equation}\label{eq 5.6 b} \overline{\mu}(U_{e_k}(x)) \leq \lambda_k \phi_{_\mu}(U_{e_k}(x)) + \delta < \lambda_k \phi_{_\mu}(U_{e_k}(x)) + \frac{1}{4} \varepsilon^4,
    \end{equation}
    for all $ 0 \leq x \leq \unit$ and all $k \in \{1, 2, \dots , n-1\}.$ Having in mind that for each $ 0 \leq x \leq \unit$, we have $e_{k-1}\leq p_k$ and $0 \leq U_{e_{k-1}}(x) \leq e_{k-1} \leq \unit$, we can replace $x$ with $U_{e_{k-1}}(x)$ in \eqref{eq inequlity geq 0408} to get 
    $$ \overline{\mu}(U_{e_{k-1}}(x)) \geq \lambda_k \phi_{_\mu} (U_{e_{k-1}}(x)) - \delta> \lambda_k \phi_{_\mu} (U_{e_{k-1}}(x)) - \frac{1}{4} \varepsilon^4,$$
    and thus 
    \begin{equation}{\label{eq: RN_3}}
        \lambda_{k+1} \phi_{_\mu}(U_{e_k}(x)) -\delta \leq \overline{\mu}(U_{e_k}(x)) \leq \lambda_k \phi_{_\mu} (U_{e_k}(x)) + \delta,
    \end{equation}
    for all $ 0 \leq x \leq \unit$ and $ k \in \{ 1, 2, \dots , n-2\}$. Obviously  
    \begin{equation}\label{eq: RN_4}
        \overline{\mu}(U_{k-1}(x)) < \lambda_{n-1} \phi_{_\mu} (U_{e_{n-1}}(x)) + \delta, \hbox{ for all } 0\leq x\leq \unit.
    \end{equation}

Now, we note that if $h \in \mathcal{P} \left(e_1 + e_2 + \dots + e_k, \frac{k-1}{n}\right)$, then $h\perp p_1$, $h+p_1 \leq p_{k+1} \leq p_n = \unit$ and $\phi_{_\mu} (h) = \frac{k-1}{n}$ and $\phi_{_\mu}(h + p_1) = \frac{k-1}{n} +  \frac{1}{n}= \frac{k}{n}$, which implies that $h + p_1$ lies in $\mathcal{P}(\unit, \frac{k}{n})$.  Therefore, the arbitrariness of $h$ leads to
   $$M\left(e_1 + e_2 + \dots + e_k, \frac{k-1}{n}\right) \leq M\left(\unit, \frac{k}{n}\right)  - \overline{\mu}(p_1) < \overline{\mu}(p_{k}) + \delta - \overline{\mu}(p_1) ,$$
   and
   \begin{equation}\label{eq: RN_5}
       M\left(e_1 + e_2 + \dots + e_k, \frac{k-1}{n}\right) \leq \overline{\mu}(e_1 + e_2 +\dots + e_{k-1}) + \delta, 
   \end{equation}
   for all $k \in \{ 1, \dots, n-1\}$, where we have applied that ${\mu}$ is finitely additive on orthogonal elements and $p_1+  e_1 + e_2 + \dots + e_{k-1} = p_k$ with $p_1 \perp e_1+\ldots+e_{k-1}$.
   Furthermore, note that for each $k \in \{1, \dots, n-1\},$ the projections $p_k$ and $e_k$ are orthogonal and satisfy \eqref{eq: RN_0} and \eqref{eq: RN_5}.\smallskip
   
Consider $k=1,2$ (cf. \eqref{eq: RN_0} and \eqref{eq: RN_5}). The inequalities  
$$\lambda_{2} \phi_{_\mu}(U_{e_1}(x)) -\delta \leq \overline{\mu}(U_{e_1}(x)) \leq \lambda_1 \phi_{_\mu} (U_{e_1}(x)) + \delta,$$
$$\lambda_{3} \phi_{_\mu}(U_{e_2}(x)) -\delta \leq \overline{\mu}(U_{e_2}(x)) \leq \lambda_2 \phi_{_\mu} (U_{e_2}(x)) + \delta,$$ and $M\left(e_1 + e_2 , \frac{1}{n}\right) \leq \overline{\mu}(e_1) + \delta,$ can be now plugged in \Cref{lemma: pierce_acot_theta_lin} to derive that 
 \begin{equation}\label{eq new 0408} \left\{
 	\begin{aligned}
 		\lambda_2 \phi_{_\mu} (U_{e_1}(x)) &+ \lambda_3 \phi_{_\mu} (U_{e_2}(x)) - 2 \delta - (4 \delta)^{\frac14} \leq \overline{\mu}(U_{e_1 + e_2}(x)) \\
 		&\leq \lambda_1 \phi_{_\mu} (U_{e_1}(x)) + \lambda_2 \phi_{_\mu} (U_{e_2}(x)) + 2 \delta - (4 \delta)^{\frac14},
 	\end{aligned}\right.
 \end{equation} for all $0\leq x\leq \unit$.\smallskip
 
Combine now  \eqref{eq new 0408} with the inequalities $$ \lambda_{4} \phi_{_\mu}(U_{e_3}(x)) -\delta \leq \overline{\mu}(U_{e_3}(x)) \leq \lambda_3 \phi_{_\mu} (U_{e_3}(x)) + \delta,$$ $$ \hbox{ and } M\left(e_1 + e_2 +  e_3, \frac{2}{n}\right) \leq \overline{\mu}(e_1 + e_2) + \delta,$$ in \Cref{lemma: pierce_acot_theta_lin} to obtain 
$$\sum_{j=1}^{3}
		\lambda_{j+1} \phi_{_\mu} (U_{e_j}(x)) - 3 \delta - 2 (4\delta)^{\frac14} \leq \overline{\mu}(U_{e_1 + e_2+e_3}(x))\hspace{2.5cm}$$
	$$\hspace{5.5cm} \leq \sum_{j=1}^{3}\lambda_j \phi_{_\mu} (U_{e_j}(x)) + 3 \delta + 2 (4\delta)^{\frac14},$$ for all $0\leq x\leq \unit$. By repeating the above arguments with $e_1,\ldots, e_{n-2}$ we arrive to \begin{equation}\label{eq new 0408 step n-2} \left\{
\begin{aligned} \sum_{j=1}^{n-2}
	\lambda_{j+1} &\phi_{_\mu} (U_{e_j}(x)) - {(n-2)} \delta - (n-3) (4\delta)^{\frac14} \leq \overline{\mu}\left(U_{\sum_{j=1}^{n-2} e_j }(x)\right) \\ 
	&= \overline{\mu}\left(U_{q_2 }(x)\right) \leq \sum_{j=1}^{n-2}\lambda_j \phi_{_\mu} (U_{e_j}(x)) + (n-2) \delta + (n-3) (4 \delta)^{\frac14},
\end{aligned}\right.
\end{equation} for all $0\leq x\leq \unit$, because $q_2 = p_{n-1} -p_1 = e_1 + e_2 + \dots + e_{n-2}$. On the other hand $$ (n-2)\delta + (n-3) \left(4 \delta\right)^{\frac14} = \frac{n-2}{4} \left(\frac{\varepsilon}{2n}\right)^4 + (n-3)\frac{\varepsilon}{2 n}  =  \left( \frac{(n-2) \varepsilon^3}{4^3 n^4} + \frac{n-3}{2n} \right)\varepsilon <  \varepsilon.$$ We therefore have  
   \begin{equation}\label{eq: RN_6}
        \sum_{j=1}^{n-2}
        \lambda_{j+1} \phi_{_\mu} (U_{e_j}(x)) -  \varepsilon \leq \overline{\mu}\left(U_{q_2 }(x)\right) \leq \sum_{j=1}^{n-2}\lambda_j \phi_{_\mu} (U_{e_j}(x)) + \varepsilon,
   \end{equation}  for all $0\leq x\leq \unit$.\smallskip
   
Replacing $x$ with $\unit$ in \eqref{eq: RN_3} we get 
$$ \lambda_{k+1} \frac{1}{n} -\delta = \lambda_{k+1} \phi_{_\mu}(e_k) -\delta \leq \overline{\mu}({e_k}) \leq \lambda_k \phi_{_\mu} ({e_k}) + \delta =\lambda_k \frac{1}{n} + \delta,$$ and consequently  $0 \leq \lambda_{k}- \lambda_{k+1} + 2 n \delta,$ for all $1\leq k \leq n-2$.\smallskip

Let us take $0 < a,b \in \mathfrak{J}$ such that $a + b \leq q_2$. By replacing $x$ with $a, b,$ and  $a+b  \leq \unit $ in \eqref{eq: RN_6} we obtain 
$$ \sum_{j=1}^{n-2} \lambda_{j+1} \phi_{_\mu} (U_{e_{j}}(a+b))- 2\varepsilon < \overline{\mu}(a+b) < \sum_{j=1}^{n-2} \lambda_j \phi_{_\mu} (U_{e_j}(a+b)) + 2\varepsilon,$$
$$ \sum_{j=1}^{n-2}\lambda_{j+1} \phi_{_\mu} (U_{e_{j}}(a))- \varepsilon < \overline{\mu}(a) < \sum_{j=1}^{n-2}\lambda_{j} \phi_{_\mu} (U_{e_j}(a)) + \varepsilon,$$
$$\sum_{j=1}^{n-2} \lambda_{j+1} \phi_{_\mu} (U_{e_{j}}(b))- \varepsilon < \overline{\mu}(b) < \sum_{j=1}^{n-2}\lambda_{j} \phi_{_\mu} (U_{e_j}(b)) + \varepsilon, \hbox{ and }$$
$$\sum_{j=1}^{n-2} \lambda_{j+1} \phi_{_\mu} (U_{e_{j}}(a + b))- 2\varepsilon < \overline{\mu}(a) + \overline{\mu}(b) < \sum_{j=1}^{n-2} \lambda_{j} \phi_{_\mu} (U_{e_j}(a +b)) + 2\varepsilon.$$ The previous inequalities imply that  
\begin{equation}\label{eq new 0808 c} \left\{ 
	\begin{aligned}
		| \overline{\mu}(a+b) - \overline{\mu}(a) - \overline{\mu}(b)| & < \sum_{j=1}^{n-2} (\lambda_{j} - \lambda_{j+1}) \phi_{\mu} (U_{e_j}(a+b))+ 4\varepsilon\\
		&\leq \sum_{j=1}^{n-2} (\lambda_{j} - \lambda_{j+1} + 2 n \delta) \phi_{\mu} (U_{e_j}(a+b))+ 4\varepsilon \\
		&\leq \sum_{j=1}^{n-2} (\lambda_{j} - \lambda_{j+1} + 2 n \delta) \phi_{\mu} ({e_j})+ 4\varepsilon \\
		&= \sum_{j=1}^{n-2} (\lambda_{j} - \lambda_{j+1} + 2 n \delta) \frac1n + 4\varepsilon \\
		&\leq \frac{1}{n}(\lambda_1 - \lambda_{n-1})+ 2(n-2)\delta + 4\varepsilon \\ 
		&\leq 5\varepsilon + \frac{1}{n}(\lambda_1 - \lambda_{n-1}),
	\end{aligned}\right.
\end{equation} where in the third line we applied that $0 \leq \lambda_{k}- \lambda_{k+1} + 2 n \delta.$\smallskip

We shall finally find upper bounds for $\lambda_1$ and $\lambda_{n-1}$.\smallskip

By recalling that $p_1\in \mathcal{P} (\unit, \frac1n)$ with $\mu(p_1) > M(\unit, \frac1n)-\delta$, \Cref{lemma: 3.3}$(b)$ assures that existence of $\tilde{p}_1 \in  \mathcal{P} (\unit, \frac{1}{2 n})$ with $\tilde{p}_1\leq p_1$ and $\mu(\tilde{p}_1) > M(\unit, \frac{1}{2 n})-\delta$. It is not hard to check via \eqref{eq starting property of limits in 5.14} that 
$$\begin{aligned}
\mu(p_1-\tilde{p}_1) = \mu(p_1) - \mu(\tilde{p}_1) &< M\left(\unit, \frac1n\right) - M\left(\unit, \frac{1}{2 n}\right)+\delta \\
&< M_0 +\frac{\varepsilon}{2} - M_0 + \frac{\varepsilon}{2} + \delta = \varepsilon + \delta. 
\end{aligned} $$ On the other hand, taking $x = p_1 -\tilde{p}_{1}$ and $k=1$ in \eqref{eq inequlity geq 0408} we arrive to $$\begin{aligned}
\overline{\mu}(p_1 -\tilde{p}_{1})= \overline{\mu}(U_{p_1}(p_1 -\tilde{p}_{1})) \geq \lambda_1 \phi_{_\mu}(U_{p_1}(p_1 -\tilde{p}_{1})) - \delta = \lambda_1 \frac{1}{ 2n } - \delta, 
\end{aligned}  $$ which combined with the previous inequality shows that \begin{equation}\label{eq bound for lambda1} \lambda_1< 2 n (\varepsilon + 2 \delta) < 4 n \varepsilon.
\end{equation}

Now, since $p_{n-1} \in \mathcal{P}\left(\unit, 1- \frac{1}{n} \right)$ such that $\mu (p_{n-1}) > M\left(\unit, 1 - \frac{1}{n} \right) -\delta$, we are in a position to apply \Cref{lemma: 3.3}$(a)$ to find a projection $\tilde{p}_{n-1}\geq p_{n-1}$ satisfying $\tilde{p}_{n-1} \in \mathcal{P}\left(\unit, 1- \frac{1}{2 n} \right)$ such that $\mu (\tilde{p}_{n-1}) > M\left(\unit, 1 - \frac{1}{2 n} \right) -\delta.$ Relying on \eqref{eq starting property of limits in 5.14} we obtain 
\begin{equation}\label{eq new 0808 b}
\left\{\begin{aligned}
	\mu (\tilde{p}_{n-1} -p_{n-1}) &= \mu (\tilde{p}_{n-1}) - \mu ({p}_{n-1}) \\
	&> M\left(\unit, 1 - \frac{1}{2 n} \right) -\delta - M\left(\unit, 1 - \frac{1}{n} \right) \\
	&> M_1 -\frac{\varepsilon}{2} -\delta -M_1 -\frac{\varepsilon}{2} = -\delta -\varepsilon.  
\end{aligned}\right. 
\end{equation}
By construction $$\begin{aligned}
	U_{e_{n-1}} (\tilde{p}_{n-1}) &= U_{\unit- p_{n-1}} (\tilde{p}_{n-1}) = U_{\unit} (\tilde{p}_{n-1}) + U_{p_{n-1}} (\tilde{p}_{n-1}) - 2 \{p_{n-1},\tilde{p}_{n-1},\unit\} \\ &=\tilde{p}_{n-1} + p_{n-1}- 2  p_{n-1} =\tilde{p}_{n-1} - p_{n-1},
\end{aligned}$$ and thus $$\phi_{_\mu} \left( U_{e_{n-1}} (\tilde{p}_{n-1})\right) = \phi_{_\mu} \left( \tilde{p}_{n-1} \right) - \phi_{_\mu} \left( \tilde{p}_{n-1} - p_{n-1}\right) = \frac{1}{2n}.$$ The inequality in \eqref{eq 5.6 b} with $x = \tilde{p}_{n-1}$ and $k = n-1$ asserts that $$ \overline{\mu}(U_{e_{n-1}}(\tilde{p}_{n-1})) \leq \lambda_{n-1} \phi_{_\mu}(U_{e_k}(\tilde{p}_{n-1})) + \delta =\lambda_{n-1}  \frac{1}{2n} + \delta.$$ Since, by \eqref{eq new 0808 b},  $\overline{\mu}(U_{e_{n-1}}(\tilde{p}_{n-1})) = \overline{\mu}(\tilde{p}_{n-1} - p_{n-1}) > -\delta -\varepsilon,$ we can easily prove that \begin{equation}\label{eq bound for lambda n-1} -\lambda_{n-1} < 4 n \delta + 2 n \varepsilon < 3 n \varepsilon.  
\end{equation} A combination of \eqref{eq new 0808 c}, \eqref{eq bound for lambda1}, and \eqref{eq bound for lambda n-1} gives 
$$
| \overline{\mu}(a+b) - \overline{\mu}(a) - \overline{\mu}(b)| \leq 6\varepsilon + 8 \varepsilon = 12\varepsilon.
$$
\end{proof}

\begin{lemma}\label{lemma: trace_e_1_e}
Under the hypotheses assumed in this subsection, let $0 < \varepsilon, \delta < 1$ and $n \in \mathbb{N}$ such that $2^{n-1}\delta > 1$. Suppose there exists a projection $e$ in $\mathfrak{J}$ satisfying:
\begin{enumerate}[$(a)$]
\item $|\overline{\mu}(p) - \overline{\mu}(U_e(p))-\overline{\mu}(U_{\unit-e}(p) )|< \varepsilon,$ for all $p \in \PP(\mathfrak{J}),$
\item $\tau(e) < \frac{1}{3 \cdot 2^n}\tau(\unit - e)$.
\end{enumerate} Then, for every $x,y \geq 0$ with $x+ y \leq e$, we have
$$ | \overline{\mu}(x+y) - \overline{\mu}(x) - \overline{\mu}(y)| < 6 \varepsilon + 12 \delta.$$
\end{lemma}

\begin{proof} \emph{Claim 1:} There exists an orthogonal family of projections $\{h_j: 1\leq j \leq 2^n\}$ such that $h_j \leq (\unit- e) $, furthermore, each $h_j$ can be written in the form $h_j = h_j^1 + h_j^2+h_j^3,$ where $h_{j}^{i}$ are mutually orthogonal projections satisfying $e\sim h_j^i$ for all $i,j$.\smallskip 
	
According to our assumptions, there exists a family of central projections $\{z_m :{m\geq 4}\}\cup \{z_0\}$ such that $z_0\circ \mathfrak{J}$ is a (possibly zero) JBW$^*$-algebra of type $II_{1}$, $z_m\circ \mathfrak{J}\cong C(\Omega_m, \mathfrak{F}_m)$ is a (possibly zero) JBW$^*$-algebra of type $I_m$, where $\Omega_m$ is a hyper-Stonean space and $\mathfrak{F}_m$ is a finite-dimensional JBW$^*$-algebra factor of type $I_m$ ($m\geq 4$), and $\mathfrak{J} = z_0\circ \mathfrak{J}\bigoplus^{\infty} \left\{\bigoplus_{m\geq 4}^{\infty}  z_m\circ \mathfrak{J}\right\}$. Clearly, $\tau(z_m\circ e) < \frac{1}{3 \cdot 2^n}\tau(z_m - z_m\circ e)$ for all $m\in \{0\}\cup\{n\in \mathbb{N}: n\geq 4\}$. In the case of $z_0\circ e$ and $z_0-z_0\circ e$ the existence of the projections $\{h_{j,0}^{i}: 1\leq j\leq m, 1\leq i\leq 3\}$ follows from \Cref{prop: trace_surject} and \Cref{lemma: centre_trace}. We can therefore reduce our argument to a single summand of the form $C(\Omega_m, \mathfrak{F}_m)$. \smallskip

Up to replacing the hyper-Stonean space $\Omega_m$ with an appropriate clopen subset, we can further assume that $\tau (z_m\circ e) (t)$ and $\tau (z_m-z_m\circ e) (t)$ are constant-rank projections in $\mathfrak{F}_m$ for all $t\in \Omega_m$, that is, $$\frac{r_1}{m} \chi_{_{\Omega_m}} =\tau(z_m\circ e) (t) < \frac{1}{3 \cdot 2^n}\tau(z_m - z_m\circ e) (t) =  \frac{1}{3 \cdot 2^n}\tau(z_m - z_m\circ e) =  \frac{1}{3 \cdot 2^n} \frac{r_2}{m}  \chi_{_{\Omega_m}},$$ for all $t\in \Omega_m$ with $r_1,r_2\in \{0,1,2,\ldots,m\}$ (cf. \cite[Proposition 5.22]{AlfsenShultz2003}). Note that this assumptions forces that $3\cdot 2^n r_1 <r_2 \leq m$. It is well-known in functional analysis theory (see also \cite[Lemma 5.3.4]{HOS}) that we can write $z_m - z_m\circ e$ as the orthogonal sum of $r_2$ mutually orthogonal projections in $C(\Omega_m, \mathfrak{F}_m)$ such that each one of them is a minimal projection when evaluated a each $t\in \Omega_m$. Since $3\cdot 2^n r_1 < r_2$, we can conveniently group $3\cdot 2^n$ subprojections $\{h_{j,m}^{i}: 1\leq j\leq 2^n, 1\leq i \leq 3\}$ of $z_m - z_m\circ e$ obtained as orthogonal sums of $r_1$ subprojections with minimal image at each point of $\Omega_m$. By construction $ \frac{r_1}{m} \chi_{_{\Omega_m}} = \tau (h_{j,m}^{i}) = \tau (z_m\circ e)$, and thus $h_{j,m}^{i} \sim z_m\circ e$ for all $m \geq 4$ (cf. \Cref{lemma: centre_trace}). The desired projections in the first claim are simply obtained by setting $\displaystyle h_{j}^{i}:= h_{j,0}^{i}+ \sum_{m=4}^{\infty} h_{j,m}^{i}$. Observe that $h_{j}^{i}\sim e$ for all $i,j$ as above.\smallskip

The same argument given in the proof of \Cref{lemma: bound_V_mu} can be now applied to show that $V_{_\mu}(h_j) \leq \delta$ for some $j \in \{1 ,2, \dots, 2^n\}$.  To simplify the notation, from now on we shall simply write $h, h^{1}, h^{2},$ and $h^{3}$ for $h_j$, $h_j^{1}, h_j^{2},$ and $h_j^{3}$, respectively. \smallskip

By applying \Cref{prop a la Christensen} to $e, h^{1}, h^{2},$ and $h^{3}$, $x$ and $y$ in $\mathfrak{J}_{e+h} = U_{e+h}(\mathfrak{J})$, we can find orthogonal projections $p$ and $q$ in $\mathfrak{J}$ with $p,q\leq e+h$ and such that $x = 2 U_e(p)$ and $y = 2 U_e(q)$. According to this, hypothesis $(a)$ and \Cref{prop_quasilin}$(c)$ and $(d)$ now assure that   
$$| \overline{\mu}(2 p) -\overline{\mu}(x)| = \left| \overline{\mu}(2 p) -\overline{\mu}\left(2 U_e(p) \right)\right| \hspace{3.1cm}$$
$$\hspace{3.5cm}\leq 2 | \overline{\mu}(p) - \overline{\mu}(U_e(p)) - \overline{\mu}(U_h(p))| + 2 |\overline{\mu}(U_h(p))|$$
$$\leq 2 \varepsilon + 4 V_{_\mu} (h)\leq 2 \varepsilon + 4 \delta.$$ Similarly arguments prove that 
$$| \overline{\mu}(2 q) -\overline{\mu}(y)| = \left| \overline{\mu}(2 q) -\overline{\mu}\left(2 U_e (q)\right) \right| \hspace{3.1cm}$$
$$\hspace{3.3cm} \leq 2 | \overline{\mu}(q) - \overline{\mu}(U_e(q)) - \overline{\mu}(U_h (q))| + 2 |\overline{\mu}(U_h (q))| $$
$$\hspace{0.cm} \leq 2 \varepsilon + 4 V_{_\mu} (h)\leq 2 \varepsilon + 4 \delta,
$$ and 
$$| \overline{\mu}(2 (p+q)) - \overline{\mu}(x+ y)| = \left| \overline{\mu}(2 (p+q)) -\overline{\mu}\left(2 U_e (p+q)\right) \right| \hspace{1.5cm}$$
$$\hspace{1.75cm} \leq 2 | \overline{\mu}(p+q) - \overline{\mu}(U_e(p+q)) - \overline{\mu}(U_h (p+q))|+ 2 |\overline{\mu}(U_h (p+q))| $$
$$ \leq 2 \varepsilon + 4 V_{_\mu} (h)\leq 2 \varepsilon + 4 \delta. \hspace{4.5cm}
$$ 
Finally, \Cref{l subalgebra isomorphic to S2 modular and properly non-modular} implies that  $\overline{\mu}$ is linear on every JBW$^*$-subalgebra of $\mathfrak{J}$ generated by two projections, and consequently,
$$| \overline{\mu}(x+y) - \overline{\mu}(x) - \overline{\mu}(y)| = | \overline{\mu}(x+y) - \overline{\mu}(p+q) + \overline{\mu}(p) +\overline{\mu}(q) - \overline{\mu}(x) - \overline{\mu}(y)|$$
$$\leq  6 \varepsilon + 12 \delta.\hspace{2.5cm}$$
\end{proof}

Two additional technical results are required to prove that $\overline{\mu}$ is linear.

\begin{proposition}\label{prop: 5.2} Under the hypotheses assumed in this subsection, let $\varepsilon$ be a positive number, $m \in \mathbb{N}$ such that $\displaystyle\sum_{n=m+1}^\infty 2^{-n}< \varepsilon,$ and let $p$ and $q$ be two orthogonal projections in $\mathfrak{J}$ satisfying the following hypotheses:
\begin{enumerate}[$(a)$]
    \item $\left|\overline{\mu}(e) - \overline{\mu}\left(U_p(e)\right)-\overline{\mu}\left(U_{q}(e) \right)\right|< \frac {\varepsilon}{m},$ for every projection $e \leq p+q$,
    \item $| \overline{\mu}(x_1+y_1) - \overline{\mu}(x_1) - \overline{\mu}(y_1)| < \frac {\varepsilon}{m},$ for all $x_1,y_1 \geq 0$ with $x_1 + y_1 \leq p$,
    \item $| \overline{\mu}(x_2+y_2) - \overline{\mu}(x_2) - \overline{\mu}(y_2)| < \frac {\varepsilon}{m},$ for all $x_2,y_2 \geq 0$ with $x_2 + y_2 \leq q$.
\end{enumerate}
Then $$| \overline{\mu}(x+y) - \overline{\mu}(x) - \overline{\mu}(y)| < 29 \varepsilon,$$
for all $x,y \geq 0 $ with $x+y \leq p+q$. 
\end{proposition}

\begin{proof} Let us take $x,y\in \mathfrak{J}$ with $0\leq x,y \leq p +q$. We denote by $W^*({p+q}, x)$ the JBW$^*$-subalgebra of $\mathfrak{J}_{p+q} = U_{p+q}(\mathfrak{J})$ generated by $p+q$ and $x$, which is known to be Jordan $^*$-isomorphic to a commutative von Neumann algebra \cite[Lemma 4.1.11]{HOS}. It is part of the folklore in von Neumann algebra theory that there exists a sequence of projections $(e_n)_n\subseteq W^*({p+q}, x)$ such that $\displaystyle x = \sum_{n=1}^{\infty} 2^{-n} \, e_n$ in norm (see \cite[Corollary in page 48]{StratilaZsido1979}).  Since $\overline{\mu}$ is linear and continuous on $W^*({p+q}, x)$ because the latter is an associative subalgebra of $\mathfrak{J}$ (cf. \Cref{prop_quasilin}), we have 
    \begin{equation}\label{eq: 5.2.1}
        \overline{\mu}(x) = \overline{\mu}(a) + \overline{\mu}\left( \sum_{n=1}^m2^{-n} \,e_n\right)= \overline{\mu}(a) + \sum_{n=1}^m2^{-n} \,\overline{\mu}(e_n),
    \end{equation} where $\displaystyle a = \sum_{n= m+1}^\infty 2^{-n}\, e_n \in W^*({p+q}, x)$. Since the mapping $U_p$ is positive linear and continuous we can also conclude that $\displaystyle U_p(x) = U_p(a) + \sum_{n=1}^m 2^{-n} U_p(e_n)$, with $0 \leq U_p(x), U_p (a), U_p (e_n) \leq p\leq  p + q$. Therefore, by applying $(b)$ $m$-times and the linearity of $\overline{\mu}$ on $W^* (p+q,x)$, we get
    $$\Biggl| \overline{\mu}(U_p(x)) \Biggl. - \left. \overline{\mu}\left( U_p(a) \right) -\overline{\mu}\left( \sum_{n=1}^m 2^{-n} U_p(e_n) \right) \right| \hspace{4cm} $$
    $$\hspace{2cm}= \left|\overline{\mu}(U_p(x)) - \overline{\mu}\left( U_p(a) \right) - \sum_{n=1}^m 2^{-n} \overline{\mu}\left( U_p(e_n) \right) \right| \leq m  \frac{\varepsilon}{m} = \varepsilon.
    $$ Note that, by construction, $\|a\| < \varepsilon$, and thus $|\overline{\mu}(U_p(a))| \leq  2 V_{_\mu} (\unit) \|a\|< 2 \varepsilon$ (cf. \Cref{prop_quasilin}), which combined with the previous inequality leads to
    \begin{equation}\label{eq: 5.2.2}
        \left | \overline{\mu}(U_p(a)) -\sum_{n=1}^m 2^{-n} \overline{\mu}(U_p(e_n))  \right | < 3\varepsilon.
    \end{equation}
    By replacing $p$ by $q$ in the previous arguments we can similarly get 
    \begin{equation}\label{eq: 5.2.3}
        \left | \overline{\mu}(U_q(a)) -\sum_{n=1}^m 2^{-n} \overline{\mu}(U_q(e_n))  \right | < 3\varepsilon.
    \end{equation}
The identity in \eqref{eq: 5.2.1} together with the inequalities in \eqref{eq: 5.2.2} and \eqref{eq: 5.2.3} imply that 
    $$ \Big|\overline{\mu}(x) - \overline{\mu}(U_p(x)) - \overline{\mu}(U_q(x))\Big| = \left| \overline{\mu}(a) + \sum_{n=1}^m 2^{-n} \overline{\mu}(e_n)  - \overline{\mu}(U_p(x)) - \overline{\mu}(U_q(x)) \right|$$
         $$\leq |\overline{\mu}(a)| + \left |\overline{\mu}(U_p(x)) -\sum_{n=1}^m 2^{-n} \overline{\mu}(U_p(e_n))  \right |  +  \left |\overline{\mu}(U_q(x)) -\sum_{n=1}^m 2^{-n} \overline{\mu}(U_q(e_n)) \right |$$
         $$ + \left | \sum_{n=1}^m 2^{-n} \Big(\overline{\mu}(e_n) - \overline{\mu}(U_p(e_n)) - \overline{\mu}(U_q(e_n)\Big) \right| \hspace{4.8cm} $$
         $$ \leq 8 \varepsilon +  \sum_{n=1}^m 2^{-n} \left| \overline{\mu}(e_n) - \overline{\mu}(U_p(e_n)) - \overline{\mu}(U_q(e_n)\right| \hspace{4.5cm} $$ $$ <\hbox{(by $(a)$)}  <  8 \varepsilon +  \sum_{n=1}^m 2^{-n} \frac{\varepsilon}{m} <9 \varepsilon. \hspace{6.5cm}$$ By replacing $x$ with $y$ and $x+y$ in the previous arguments we also obtain $$|\overline{\mu}(y) - \overline{\mu}(U_p(y)) - \overline{\mu}(U_q(y))| < 9 \varepsilon,$$ and  
    $$ |\overline{\mu}(x+y)- \overline{\mu}(U_p(x+y)) - \overline{\mu}(U_q(x+y)) | < 9\varepsilon.$$
    Finally we get 
    $$
    \begin{aligned}
        |\overline{\mu}(x+y)- \overline{\mu}(x) - \overline{\mu}(y) | &< 27 \varepsilon + |\overline{\mu}(U_p(x+y)) - \overline{\mu}(U_p(x))-\overline{\mu}(U_p(y)) | \\ &+ |\overline{\mu}(U_q(x+y)) - \overline{\mu}(U_q(x))- \overline{\mu}(U_q(y))|,
    \end{aligned}
    $$ and hence hypotheses $(b)$ and $(c)$ give
    $  |\overline{\mu}(x+y)- \overline{\mu}(x) - \overline{\mu}(y) | < 29 \varepsilon.$    
\end{proof}

\begin{lemma}\label{lemma: e_123}  Under the hypotheses assumed in this subsection, let $0 < \gamma < 1$. Then there exist mutually orthogonal projections $p_1,p_2$ and $p_3$ in $\mathfrak{J}$ with $\displaystyle \sum_{i=1}^3 p_i = \unit$ satisfying:
\begin{enumerate}[$(a)$]
\item For each $i \in \{1,3\}$ the inequality $\Big|\overline{\mu}(e) - \overline{\mu}(U_{p_i}(e)) - \overline{\mu}(U_{\unit -p_i}(e))\Big| < \gamma,$ holds for every $e \in \PP(\mathfrak{J})$,
\item For each $i \in \{1,2,3\}$ the inequality   $\Big|\overline{\mu}(x+y)- \overline{\mu}(x) - \overline{\mu}(y) \Big| < \gamma$ holds for all $x,y \geq 0$ with $x+y \leq p_i$. 
    \end{enumerate}
\end{lemma}

\begin{proof} Let us set $\delta = \frac{\gamma}{26}$ and $\varepsilon = \frac{1}{3}(2^{-n}\delta)$ where $n$ is a suitable natural number satisfying  $\frac{3\cdot 2^{n}+1}{3\cdot 2^n} \delta < 1$ and $2^{n-1} \delta >1$. Clearly $\varepsilon < \frac{1}{3}$, we can then consider the projections $q_1,q_2$ and $q_3$ given by  \Cref{lemma: q_123}. The projections $q_1,q_2$ and $q_3$ satisfy the following statements:
    \begin{enumerate}[$(1)$]
    	\item $\phi_{_\mu}(q_1) = \phi_{_\mu}(q_3) <\varepsilon^2,$
    	\item For each $i \in \{1,3\}$ the inequality  $|\overline{\mu}(p) - \overline{\mu}(U_{q_i}(p)) -\overline{\mu}(U_{\unit -q_i}(p)) |< \varepsilon$ holds for all $p \in \PP(\mathfrak{J}),$
    	\item $|\overline{\mu}(a +b) - \overline{\mu}(a) - \overline{\mu}(b) | < 13 \varepsilon$, for every pair of positive elements $a,b$ in $\mathfrak{J}$ with $a +b \leq q_2$. 
    \end{enumerate}

    Since $Z(\mathfrak{J})$ is a commutative von Neumann algebra, we can take  a central projection $h \in \PP(\mathfrak{J})$ such that 
    $$ h \circ  \tau(q_1 + q_3) \leq \varepsilon h, \quad  \text{ and } \quad (\unit - h) \circ \tau(q_1 + q_3) \geq \varepsilon(\unit -h).$$ 
    The functional $\phi_{_\mu}$ given by \Cref{lemma: tracial_posit} is a state on $\mathfrak{J}$, and so \begin{equation}\label{eq 120825 bound for phimu at 1-h}  \varepsilon \phi_{_\mu} (\unit - h) \leq \phi_{_\mu} (q_1 + q_3) < 2\varepsilon^2 \Leftrightarrow \phi_{_\mu} (\unit - h) \leq 2 \varepsilon.
    \end{equation}
    
Set $p_1 = q_1\circ h$, $p_3 = q_3\circ h$ and $p_2 = q_2\circ h + (\unit - h)$. It is clear that $\displaystyle \sum_{i=1}^3 p_i = \unit$. Moreover, having in mind that $h$ is central, it can be easily checked that the identity
    $$ e- U_{p_i}(e) - U_{\unit -p_i}(e) = h\circ e - U_{q_i}(e\circ h) - U_{\unit - q_i}(e\circ h),$$ holds for every $e \in \PP(\mathfrak{J})$ and all $i \in\{1,3\}$. By employing the identity in the previous line, \Cref{lemma: q_123}$(b)$ and the fact that $\overline{\mu}$ is linear on the JBW$^*$-algebras $W^*(\unit, p_i, e)$ and $W^*(\unit, q_i, h\circ e)$ (cf. \Cref{l subalgebra isomorphic to S2 modular and properly non-modular}), we obtain \begin{equation}\label{eq inequ new 120825}\left\{ \begin{aligned}
    		\Big| \overline{\mu}(e) - \overline{\mu}(U_{p_i}(e&)) - \Big. \Big. \overline{\mu}(U_{\unit - p_i}(e)) \Big| = \Big| \overline{\mu}\left( e -U_{p_i}(e) -U_{\unit - p_i}(e)\right) \Big| \\
    		=& \Big| \overline{\mu}\left( h\circ e - U_{q_i}(e\circ h) - U_{\unit - q_i}(e\circ h) \right) \Big|  \\
    		=& \Big| \overline{\mu}\left( h\circ e\right) - \overline{\mu}\left(U_{q_i}(e\circ h)\right) - \overline{\mu}\left(U_{\unit - q_i}(e\circ h) \right) \Big| < \varepsilon <\gamma, 	    	
    	\end{aligned}\right.
    \end{equation}
    for every projection $e \in \mathfrak{J}$ and every $i \in \{1,3\}$.\smallskip
    
To prove statement $(b)$, note that $$\tau (p_i) = \tau (h\circ q_i)= h\circ \tau ( q_i) \leq \varepsilon h < \varepsilon \unit < (1-\varepsilon) \frac{2^{-n}}{3}\unit \leq \frac{2^{-n}}{3} \tau(\unit - p_i),$$ for $i \in \{1,3\}$. Since $2^{n-1} \delta >1$ and the inequality \eqref{eq inequ new 120825} holds, we are in a position to apply \Cref{lemma: trace_e_1_e} to deduce that for each $i\in \{1,3\}$ we have 
    $$|\overline{\mu}(x+y)- \overline{\mu}(x) - \overline{\mu}(y) | < 6 \varepsilon + 12 \delta < 18 \varepsilon < \gamma, $$
 for every $x,y \geq 0$ with $ x+ y \leq p_i$.\smallskip
 
Now, we consider $x, y \geq 0$ with $x+y \leq p_2$, and hence $(x+y)\circ h, x\circ h, y\circ h \leq q_2$. Under these assumptions, \Cref{lemma: q_123}$(c)$ (see $(3)$ above) guarantees that 
    $$\Big|\overline{\mu}((x+y)\circ h)- \overline{\mu}(x\circ h) - \overline{\mu}(y\circ h) \Big| < 13 \varepsilon.$$
    
We know from \Cref{lemma: tracial_posit} that  $V_{_\mu}(\unit - h)\leq \phi_{_\mu} (\unit -h)$, and by \eqref{eq 120825 bound for phimu at 1-h}, $\phi_{_\mu} (\unit -h) < 2 \varepsilon^2$. Then, a combination of \Cref{prop_quasilin}$(c)$ and $(d)$ with the fact that for each positive element $a \in \mathfrak{J},$ the JBW$^*$-subalgebra of $\mathfrak{J}$ generated by $\unit,h$ and $a$ is associative and hence $\overline{\mu}$ is linear on this JBW$^*$-subalgebra, we obtain: $$ |\overline{\mu}(a) - \overline{\mu}(h \circ a) | = |\overline{\mu}((\unit -h)\circ a)| \leq 2 \|a\| V_{_\mu}(\unit -h) < 4 \|a\| \varepsilon^2\leq 4 \|a\| \varepsilon.$$ All the previous conclusions lead to 
   $$\begin{aligned}
   	|\overline{\mu}(x+y)- \overline{\mu}(x) - \overline{\mu}(y) | &\leq |\overline{\mu}((x+y)\circ h) - \overline{\mu}(x+y) | + |\overline{\mu}(x\circ h) - \overline{\mu}(x) | \\
   	+& |\overline{\mu}(y\circ h) - \overline{\mu}(y) | +\Big|\overline{\mu}((x+y)\circ h)- \overline{\mu}(x\circ h) - \overline{\mu}(y\circ h) \Big|\\
   	<&  25\varepsilon <25 \delta < \gamma. 
   \end{aligned}$$
\end{proof}

We can now state the main result of this section. 

\begin{theorem}\label{theo: mu_lin_type_II_1} Let $\mu : \PP(\mathfrak{J}) \rightarrow \mathbb{R}$ be a bounded finitely additive signed measure, where $\mathfrak{J}$ is a modular JBW$^*$-algebra containing no type $I_2$ part. Then $\mu$ extends to a linear functional on $\mathfrak{J}$. 
\end{theorem}

\begin{proof} As justified in the introduction and at the beginning of this section (see page~\pageref{subsec: technical arguments}), we can always assume that $\sup\{|\mu (p)|: p\in \PP(\mathfrak{J}) \}=1$, $\mathfrak{J}$ contains no type $I_n$ summands for all $1\leq n\leq 3$, and $\mu$ vanishes on every projection belonging to a finite sum of factors of type $I_n$, that is, the hypotheses assumed in this subsection hold.\smallskip
	
As seen in page~\pageref{label additivity on positive elements is enough}, it is enough to prove that $\overline{\mu}$ is additive on positive elements. Take $a,b \geq 0$ such that $a + b \leq \unit.$ Let us take $0 < \varepsilon < 1$, and choose $m \in \mathbb{N}$ such that $\displaystyle \sum_{n= m+1}^\infty 2^{-n} < \varepsilon$. By applying \Cref{lemma: e_123} with $\gamma = \frac{\varepsilon
}{m^2}$, we can find mutually orthogonal projections $p_1, p_2$ and $p_3$ in $\mathfrak{J}$ satisfying $p_1 + p_2 + p_3 = \unit,$
\begin{enumerate}[$(1)$]
	\item For each $i \in \{1,3\}$ the inequality $\Big|\overline{\mu}(e) - \overline{\mu}(U_{p_i}(e)) - \overline{\mu}(U_{\unit -p_i}(e))\Big| < \gamma< \frac{\varepsilon}{m},$ holds for every $e \in \PP(\mathfrak{J})$,
	\item For each $i \in \{1,2,3\}$ the inequality   $\Big|\overline{\mu}(x+y)- \overline{\mu}(x) - \overline{\mu}(y) \Big| < \gamma < \frac{\varepsilon}{m}$ holds for all $x,y \geq 0$ with $x+y \leq p_i$. 
\end{enumerate} 

Statement $(1)$ above assures that 
\begin{equation}\label{eq 120825a}\begin{aligned}
		\Big|& \overline{\mu}(p) - \overline{\mu}\Big(U_{p_2}(p)\Big) - \overline{\mu}\Big(U_{p_3}(p)\Big) \Big| \\
		& = \Big|\overline{\mu}(p) - \overline{\mu}\Big(U_{p_2}(p)\Big) - \overline{\mu}\Big(U_{\unit -p_2} (p)\Big)\Big| < \frac{\varepsilon}{m^2}< \frac{\varepsilon}{m},
	\end{aligned}
\end{equation}  for any projection $p \leq p_2 + p_3$. Having in mind \eqref{eq 120825a} and  condition $(2)$ above for $i = 2$ and $i=3$, we see that we are in a position to apply \Cref{prop: 5.2} with $p = p_2$ and $q= p_3$, and thus, the inequality \begin{equation}\label{eq 12082025 b}  \Big| \overline{\mu}(a+b) - \overline{\mu}(a) - \overline{\mu}(b)| < 29 \frac{\varepsilon}{m},
\end{equation} holds for all $a,b \leq 0$ with $a+b \leq p_2 +p_3$.\smallskip

Condition $(1)$ and $(2)$ above with $i=1$ give 
$$\Big|\overline{\mu}(p) - \overline{\mu}(U_{p_1}(p)) - \overline{\mu}(U_{p_2+p_3}(p))\Big| < \gamma = \frac{\varepsilon}{m^2}  < 29 \frac{\varepsilon}{m},$$ for all $p \in \PP(\mathfrak{J})$ and 
$$\Big|\overline{\mu}(x+y)- \overline{\mu}(x) - \overline{\mu}(y) \Big| < \gamma= \frac{\varepsilon}{m^2}< 29 \frac{\varepsilon}{m},$$ for all $x,y \geq 0$ with $x+y \leq p_1$. The last two inequalities combined with \eqref{eq 12082025 b} allow us to apply \Cref{prop: 5.2} with $p= p_1$ and $q= p_2+p_3$ to deduce that 
$$ \Big| \overline{\mu}(a+b) - \overline{\mu}(a) - \overline{\mu}(b)\Big| < 29^2\varepsilon,$$
for all $a,b \geq 0$ with $a + b \leq \unit$. The proof concludes by the arbitrariness of the positive $\varepsilon$.
\end{proof}

\section{A Mackey-Gleason-Bunce-Wright theorem for \texorpdfstring{JBW$^*$-}{JBW*-}algebras}

The main goal of this research culminates now with a Mackey-Gleason-Bunce-Wright theorem (MGBW theorem in short) for signed and vector-valued finitely additive measures on the lattice of projections of a JBW$^*$-algebra without type $I_2$ summands, a result which has been pursued for decades.

\begin{theorem}\label{theo: M_G_JBW_star} Let $\mathfrak{J}$ be a JBW$^*$-algebra with no type $I_2$ direct summand. Then every bounded finitely additive measure $\mu: \PP(\mathfrak{J}) \rightarrow \mathbb{R}$ extends to a bounded linear functional on $\mathfrak{J}_{sa}$, equivalently, to a self-adjoint functional on $\mathfrak{J}$. Consequently, every bounded finitely additive measure $\mu: \PP(\mathfrak{J}) \rightarrow \mathbb{C}$ extends to a bounded linear functional on $\mathfrak{J}$.
\end{theorem}

\begin{proof} The conclusion for real-valued measures is a direct consequence of \Cref{theo: lin_type_I}, \Cref{theo: mu_lin_prop_inf} and \Cref{theo: mu_lin_type_II_1} and the structure theory of JBW$^*$-algebras recalled at subsection~\ref{subsec: structure}. If $\mu$ is a complex-valued measure we can apply the previous conclusion to $\Re\hbox{e}\mu$ and $\Im\hbox{m}\mu$. 
\end{proof}

As in the case of bounded measures on the lattice of projections of a von Neumann algebra without type $I_2$ direct summand (see \cite{BunceWright1994}), a vector-valued version of the MGBW theorem can be easily derived via Hahn-Banach theorem from \Cref{theo: M_G_JBW_star} above. The proof given in \cite[Lemma 1.1]{BunceWright1994} remains valid for JBW$^*$-algebras, and hence details are omitted. 

\begin{theorem}\label{t BWMG theorem for vector-valued measures} Let $\mathfrak{J}$ be a JBW$^*$-algebra with no type $I_2$ direct summand, and let $X$ be a Banach space. Then every bounded finitely additive measure $\mu: \PP(\mathfrak{J}) \rightarrow X$ admits a unique extension to a bounded linear operator from $\mathfrak{J}$ to $X$. 
\end{theorem}

We devote some final words concerning the hypothesis related to the type $I_2$ part. This concrete part is deeply connected with spin factors. As observed by Topping in \cite[Theorem 27]{Topping65}, and contrary to the case of von Neumann algebras, infinite dimensional type $I$ modular JW$^*$-factors exist. There are spin factors of arbitrary dimension $\geq 3$ (see \cite[Proposition 6.1.5]{HOS}). A JBW$^*$-algebra is a JBW$^*$-factor of type $I_2$ if and only if it is a spin factor \cite[Theorem 6.1.8]{HOS}. We recently included a brief review on the basic structure of spin factors in \cite[\S 1.4]{EscoPeVi2025}. The reader is referred to this source and the references therein, for all fine details. We shall simply recall that every spin factor has rank two (i.e. the cardinality of any family of mutually orthogonal projections is at most two), and the smallest spin factor is the three dimensional, which can be identified with the space $S_3(\mathbb{C})$ of all complex $3\times 3$ symmetric matrices. Kadison's original counterexample for $M_2 (\mathbb{C})$ also works in this case. Given a spin factor $\mathcal{S}$, we write $\PP_1 (\mathcal{S})$ for the set of all rank-one projections in $\mathcal{S}$. Set $p_1 = \left(\begin{matrix}
	1 & 0 \\
	0&0
\end{matrix} \right),$ and $p_2 = \left(\begin{matrix}
0 & 0 \\
0&1 
\end{matrix} \right)$ and define $\mu : \PP (S_3(\mathbb{C})) \to \mathbb{R}$ by $$\mu (\unit) =1 = \mu (p_1), \mu (0) =0 = \mu (p_2), \hbox{ and } \mu (p) =\frac{1}{2}, \forall p\in \PP_1 (S_3(\mathbb{C}))\backslash\{p_1,p_2\}.$$ It is easy to check that $\mu$ is a positive finitely additive measure which does not admit a linear extension to $S_3(\mathbb{C})$. Namely, every linear functional $\phi: S_3(\mathbb{C}) \to \mathbb{C}$ is represented in the form $\phi \left( \left(\begin{matrix}
a & b \\
b & c
\end{matrix} \right) \right) = \alpha a + \beta b + \gamma c$, for some $\alpha, \beta, \gamma\in \mathbb{C}$. If $\phi|_{\PP (S_3(\mathbb{C}))} = \mu$, the conditions $1 = \phi (p_1)$, $\phi (p_2) =0$, and $\phi \left( \left(\begin{matrix}
1/2 & 1/2 \\
1/2 & 1/2
\end{matrix} \right) \right) =\phi \left( \left(\begin{matrix}
1/3 & \sqrt{2}/3 \\
\sqrt{2}/3 & 2/3
\end{matrix} \right) \right) = \frac12$, lead to $a = 1$, $c=0,$ $b = 0$ and $1/3 = 1/2,$ which is impossible.\smallskip

On the other hand, it is known that every spin factor $\mathcal{S}$ admits a unital JBW$^*$-subalgebra Jordan $^*$-isomorphic to $S_3 (\mathbb{C})$. Under these conditions, let $\mu$ be the measure on $S_3(\mathbb{C})$ defined in the above paragraph. The mapping $\tilde{\mu} : \PP (\mathcal{S}) \to \mathbb{R}$ by $\tilde{\mu} (p) = \mu (p),$ if  $p\in S_3(\mathbb{C}),$ and $\tilde{\mu} (p) = \frac12$ otherwise. Having in mind that $S_3(\mathbb{C})$ is a unital JBW$^*$-subalgebra of $\mathcal{S}$ and the latter has rank $2$, it is not hard to check that $\tilde{\mu}$ is a positive finitely additive measure. Clearly $\tilde{\mu}$ does not admit an extension to a bounded functional on $\mathcal{S}$, since its restriction to $S_3(\mathbb{C})$ is precisely $\mu$.\smallskip

We finally deal with type $I_2$ JBW$^*$-algebras. Given a Radon measure $\nu$ on a locally compact Hausdorff space $\Omega$, a mapping $f$ from $\Omega$ to a Banach space $X$ is called \emph{strongly $\nu$-measurable} if for each positive $\epsilon$ and each compact subset $K\subseteq \Omega$, there exists a compact subset $K_{\varepsilon} \subseteq K$ such that $\nu (K\backslash K_{\varepsilon}) < \varepsilon$ and $f$ is continuous on $K_{\varepsilon}$. We denote by $L_{\infty} (\Omega, \nu, X)$ the Banach space
of all equivalence classes of bounded strongly $\nu$-measurable functions from $\Omega$ to $X$ under the equivalence relation of equality locally $\nu$-almost everywhere. Let $\mathfrak{J}$ be a JB$^*$-algebra. By \cite[Theorems 1 and 2]{Stacey82}, there exists a locally compact Hausdorff space $\Omega$, a Radon measure $\nu$ on $\Omega$, and a spin factor $\mathcal{S}$, such that $L_{\infty} (\Omega, \nu, \mathcal{S})$ is a summand in $\mathfrak{J}$. Clearly $\mathcal{S}$ embeds as the subalgebra of all constant functions inside $L_{\infty} (\Omega, \nu, \mathcal{S})$. Arguing as in the proof of \cite[Corollary in page 124]{Stacey82}, we can find a Jordan $^*$-homomorphism $\pi_1: L_{\infty} (\Omega, \nu, \mathcal{S}) \to \mathcal{S}$ such that $\pi_1 (a) = a$ for all $a\in \mathcal{S}$. Having in mind that $L_{\infty} (\Omega, \nu, \mathcal{S}) $ is a summand in $\mathfrak{J}$, we can define Jordan $^*$-homomorphisms $\iota : \mathcal{S}\hookrightarrow  \mathfrak{J}$, and $\pi: \mathfrak{J} \to \mathcal{S}$ such that  $\iota$ is an isometric embedding, $\pi (\iota (a)) = a,$ for all $a\in \mathcal{S}$. Let $\tilde{\mu}: \PP(\mathcal{S})\to \mathbb{R}$ be the measure built in the previous paragraph. Define $\tilde{\tilde{\mu}}: \PP(\mathfrak{J})\to \mathbb{R}$ the measure given by $\tilde{\tilde{\mu}} (p) = {\tilde{\mu}} (\pi (p))$. It is easy to check that $\tilde{\tilde{\mu}}$ is positive and finitely additive, since ${\tilde{\mu}}$ is and $\pi$ is a Jordan $^*$-homomorphism. If there exists a bounded linear functional $\phi :\mathfrak{J}\to \mathbb{C}$ whose restriction to $\PP (\mathfrak{J})$ is $\tilde{\tilde{\mu}}$, the functional $\tilde{\phi} =\phi\iota : \mathcal{S}\to \mathbb{C}$ is linear and continuous and satisfies $$\tilde{\phi}(p) = \phi (\iota (p)) = \tilde{\tilde{\mu}} (\iota (p)) ={\tilde{\mu}}  (\pi(\iota(p))) = {\tilde{\mu}} (p), \ \forall p\in \PP(\mathcal{S}),$$ which contradicts that ${\tilde{\mu}}$ does not admit a linear extension to $\mathcal{S}$. \medskip\medskip

\textbf{Acknowledgements}\medskip

\noindent Authors supported by MICIU/AEI/10.13039/501100011033 and ERDF/EU grant PID2021-122126NB-C31, by ``Maria de Maeztu'' Excellence Unit IMAG, reference CEX2020-001105-M funded by MICIU/AEI/10.13039/501100011033, and by Junta de Andaluc{\'i}a grants FQM185 and FQM375. First author supported by grant FPU21/00617 at University of Granada founded by Ministerio de Universidades (Spain). Second author partially supported by MOST China project number G2023125007L.

\smallskip\smallskip

\noindent\textbf{Data Availability} Statement Data sharing is not applicable to this article as no datasets were generated or analyzed during the preparation of the paper.\smallskip\smallskip

\noindent\textbf{Declarations} 
\smallskip\smallskip

\noindent\textbf{Conflict of interest} The authors declare that he has no conflict of interest.

\end{document}